%% file: Hitchin_metric.tex
\documentclass[a4,10pt]{article}
\usepackage{amsmath,amscd,amssymb}
\usepackage[margin=.8in]{geometry}

\input{notation}

\input{new_theorem}

\begin{document}

\title{Asymptotic behaviour of the Hitchin metric\\
on the moduli space of Higgs bundles}

\author{Takuro Mochizuki\thanks{Research Institute for Mathematical Sciences, Kyoto University, Kyoto 606-8512, Japan, takuro@kurims.kyoto-u.ac.jp}}
\date{}
\maketitle

\begin{abstract}
The moduli space of stable Higgs bundles of degree $0$
is equipped with the hyperk\"ahler metric,
called the Hitchin metric.
On the locus where the spectral curves are smooth,
there is the  hyperk\"ahler metric called the semi-flat metric,
associated with the algebraic integrable systems with the Hitchin section. 
We prove the exponentially rapid decay
of the difference between the Hitchin metric
and the semi-flat metric along the ray
$(E,t\theta)$ as $t\to\infty$.

\vspace{.1in}
\noindent
MSC: 53C07, 58E15, 14D21, 81T13.
\\
Keywords:
harmonic bundle,
moduli space of Higgs bundles,
Hitchin metric,
semi-flat metric
\end{abstract}

\section{Introduction}

\subsection{The Hitchin metric and the semi-flat metric of
the moduli space of Higgs bundles}

\subsubsection{Hitchin metric}

Let $X$ be a compact Riemann surface
with a conformal metric $g_X$.
In this paper,
a Higgs bundle on $X$ is
a holomorphic vector bundle $E$
with a holomorphic section $\theta$
of $\End(E)\otimes K_X$,
where $K_X$ denotes the canonical bundle of $X$.
The study of the moduli space of Higgs bundles
was initiated by Hitchin in \cite{Hitchin-self-duality}.
Let $\nbigm_H$ denote the moduli space of
stable Higgs bundles of degree $0$ and rank $n$ on $X$.
It is naturally equipped with
the integrable almost complex structure $I$.
Moreover, $\nbigm_H$ is equipped with
the naturally defined hyperk\"ahler structure,
i.e.,
the Riemannian metric $g_{H,\real}$
and the three integrable almost complex structures $I$, $J$ and $K$
such that
(i) the associated $2$-forms
$\omega_{T}(\cdot,\cdot)=g_{H,\real}(\cdot,T\cdot)$
$(T=I,J,K)$ are closed,
(ii) $IJK=-1$.
Let $g_H$ denote the associated K\"ahler metric
of $(\nbigm_H,I)$,
i.e.,
$g_{H}=g_{H,\real}+\sqrt{-1}\omega_{I}$.
The Riemannian metric $g_{H,\real}$ and
the K\"ahler metric $g_{H}$
are called the Hitchin metric.

By the theorem of Hitchin and Simpson \cite{Hitchin-self-duality, s1},
because $(E,\theta)\in\nbigm_H$ is stable of degree $0$,
it carries a harmonic metric $h$,
i.e.,
$h$ is a Hermitian metric of $E$
satisfying
\[
 R(h)+[\theta,\theta^{\dagger}_h]=0.
\]
Here, $R(h)$ denotes the curvature of the Chern connection of $(E,h)$,
and $\theta^{\dagger}_h$ denotes the adjoint of $\theta$
with respect to $h$.
The induced Hermitian metric on $\End(E)$ is also denoted by $h$.
The tangent space $T_{(E,\theta)}\nbigm$
is identified with
the space of harmonic $1$-forms of
the harmonic bundle $(\End(E),\ad\theta,h)$.
The Hermitian metric $(g_H)_{|(E,\theta)}$ is identified with
the $L^2$-metric on the space of harmonic $1$-forms.
(For example, see \S\ref{subsection;23.5.16.40}.)
Namely, for harmonic $1$-forms $\tau$ of
$(\End(E),\ad\theta,h)$,
let $[\tau]\in T_{(E,\theta)}\nbigm_H$
denote the corresponding element,
and then we have
\[
 g_{H|(E,\theta)}([\tau],[\tau])
 =2\sqrt{-1}\int_X\Bigl(
 h(\tau^{1,0},\tau^{1,0})
 -h(\tau^{0,1},\tau^{0,1})
 \Bigr)
 =2\sqrt{-1} \int_X
 \tr\Bigl(
 \tau^{1,0}\wedge (\tau^{1,0})^{\dagger}_h
-\tau^{0,1}\wedge(\tau^{0,1})^{\dagger}_h
 \Bigr).
\]
It is attractive and challenging to study
the asymptotic behaviour of $g_H$ around infinity
of the non-compact space $\nbigm_H$.

\subsubsection{Semi-flat metric}

We set
$\nbiga_H=\bigoplus_{j=1}^n H^0(X,K_X^j)$.
Let $\pi_0:T^{\ast}X\to X$ denote the projection.
Let $\eta$ denote the tautological section of
$\pi_0^{\ast}(K_X)$.
For any $a=(a_j)\in\nbiga_H$,
we obtain the section
$\xi_a=\eta^n+\sum_{j=1}^{n} (-1)^ja_j\eta^{n-j}$
of $\pi_0^{\ast}(K_X^{\otimes n})$.
Let $\Sigma_a$ denote the $0$-scheme of the section $\xi_a$.
Such a curve is called a spectral curve.
Let $\nbiga_H'\subset \nbiga_H$ denote
the Zariski open subset of
$a\in\nbiga_H$
such that $\Sigma_a$ are smooth.
It is known that $\nbiga_H'$ is non-empty.

For any Higgs bundle $(E,\theta)$ of rank $n$,
we obtain
$\Phi_H(E,\theta)=(a_j(E,\theta)\,|\,j=1,\ldots,n)\in\nbiga_H$
determined by
\[
\det(\eta\id_E-\theta)
=\eta^n+\sum_{j=1}^n(-1)^ja_j(E,\theta)\eta^{n-j}.
\]
It induces the morphism $\Phi_H:\nbigm_H\to \nbiga_H$
called the Hitchin fibration.
There exists the Hitchin section of $\Phi_H$
depending on a root $K_X^{1/2}$ of $K_X$.
For any $(E,\theta)\in\nbigm_H$,
the curve $\Sigma_{\theta}=\Sigma_{\Phi_H(E,\theta)}$
is called the spectral curve of $(E,\theta)$.

Let $\nbigm_H'$ denote the Zariski open subset of
$(E,\theta)\in\nbigm_H$
such that the spectral curves
$\Sigma_{\theta}$ are smooth,
i.e.,  $\nbigm_H'=\Phi_H^{-1}(\nbiga_H')$.
We obtain the induced smooth fibration
$\Phi_H':\nbigm_H'\to\nbiga_H'$
which is an algebraic integrable system
in the sense of \cite{Freed}.
Each fiber $(\Phi_H')^{-1}(a)$ $(a\in\nbiga_H')$
is naturally polarized.
According to \cite[Theorem 3.8]{Freed},
$\nbigm_H'$ carries a canonical hyperk\"ahler structure,
i.e.,
the Riemannian metric $g_{\semiflat,\real}$
with the integrable almost complex structures
$I,J',K'$ satisfying the axioms of hyperk\"ahler structure.
We set $g_{\semiflat}=g_{\semiflat,\real}+\sqrt{-1}\omega_I$.
The Riemannian metric $g_{\semiflat,\real}$
and the K\"ahler metric $g_{\semiflat}$ are called
the semi-flat metric.

Let us recall how $g_{\semiflat}$ is constructed.
For $a=(a_i)\in\nbiga_{H,n}'$,
let $U_{a}$ be a simply connected neighbourhood
of $a$ in $\nbiga_{H,n}'$.
We may regard $\Phi'_{H|U_a}:\Phi_H^{-1}(U_{a})\to U_a$
as a principal torus bundle.
By using a Hitchin section,
we define the integrable connection
of the principal torus bundle,
which induces the decomposition of
the tangent space $T_{(E,\theta)}\nbigm_H'$
into the horizontal direction and the vertical direction:
\[
 T_{(E,\theta)}\nbigm_H'
 =(T_{(E,\theta)}\nbigm_H')^{\hor}
 \oplus
 (T_{(E,\theta)}\nbigm_H')^{\ver}.
\]
It is independent of the choice of a Hitchin section.
Because the fiber $\Phi_H^{-1}(\Phi_H(E,\theta))$
is identified with a connected component of
the Picard variety of $\Sigma_{\theta}$,
$(T_{(E,\theta)}\nbigm_H')^{\ver}$ is naturally identified with
$H^1(\Sigma_{\theta},\nbigo_{\Sigma_{\theta}})$.
The moduli space $\nbigm_{H}'$ is equipped
with a complex symplectic structure.
Because $\Phi_H^{-1}(\Phi_H(E,\theta))$ is Lagrangian,
$(T_{(E,\theta)}\nbigm_H')^{\ver}$
and
$(T_{(E,\theta)}\nbigm_H')^{\hor}$
are mutually dual.
Hence, $(T_{(E,\theta)}\nbigm_H')^{\hor}$
is naturally identified with
$H^0(\Sigma_{\theta},K_{\Sigma_{\theta}})$.
The space 
$H^1(\Sigma_{\theta},\nbigo_{\Sigma_{\theta}})
\oplus
H^0(\Sigma_{\theta},K_{\Sigma_{\theta}})$
is identified with the space of
harmonic $1$-forms on $\Sigma_{\theta}$.
The Hermitian metric
$(g_{\semiflat})_{(E,\theta)}$
is identified with
the $L^2$-metric
on the space of harmonic $1$-forms on $\Sigma_{\theta}$.
(See Proposition \ref{prop;23.5.16.1}.)
Namely,
for harmonic $1$-forms $\tau$ on $\Sigma_{\theta}$,
let $[\tau]\in
H^1(\Sigma_{\theta},\nbigo_{\Sigma_{\theta}})
\oplus
H^0(\Sigma_{\theta},K_{\Sigma_{\theta}})$
denote the corresponding element,
and then we have
\[
 g_{\semiflat|(E,\theta)}([\tau],[\tau])
=
 2\sqrt{-1}\int_{\Sigma_{\theta}}
 \Bigl(
 \tau^{1,0}\wedge \overline{\tau^{1,0}}
-\tau^{0,1}\wedge \overline{\tau^{0,1}}
 \Bigr).
\]

\subsubsection{A weak version of a conjecture of Gaiotto-Moore-Neitzke}

A conjecture of
Gaiotto-Moore-Neitzke
\cite{Gaiotto-Moore-Neitzke, Gaiotto-Moore-Neitzke2}
implies the following estimate for any $(E,\theta)\in\nbigm_H'$
with respect to $(g_{\semiflat})_{(E,t\theta)}$
as $t\to\infty$:
\begin{equation}
\label{eq;23.5.8.10}
(g_H-g_{\semiflat})_{(E,t\theta)}=O(e^{-\beta t}).
\end{equation}
Here, $\beta$ is a positive constant depending on $(E,\theta)$.
It is our purpose in this paper
to prove this weak version of GMN-conjecture.
Note that their original conjecture is much more precise.

A previous work was done by Fredrickson in \cite{Fredrickson2}
on the basis of
her another work \cite{Fredrickson-Generic-Ends},
the works of Mazzeo-Swoboda-Weiss-Witt
\cite{Mazzeo-Swoboda-Weiss-Witt, MSWW2}
and interesting ideas of Dumas-Neitzke
\cite{Dumas-Neitzke}.
More precisely,
for each $(E,\theta)\in\nbigm_H'$,
let $\Crit(\theta)\subset\Sigma_{\theta}$
denote the set of the critical points
of the projection $\pi:\Sigma_{\theta}\to X$.
For each $Q\in \Crit(\theta)$,
let $r_Q$ denote the ramification index of $\pi$ at $Q$.
We obtain the Zariski open subset 
$\nbigm_H''$ of $(E,\theta)\in\nbigm_H'$
such that $r_Q\leq 2$ for any $Q\in\Crit(\theta)$.
We have $\nbigm'_H=\nbigm_H''$ if $n=2$
but $\nbigm_H'\neq\nbigm_H''$ if $n\geq 3$, in general.
The conjectural estimate (\ref{eq;23.5.8.10}) has been studied
for any $(E,\theta)\in\nbigm_H''$ in \cite{Fredrickson2}.
We should also note that
a similar estimate
has been studied by Fredrickson-Mazzeo-Swoboda-Weiss
\cite{FMSW}
in the case of parabolic Higgs bundles of rank $2$.
Chen-Li \cite{Chen-Li} studied
the rank $2$ but more singular case
by using the methods in \cite{Fredrickson2,FMSW}
and this paper.
In \cite{Holdt}, Holdt studied the original approach of
\cite{Gaiotto-Moore-Neitzke, Gaiotto-Moore-Neitzke2}
in the case that a Higgs bundle has a non-trivial parabolic structure.
Note that because it is based on Fock-Goncharov coordinates,
it is not known whether the approach can be applied
to the non-singular case.

We shall prove
the conjectural estimate (\ref{eq;23.5.8.10})
for any $(E,\theta)\in\nbigm_H'$
by using a different method.
More precisely,
let $s$ be the automorphism of
the complex vector bundle $T\nbigm_H'$
determined by $g_{H}=g_{\semiflat}\cdot s$.
\begin{thm}[Theorem
\ref{thm;23.5.11.20}]
\label{thm;23.5.12.40}
For any compact subset $K\subset\nbigm_H'$,
there exist positive constants $B_i$ $(i=1,2)$
such that
\[
|(s-\id)_{(E,t\theta)}|_{g_{\semiflat}}\leq B_1\exp(-B_2t)
\]
for any $(E,\theta)\in K$. 
\end{thm}

\begin{rem}
\label{rem;24.10.23.20}
Regrettably,
to the best of the author's understanding,
some more explanations or modifications
are necessary for
the arguments in {\rm\cite{Fredrickson2,FMSW,MSWW2}}
\footnote{The author thanks Laura Fredrickson
for sincerely replying to the author's questions.
The author also thanks Rafe Mazzeo
for answering to the inquiry regarding this remark.}.
\hfill\qed
\end{rem}

\begin{rem}
In {\rm\cite{Mochizuki-rank-two}},
we study the constant $B_2$ more closely for the rank $2$ case.
In {\rm\cite{Mochizuki-non-zero-degree}},
we study a generalization to the non-zero degree case. 
\hfill\qed
\end{rem}

\subsection{Large-scale solutions of the Hitchin equation}

\subsubsection{Large-scale solutions}

Let $(E,\theta)\in\nbigm_H'$.
We fix a flat metric $h_{\det(E)}$ of $\det(E)$.
For any $t>0$,
we obtain harmonic metrics
$h_t$ of $(E,\theta)$ such that $\det(h_t)=h_{\det(E)}$.
For large $t$,
they are also called large-scale solutions
of the Hitchin equation.
To study the asymptotic behaviour of $g_{H}$ along the ray $(E,t\theta)$,
it is basic to study the asymptotic behaviour of
large-scale solutions $h_t$ as $t\to\infty$.

\subsubsection{Limiting configuration}

Let $\pi:\Sigma_{\theta}\to X$ denote the projection.
According to Hitchin \cite{Hitchin-self-duality}
and Beauville-Narasimhan-Ramanan \cite{Beauville-Narasimhan-Ramanan},
there exists a line bundle $L$ on $\Sigma_{\theta}$
with an isomorphism
\[
 (\pi_{\ast}(L),\pi_{\ast}(\eta))\simeq (E,\theta),
\]
where $\pi_{\ast}(\eta)$ is the Higgs field
induced by the tautological $1$-form $\eta$ on $T^{\ast}X$
as the push-forward of
$\eta:L\to L\otimes \pi^{\ast}(K_X)$.
We set $D(\theta)=\pi(\Crit(\theta))$
and $X^{\circ}=X\setminus D(\theta)$.
For any subset $S\subset X$,
we set $\Sigma_{\theta|K}:=\pi^{-1}(S)$.
There exists a Hermitian metric $h_L$
of $L_{|\Sigma_{\theta|X^{\circ}}}$
such that the following holds.
\begin{itemize}
 \item $h_L$ is flat, i.e.,
       $R(h_L)=0$.
 \item For each $Q\in \pi^{-1}(D(\theta))$,
       there exist a local holomorphic coordinate $\zeta_Q$
       and a local frame $v_Q$ around $Q$
       such that $\zeta_Q(Q)=0$
       and $h_L(v_Q,v_Q)=|\zeta_Q|^{-r_Q+1}$.
\end{itemize}
By taking the push-forward of $h_L$ with respect to $\pi$,
we obtain a Hermitian metric $h_{\infty}$ of $E_{|X^{\circ}}$.
It is a decoupled harmonic metric in the sense
\[
R(h_{\infty})=[\theta,\theta^{\dagger}_{h_{\infty}}]=0.
\]
We also note that $\det(h_{\infty})$ induces a $C^{\infty}$-metric of
$\det(E)$ on $X$ whose curvature is $0$,
and hence we may assume that $\det(h_{\infty})=h_{\det(E)}$
by adjusting a positive constant.
It is called the limiting configuration.
It was introduced by Mazzeo-Swoboda-Weiss-Witt
\cite{Mazzeo-Swoboda-Weiss-Witt} in the rank $2$ case,
and by Fredrickson \cite{Fredrickson-Generic-Ends}
in the higher rank case.

\subsubsection{Convergence}
\label{subsection;24.10.23.1}

In \cite{Mazzeo-Swoboda-Weiss-Witt},
Mazzeo-Swoboda-Weiss-Witt
proved the exponentially rapid convergence of
the sequence $h_t$ to $h_{\infty}$ 
as $t\to\infty$
locally on $X^{\circ}$.
It is generalized to the higher rank case
by Fredrickson in \cite{Fredrickson-Generic-Ends}
under some mild assumption
on the ramification of the projection $\Sigma_{\theta}\to X$
(see \cite[Remark 1.2]{Mochizuki-Szabo}),
and by Szabo and the author
in \cite{Mochizuki-Szabo} more generally.

We recall a very rough strategy
in \cite{Mazzeo-Swoboda-Weiss-Witt}.
For simplicity, we assume $\tr(\theta)=0$,
$\det(E)=\nbigo_X$ and $h_{\det(E)}=1$.
For any point $P\in D(\theta)$,
let $(X_P,z_P)$ be a holomorphic coordinate neighbourhood around $P$
such that
(i) $z_P(P)=0$,
(ii) $X_P\simeq \{|z_P|<2\}$,
(iii) $X_P\cap D(\theta)=\{P\}$.
We set $(E_P,\theta_P):=(E,\theta)_{|X_P}$.
For any $0<r\leq 2$,
we set $X_P(r)=\{|z_P|<r\}$.
\begin{description}
 \item[Step 1] 
       We construct harmonic metrics
       $h_{P,t}$ $(t>0)$ of
	    $(E_P,t\theta_P)$
	    such that
	    (i) $\det(h_{P,t})=h_{\det(E)|X_P}$,
	    (ii)
       $h_{P,t}-h_{\infty}$
       and their derivatives are
       $O(e^{-\epsilon t})$
       on $X(1)\setminus X(1/2)$
       for some $\epsilon>0$.
 \item[Step 2] We patch $h_{\infty}$
       and $h_{P,t}$ $(P\in D(\theta))$
       by using a partition of unity
       to construct a family of Hermitian metrics
       $\htilde_t$ of $E$
	    such that 
	    (i) $\det(\htilde_t)=h_{\det(E)}$,
	    (ii) $\htilde_t-h_{\infty}=O(e^{-\epsilon't})$
	    on $X\setminus \bigcup_{P\in D(\theta)}X_P(1/2)$,
	    (iii) $\htilde_t-h_{P,t}=O(e^{-\epsilon't})$ on $X_P$,
	    (iv) $R(\htilde_t)+t^2[\theta,\theta^{\dagger}_{\htilde_t}]
       =O(e^{-\epsilon't})$
	    holds on $X$ for some $\epsilon'>0$.
	    In particular,
       we have the exponentially rapid convergence of
       the sequence $\htilde_t$
       to $h_{\infty}$
       locally on $X^{\circ}$.
 \item[Step 3] We show
       $\htilde_t-h_t=O(e^{-\delta t})$
       for some 
       $\delta>0$ on $X$.
\end{description}

In \cite{Mazzeo-Swoboda-Weiss-Witt},
the so called fiducial solutions play important roles
for the construction of $h_{P,t}$ in Step 1.
In the rank $2$ case,
the Higgs bundle $(E_P,\theta_P)$ is essentially uniquely determined,
and it is highly symmetric.
Namely, we may normalize the coordinate $z_P$ such that
$\det(\theta_P)=-z_P(dz_P)^2$,
and there exists a frame
$v_1,v_2$ of $E_P$
such that
\[
 \theta_P(v_1)=v_2\,dz_P,
 \quad
 \theta_P(v_2)=v_1\,z_P\,dz_P,
 \quad
 v_1\wedge v_2=1.
\]
There exists a harmonic metric $h_{P,t}$ of $(E_P,\theta_P)$
such that
(i) $h_{P,t}(v_i,v_i)$ depends only on $|z_P|$,
(ii) $h_{P,t}(v_i,v_j)=0$ $(i\neq j)$,
(iii) $\det(h_{P,t})=1$.
We can find such a harmonic metric
by solving a non-linear ordinary differential equation,
which is a special case of Painlev\'{e} III equation,
or by rescaling a wild harmonic bundle on $\cnum$.
The solutions are called fiducial solutions.
We can adjust $h_{P,t}$ so that the conditions in Step 1
are satisfied.
This method can be also useful in the higher rank case
if $r_Q\leq 2$ for any $Q\in \Crit(\theta)$,
and it was efficiently used in
\cite{Fredrickson2, Fredrickson-Generic-Ends}.
However, it is not clear whether the method is so useful
for the study of large-scale solutions
in the higher rank case
without any assumption on the ramification.

\subsubsection{Locally defined symmetric pairings and a new method}

In \cite{Mochizuki-Szabo},
we introduced another method
for the construction of $h_{P,t}$
for $(E,t\theta)$ as in Step 1,
which is available in a more general case.
There exists a non-degenerate symmetric pairing
$C_P$ of $E_P$ such that the following holds
(see \cite{Mochizuki-Szabo}):
\begin{itemize}
 \item $C_P$ is a non-degenerate symmetric pairing of $(E_P,\theta_P)$,
       i.e.,
       $\theta_P$ is self-adjoint with respect to $C_P$.
 \item $C_{P|X_P\setminus \{P\}}$ is compatible with
       $h_{\infty|X_P\setminus\{P\}}$
       in the sense
       that the induced morphism
       $\Psi_{C_P}:E_P\to E_P^{\lor}$
       is isometric with respect to
       $h_{\infty}$ and the induced metric $h_{\infty}^{\lor}$
       on $X_P\setminus\{P\}$.
\end{itemize}
By the theorem of Donaldson on the Dirichlet problem
for harmonic metrics \cite{Donaldson-boundary-value},
there exists a harmonic metric $h_{P,t}$
such that
$h_{P,t|\del X_{P}}=h_{\infty|\del X_P}$.
It is easy to observe that
$C_P$ and $h_{P,t}$ are compatible
(see \cite{Li-Mochizuki3}).
By an estimate in \cite{Mochizuki-Szabo}
on the basis of a variant of Simpson's main estimate
in \cite{Decouple},
which goes back to \cite{s2},
the conditions in Step 1 are automatically satisfied.

We remark that the argument for Step 3 in \cite{Mochizuki-Szabo}
is also different from those in
\cite{Fredrickson-Generic-Ends, Mazzeo-Swoboda-Weiss-Witt}.
In the latter,
after the detailed analysis of the fiducial solutions,
they closely study the family of Green operators 
associated with $(E,t\theta,\htilde_t)$,
and eventually they proved for sufficiently large $t$,
that there exists $h_t$ in an exponentially small neighbourhood
of $\htilde_t$ in the space of Hermitian metrics
by using the inverse function theorem.
The argument in \cite{Mochizuki-Szabo}
is based on the maximum principle of subharmonic functions.

The study
\cite{Fredrickson-Generic-Ends, Mazzeo-Swoboda-Weiss-Witt}
were motivated by the conjectural estimate (\ref{eq;23.5.8.10}),
and actually applied in \cite{Fredrickson2, MSWW2}.
We would like to check whether
our result and method in \cite{Mochizuki-Szabo}
are useful enough
for the study of (\ref{eq;23.5.8.10}),
that is our original purpose in this paper.

We note that our route to (\ref{eq;23.5.8.10})
is apparently different from those in \cite{Fredrickson2, MSWW2},
where they first show that $g_{\semiflat}$
is identified with the natural $L^2$-metric $g_{L^2}(\nbigm'_{\infty})$
on the moduli space $\nbigm_{\infty}'$
of Higgs bundles with the limiting configurations,
and they show that
$g_H-g_{L^2}(\nbigm'_{\infty})=O(e^{-\beta t})$
along the ray $(E,t\theta)$ as $t\to\infty$.
(See \cite{Fredrickson2} for more precise explanations.)
We shall compare $g_{H}$ and $g_{\semiflat}$ more directly.
An outline is explained in the next subsection.

\subsection{Outline}

\subsubsection{The horizontal direction and the vertical direction}

Let $(E,\theta)\in\nbigm_H'$.
Let $\Def(E,\theta)$ be the complex of sheaves
\[
\begin{CD}
 \End(E)
 @>{\ad\theta}>>
 \End(E)\otimes K_X.
\end{CD}
\]
Here, the first term sits in the degree $0$.
The hypercohomology group of the complex $\Def(E,\theta)$
is denoted as $H^{\ast}(X,\Def(E,\theta))$.
It is well known that
there exists the following natural isomorphism
(see \S\ref{subsection;24.10.22.20} for an explanation):
\[
T_{(E,\theta)}\nbigm_H'
\simeq
H^1(X,\Def(E,\theta)).
\]
We shall identify them by the isomorphism.
There exist natural isomorphisms
(see Proposition \ref{prop;23.3.23.10}):
\[
 \nbigh^0(\Def(E,\theta))
 \simeq
 \pi_{\ast}\nbigo_{\Sigma_{\theta}},
 \quad\quad
 \nbigh^1(\Def(E,\theta))
 \simeq
 \pi_{\ast}K_{\Sigma_{\theta}}.
\]
We obtain the following exact sequence:
\[
\begin{CD}
 0 @>>>
 H^1(\Sigma_{\theta},\nbigo_{\Sigma_{\theta}})
 @>{\iota^{\ver}}>>
 H^1(X,\Def(E,\theta))
 @>{\Upsilon_1}>>
 H^0(\Sigma_{\theta},K_{\Sigma_{\theta}})
 @>>> 0.
\end{CD}
\]
It is easy to observe that
$\iota^{\ver}$ induces the isomorphism of
$H^1(\Sigma_{\theta},\nbigo_{\Sigma_{\theta}})$
and the tangent space of
$\Phi_H^{-1}(\Phi_H(E,\theta))\subset\nbigm'_H$,
i.e.,
$(T_{(E,\theta)}\nbigm'_H)^{\ver}$.
There uniquely exists the splitting
$\iota^{\hor}:H^0(\Sigma_{\theta},K_{\Sigma_{\theta}})
\to H^1(X,\Def(E,\theta))$
such that
the image $\iota^{\hor}$ equals
$(T_{(E,\theta)}\nbigm'_H)^{\hor}$.

\subsubsection{$C^{\infty}$-representatives of
the horizontal infinitesimal deformations}

Let $A^{\bullet}(\Def(E,\theta))$ denote the complex
obtained as the space of global sections of
the Dolbeault resolution of $\Def(E,\theta)$.
There exists the line bundle $L$ on $\Sigma_{\theta}$
with an isomorphism
$(\pi_{\ast}(L),\pi_{\ast}(\eta))\simeq (E,\theta)$.
For any $\nu\in H^0(\Sigma_{\theta},K_{\Sigma_{\theta}})$,
the multiplication of $\nu$ on $L$
induces $F_{\nu}\in H^0(X,\End(E)\otimes K_X(D))$.
We would like to modify $F_{\nu}$ to
a $1$-cocycle $\tau$
in $A^{\bullet}(\Def(E,\theta))$
such that
the cohomology class $[\tau]$
equals $\iota^{\hor}(\nu)$.
We remark that
any $1$-cocycle $\tau$ in $A^1(\Def(E,\theta))$
naturally induces
a holomorphic section $\rho(\tau)$
of $\pi_{\ast}(K_{\Sigma_{\theta}})$
such that $\Upsilon_1([\tau])=\rho(\tau)$.
(See \S\ref{subsection;23.5.12.2}.)

Let $\End^{\sym}(E_P;C_P)$ (resp. $\End^{\asym}(E_P;C_P)$)
denote the sheaf of holomorphic endomorphisms of $E_P$
which are symmetric (resp. anti-symmetric)
with respect to $C_P$.
We obtain the following subcomplex
$\Def(E_P,\theta_P;C_P)$ of $\Def(E_P,\theta_P)$:
\[
\begin{CD}
\End^{\asym}(E_P;C_P)
@>{\ad\theta_P}>>
\End^{\sym}(E_P;C_P)\otimes K_{X_P}.
\end{CD}
\]
Let $A^{\bullet}(\Def(E_P,\theta_P;C_P))$ denote the complex
obtained as the space of the global sections of
the Dolbeault resolution of
$\Def(E_P,\theta_P;C_P)$.
Then, we can obtain the following sufficient condition
for a representative of $\iota^{\hor}(\nu)$
in $A^1(\Def(E,\theta))$.
\begin{prop}[Corollary
 \ref{cor;23.5.12.1}]
\label{prop;23.5.13.1}
A $1$-cocycle $\tau\in A^1(\Def(E,\theta))$
satisfies $[\tau]=\iota^{\hor}(\nu)$
if the following holds. 
\begin{itemize}
 \item For each $P\in D(\theta)$,
       there exists a relatively compact open neighbourhood $X_P'$
       of $P$ in $X_P$
       such that
       $\tau_{|X\setminus\bigcup X_P'}
       =F_{\nu|X\setminus\bigcup X_P'}$.       
 \item $\tau_{|X_P}\in A^1(\Def(E_P,\theta_P;C_P))$.
\end{itemize}
\end{prop}

\subsubsection{Harmonic representatives
for horizontal infinitesimal deformations}

Let $h_t$ be the harmonic metric of $(E,t\theta)$
such that $\det(h_t)=h_{\det(E)}$.
Let $g_X$ be a conformal metric of $X$.
Let $\delbar_E$ denote the Cauchy operator of $E$.
The induced Cauchy operator of $\End(E)$ is also denoted by $\delbar_E$.
Let $\iota^{\hor}_t:H^0(\Sigma_{\theta},K_{\Sigma_{\theta}})
\to H^1(X,\Def(E,t\theta))$
denote the splitting morphism for $(E,t\theta)$,
where we identify
$\Sigma_{t\theta}$ and $\Sigma_{\theta}$ in a natural way.
To study the harmonic representative
$\ttH(\nu,t)$ of $\iota^{\hor}_t(\nu)$,
we construct a $1$-cocycle
$\ttH'(\nu,t)\in A^1(\Def(E,t\theta))$
such that the following holds.
\begin{itemize}
 \item[(i)] $\rho(\ttH'(\nu,t))=\pi_{\ast}(\nu)$.
 \item[(ii)] $\ttH'(\nu,t)$ satisfies
	     the conditions in
	     Proposition \ref{prop;23.5.13.1}.
	     In particular,
	     $[\ttH'(\nu,t)]=\iota^{\hor}_t(\nu)$.
 \item[(iii)]
	    $(\delbar_E+\ad t\theta)^{\ast}_{h_t,g_X}\ttH'(\nu,t)
	     =O(e^{-\beta t})\|\nu\|$ for some $\beta>0$,
	     where
	     $\|\nu\|$ denotes the $L^2$-norm of $\nu$
	     on $\Sigma_{\theta}$.
 \item[(iv)]
	    $\tr\bigl(
	    (\delbar_E+\ad t\theta)^{\ast}_{h_t,g_X}\ttH'(\nu,t)
	    \bigr)=0$.
\end{itemize}

Let $P\in D(\theta)$.
There exists a holomorphic function $\alpha_P$
on $\Sigma_{\theta|X_P}$
such that $d\alpha_P=\nu_{|\Sigma_{\theta|X_P}}$
and that $\alpha_P(Q)=0$ for any $Q\in \pi^{-1}(P)$.
The multiplication of $\alpha_P$ on $L_{|\Sigma_{\theta|X_P}}$
induces an endomorphism
$F_{\alpha_P}$ of $E_P$
such that $[\theta,F_{\alpha_P}]=0$.
Let $h_{P,t}$ be the harmonic metric of $(E_P,t\theta_P)$
such that $h_{P,t|\del X_P}=h_{\infty|\del X_P}$.
Let $\nabla_{h_{P,t}}=\delbar_E+\del_{E_P,h_{P,t}}$
denote the Chern connection of $(E,h_{P,t})$.
Let $\theta^{\dagger}_{h_{P,t}}$ denote the adjoint of $\theta$
with respect to $h_{P,t}$.
Then, it is easy to see that
\[
 \ttH_P(\nu,t)=(\del_{E,h_{P,t}}+\ad t\theta^{\dagger}_{h_{P,t}})
 F_{\alpha_P}
 \in A^1(\Def(E_P,t\theta_P))
\]
is a harmonic $1$-form
for $(E_P,t\theta_P,h_{P,t})$.
(See Proposition \ref{prop;23.4.11.1}.)
Moreover, it satisfies
$\rho(\ttH_P(\nu,t))=\pi_{\ast}(\nu)_{|X_P}$
and
$\ttH_P(\nu,t)\in A^1(\Def(E_P,t\theta_P;C_P))$.

By the asymptotic orthogonality in \cite{Decouple},
we obtain
$\ttH_P(\nu,t)-F_{\nu}=O(e^{-\beta t})\|\nu\|$ for some $\beta>0$
on $X_P(1)\setminus X_P(1/2)$.
There exists a unique $C^{\infty}$-section
$\sigma_P(\nu,t)$
of $\End^{\asym}(E;C_P)$ on $X_P(1)\setminus X_P(1/2)$
such that
$\ttH_P(\nu,t)-F_{\nu}=(\delbar_E+\ad t\theta)\sigma_P(\nu,t)$ on
$X_P(1)\setminus X_P(1/2)$
and that
$|\sigma_P(\nu,t)|_{h_{P,t}}=O(e^{-\epsilon t})\|\nu\|$
for some $\epsilon>0$.
We patch $F_{\nu}$ and $\ttH_P(\nu,t)$ $(P\in D(\theta))$
by using $\sigma_P(\nu,t)$ and a partition of unity,
and we obtain a $1$-cocycle $\ttH'(\nu,t)\in A^1(\Def(E,\theta))$
satisfying the conditions (i--iv).
Note that we can check the condition (iii)
by using $h_{P,t}-h_t=O(e^{-\beta t})$ on $X_P(1)\setminus X_P(1/2)$
for some $\beta>0$
in \cite{Mochizuki-Szabo}.

For each $t>0$, there uniquely exists
$\gamma(\nu,t)\in A^0(\End(E))$ such that
$\tr\gamma(\nu,t)=0$
and
\[
 (\delbar_E+\ad t\theta)^{\ast}_{h_{t},g_X}
 \circ
 (\delbar_E+\ad t\theta)\gamma(\nu,t)
=(\delbar_E+\ad t\theta)^{\ast}_{h_{t},g_X}
\ttH'(\nu,t).
\]
Then,
$\ttH(\nu,t)=\ttH'(\nu,t)-(\delbar_E+\ad t\theta)\gamma(\nu,t)$
are the harmonic representatives of
$\iota^{\hor}_t(\nu)$.
By using the estimate for the first eigenvalues of
the Laplacian operators
$(\delbar_E+\ad t\theta)^{\ast}_{h_t,g_X}
\circ
 (\delbar_E+\ad t\theta)$ (see Proposition \ref{prop;23.4.12.1})
and the regularity for solutions of
the Poisson equation (see Proposition \ref{prop;23.5.5.61}),
we obtain
$|\gamma(\nu,t)|_{h_t}
+|(\delbar_E+\ad t\theta)\gamma(\nu,t)|_{h_t,g_X}
=O(e^{-\epsilon t})$ for some $\epsilon>0$.

\subsubsection{Estimate of the pairings in the horizontal direction}

As the restriction of the estimate (\ref{eq;23.5.8.10})
in the horizontal direction,
we would like to prove that
\begin{equation}
\label{eq;23.5.19.1}
\Bigl|
 g_{H}(\iota^{\hor}_t(\nu),\iota^{\hor}_t(\nu))
 -g_{\semiflat}(\iota^{\hor}_t(\nu),\iota^{\hor}_t(\nu))
 \Bigr|
=O(e^{-\epsilon t})
 \|\nu\|^2
\end{equation}
for some $\epsilon>0$.
It is equivalent to 
\begin{equation}
\label{eq;23.5.12.10}
 2\sqrt{-1}\int_X
 \Bigl(
 h_t\bigl(
 \ttH(\nu,t)^{1,0},
 \ttH(\nu,t)^{1,0}
 \bigr)
 - h_t\bigl(
 \ttH(\nu,t)^{0,1},
 \ttH(\nu,t)^{0,1}
 \bigr)
 \Bigr)
 -2\sqrt{-1}\int_{\Sigma_{\theta}}
 \nu\wedge \nubar
 =O\bigl(e^{-\epsilon t}\bigr)
 \|\nu\|^2.
\end{equation}

Let us explain an outline.
We set $X_1=X\setminus\bigcup_{P}X_P(1)$.
We obtain (\ref{eq;23.5.12.10})
from the following two estimates:
\begin{equation}
\label{eq;23.5.12.11}
 2\sqrt{-1}\int_{X_1}
 h_t\bigl(
 F_{\nu},F_{\nu}
 \bigr)
 -2\sqrt{-1}\int_{\Sigma_{\theta|X_1}}
 \nu\wedge \nubar
 =O\bigl(e^{-\epsilon t}\bigr)
 \|\nu\|^2,
\end{equation}
\begin{multline}
\label{eq;23.5.12.12}
 2\sqrt{-1}\int_{X_P(1)}
 \Bigl(
 h_{P,t}\bigl(
 \ttH_P(\nu,t)^{1,0},
 \ttH_P(\nu,t)^{1,0}
 \bigr)
 - h_{P,t}\bigl(
 \ttH_P(\nu,t)^{0,1},
 \ttH_P(\nu,t)^{0,1}
 \bigr)
 \Bigr)
 -2\sqrt{-1}\int_{\Sigma_{\theta|X_P(1)}}
 \nu\wedge \nubar
\\
 =O\bigl(e^{-\epsilon t}\bigr)
 \|\nu\|^2.
\end{multline}

Let $F_{\nubar}$ denote the $C^{\infty}$-section of
$\End(E)\otimes\Omega^{0,1}$ on $X^{\circ}$
induced by the multiplication of $\nubar$ on $L$.
We have
\[
 \int_{X_1}
 h_t\bigl(
 F_{\nu},F_{\nu}
 \bigr)
-\int_{\Sigma_{\theta|X_1}}
 \nu\wedge \nubar
=\int_{X_1}
 \tr\Bigl(
 F_{\nu}\cdot
 \bigl(
 (F_{\nu})^{\dagger}_{h_t}
-F_{\nubar}
 \bigr)
 \Bigr).
\]
By using the asymptotic orthogonality of $h_t$ on $X_1$
\cite{Decouple},
we obtain
\[
|F_{\nu}|_{h_t,g_X}\leq B_0\|\nu\|,
\quad\quad
\bigl|
(F_{\nu})^{\dagger}_{h_t}
-F_{\nubar}
\bigr|_{h_t,g_X}\leq B_1\exp(-B_2t)\|\nu\|
\]
for some positive constants $B_i$ $(i=0,1,2)$ on $X_1$.
Hence, we obtain the estimate (\ref{eq;23.5.12.11}).

To obtain the estimate (\ref{eq;23.5.12.12}),
we use a variant of
the key idea due to Dumas-Neitzke \cite{Dumas-Neitzke}
efficiently applied by Fredrickson in \cite{Fredrickson2}.
By applying the Stokes formula to harmonic $1$-forms,
we obtain the following equality
(Proposition \ref{prop;23.5.4.20}):
\begin{multline}
 \sqrt{-1}\int_{X_P(1)}
 \Bigl(
 h_{P,t}\bigl(
 \ttH_P(\nu,t)^{1,0},
 \ttH_P(\nu,t)^{1,0}
 \bigr)
 - h_{P,t}\bigl(
 \ttH_P(\nu,t)^{0,1},
 \ttH_P(\nu,t)^{0,1}
 \bigr)
 \Bigr)
=
 \\
 \sqrt{-1}
 \int_{\del X_P(1)}
  \Bigl(
 h_{P,t}\bigl(
 F_{\alpha_P},
 \ttH_P(\nu,t)^{1,0}
 \bigr)
 - h_{P,t}\bigl(
 F_{\alpha_P},
 \ttH_P(\nu,t)^{0,1}
 \bigr)
 \Bigr)
 =\\
 \sqrt{-1} \int_{\del X_P(1)}
 \tr\Bigl(
  F_{\alpha_P}\cdot
 \bigl(\ttH_P(\nu,t)^{1,0}\bigr)^{\dagger}_{h_{P,t}}
-F_{\alpha_P}\cdot
 \bigl(\ttH_P(\nu,t)^{0,1}\bigr)^{\dagger}_{h_{P,t}}
 \Bigr).
\end{multline}
We also have
\[
 \sqrt{-1}
 \int_{\Sigma_{\theta|X_P(1)}}
 \nu\wedge\nubar
 =\sqrt{-1}
 \int_{\del \Sigma_{\theta|X_P(1)}}
 \alpha_P\wedge\nubar
 =\sqrt{-1}
 \int_{\del X_P}
 \tr\bigl(F_{\alpha_P}\cdot F_{\nubar}\bigr).
\]
On a neighbourhood of $\del X_P(1)$
we obtain
\[
 |F_{\alpha_P}|_{h_{P,t}}
 \leq B_0\|\nu\|,
 \quad\quad
\bigl|
 (\ttH_P(\nu,t)^{1,0})^{\dagger}_{h_{P,t}}
 -F_{\nubar}
 \bigr|_{h_{P,t},g_X}
 \leq B_1\exp(-B_2t)\|\nu\|,
\]
\[
|\ttH_P(\nu,t)^{0,1}|_{h_{P,t}g_X}
 \leq
 B_3\exp(-B_4 t)\|\nu\|.
\]
Thus, we obtain the estimate (\ref{eq;23.5.12.12}).

\begin{rem}
A special case of Proposition {\rm\ref{prop;23.5.4.20}}
seems to be closely related with 
the exactness in {\rm\cite[\S9]{Dumas-Neitzke}}
though the author has not studied
the precise relation.
\hfill\qed
\end{rem}

\subsubsection{Harmonic representatives of vertical infinitesimal deformations}

Let $\tau\in
H^1(\Sigma_{\theta},\nbigo_{\Sigma_{\theta}})$.
It is uniquely expressed by
a harmonic $(0,1)$-form on $\Sigma_{\theta}$,
which is also denoted by $\tau$.
Under the natural identification
$\Sigma_{t\theta}$ and $\Sigma_{\theta}$,
we have
$\pi_{\ast}(\nbigo_{\Sigma_{\theta}})
=\nbigh^0(\Def(E,t\theta))$.
Hence, we obtain 
$\iota_t^{\ver}(\tau)\in H^1(X,\Def(E,t\theta))$.
There exists the harmonic representative
$\ttV(\tau,t)\in A^1(\Def(E,t\theta))$
of $\iota_t^{\ver}(\tau)$.

The multiplication of $\tau$ on $L$ induces
a $C^{\infty}$-section $F_{\tau}$ of
$\End(E)\otimes\Omega^{0,1}_{X^{\circ}}$.
It is helpful to prove that there exist
$B_1,B_2>0$ such that
the following holds on
$X\setminus\bigcup_{P\in D(\theta)}X_P(1/2)$:
\begin{equation}
\label{eq;23.5.12.13}
\bigl|
\ttV(\tau,t)-F_{\tau}
\bigr|_{h_t,g_X}
\leq
B_1\exp(-B_2t)\|\tau\|.
\end{equation}

In general,
for any complex manifold $Y$,
let $Y^{\dagger}$ denote the conjugate of $Y$.
There exists the natural isomorphism
$T^{\ast}(X^{\dagger})\simeq (T^{\ast}X)^{\dagger}$.

We have the holomorphic vector bundles
$E^{\dagger}_t=(E,\del_{E,h_t})$
with the Higgs field $\theta_{h_t}^{\dagger}$ on $X^{\dagger}$.
The spectral curve of $(E^{\dagger}_t,t\theta_{h_t}^{\dagger})$
equals $(\Sigma_{t\theta})^{\dagger}\subset
(T^{\ast}X)^{\dagger}=T^{\ast}(X^{\dagger})$.
Let $\pi^{\dagger}_t:(\Sigma_{t\theta})^{\dagger}\to X^{\dagger}$
denote the projection.
There exists a holomorphic line bundle
$L^{\dagger}_t$ on $(\Sigma_{t\theta})^{\dagger}$
with an isomorphism
$(\pi_t^{\dagger})_{\ast}L^{\dagger}_t\simeq E^{\dagger}_t$.

Let $P\in D(\theta)$.
There exists an anti-holomorphic function $\beta_P$ on
$\Sigma_{\theta|X(P)}$
such that
$\delbar\beta_P=\tau_{|\Sigma_{\theta|X_P}}$
and $\beta_P(Q)=0$ for any $Q\in \pi^{-1}(P)$.
We may regard $\beta_P$
as a holomorphic function on
$(\Sigma^{\dagger}_{t\theta})_{|X_P^{\dagger}}$.
The multiplication of $\beta_P$ on $L^{\dagger}_t$
induces the endomorphism
$F_{(E^{\dagger}_t,t\theta^{\dagger}_{h_t}),\beta_P}$
of $E^{\dagger}_t$.
It satisfies
\[
 \bigl(
 \del_{E,h_t}+\ad t\theta^{\dagger}_{h_t}
 \bigr)
 F_{(E^{\dagger}_t,t\theta^{\dagger}_{h_t}),\beta_P}
=0.
\]
We obtain the following harmonic $1$-form
for $(E,t\theta,h_t)$ on $X_P$:
\[
 \ttV_P(\tau,t):=
 (\delbar_E+\ad t\theta)
  F_{(E^{\dagger}_t,t\theta^{\dagger}_{h_t}),\beta_P}.
\]
By the asymptotic orthogonality in \cite{Decouple},
we obtain
$\ttV_P(\tau,t)-F_{\tau}
=O\bigl(e^{-\epsilon t}\bigr)$
on $X_P(1)\setminus X_{P}(1/2)$
for some $\epsilon>0$.
\begin{rem}
We have
$F_{(E^{\dagger}_t,t\theta^{\dagger}_{h_t}),\beta_P}
=(F_{\betabar_P})^{\dagger}_{h_t}$,
where $F_{\betabar_P}$ is the holomorphic endomorphism of 
$(E,\theta)_{|X_P}$ induced by
the holomorphic function $\betabar_P$
on $\Sigma_{\theta|X_P}$. 
\hfill\qed 
\end{rem}

For each $P\in D(\theta)$,
let $\chi_P:X\to \{0\leq a\leq 1\}$
be a $C^{\infty}$-function such that
$\chi_P=1$ on $X_P(1/2)$
and $\chi_P=0$ on $X\setminus X_P(1)$.
We define the $C^{\infty}$ $(0,1)$-form $\tau'$
on $\Sigma_{\theta}$
as follows:
\[
 \tau'=\tau-\sum_{P\in D(\theta)}
 \delbar\bigl(\chi_P\beta_P\bigr).
\]
We obtain the $C^{\infty}$-section $F_{\tau'}$
of $\End(E)\otimes\Omega^{0,1}$
such that
$[F_{\tau'}]=\iota^{\ver}_t(\tau)$.
We set
\[
 \ttV'(\tau,t)
 =F_{\tau'}
 +\sum_{P\in D(\theta)}
 (\delbar_E+\ad t\theta)
 \bigl(
  \chi_P\cdot
  F_{(E^{\dagger}_t,t\theta^{\dagger}_{h_t}),\beta_P}
 \bigr).
\]
Then, $\ttV'(\tau,t)$ satisfies (\ref{eq;23.5.12.13})
for some $B_i>0$.
We have 
$(\delbar_E+\ad t\theta)^{\ast}_{h_t,g_X}\ttV'(\tau,t)
=O\bigl(e^{-\epsilon t}\bigr)\|\tau\|$
on $X$ for some $\epsilon>0$.
Moreover,
$\tr\bigl(
(\delbar_E+\ad t\theta)^{\ast}_{h_t,g_X}\ttV'(\tau,t)\bigr)=0$ holds.
Then, for each $t>0$,
there uniquely exists $\gamma(\tau,t)\in A^0(\End(E))$
such that $\tr\gamma(\tau,t)=0$
and
\[
 (\delbar_E+\ad t\theta)^{\ast}_{h_t,g_X}
 \circ
 (\delbar_E+\ad t\theta)\gamma(\tau,t)
=(\delbar_E+\ad t\theta)^{\ast}_{h_t,g_X}\ttV'(\tau,t).
\]
Again,
by using the estimate for the first eigenvalues of
the Laplacian operators
$(\delbar_E+\ad t\theta)^{\ast}_{h_t,g_X}
\circ
(\delbar_E+\ad t\theta)$ (see Proposition \ref{prop;23.4.12.1})
and the regularity for solutions of
the Poisson equation (see Proposition \ref{prop;23.5.5.61}),
we obtain
$|\gamma(\tau,t)|_{h_t}
+|(\delbar_E+\ad t\theta)\gamma(\tau,t)|_{h_t,g_X}
\leq B_1\exp(-B_2t)\|\tau\|$
for some $B_i>0$.
Because
$\ttV(\tau,t)=\ttV'(\tau,t)-(\delbar_E+\ad t\theta)\gamma(\tau,t)$,
we obtain the estimate (\ref{eq;23.5.12.13})
for $\ttV(\tau,t)$.

\subsubsection{Estimate of pairings with vertical infinitesimal deformations}

To prove the estimate (\ref{eq;23.5.8.10}),
we would like to prove the following estimates
for $\tau\in H^1(\Sigma_{\theta},\nbigo_{\Sigma_{\theta}})$
and $\nu\in H^0(\Sigma_{\theta},K_{\Sigma_{\theta}})$:
\begin{equation}
\label{eq;23.5.19.2}
\Bigl|
 g_{H}(\iota^{\ver}_t(\tau),\iota^{\ver}_t(\tau))
 -g_{\semiflat}(\iota^{\ver}_t(\tau),\iota^{\ver}_t(\tau))
 \Bigr|
=O(e^{-\epsilon t})
 \|\tau\|^2,
\end{equation}
\begin{equation}
\label{eq;23.5.19.3}
\Bigl|
 g_{H}(\iota^{\ver}_t(\tau),\iota^{\hor}_t(\nu))
-g_{\semiflat}(\iota^{\ver}_t(\tau),\iota^{\hor}_t(\nu))
 \Bigr|
=O(e^{-\epsilon t})
 \|\tau\|\cdot\|\nu\|
\end{equation}
for some $\epsilon>0$.
They are equivalent to the following estimates,
respectively:
\begin{equation}
\label{eq;23.5.12.30}
 2\sqrt{-1}\int_X
  \Bigl(
  h_t\bigl(
  \ttV(\tau,t)^{1,0},
  \ttV(\tau,t)^{1,0}
  \bigr)
- h_t\bigl(
  \ttV(\tau,t)^{0,1},
  \ttV(\tau,t)^{0,1}
  \bigr)
  \Bigr)
  -(-2\sqrt{-1})\int_{\Sigma_{\theta}}
  \tau\taubar
=O(e^{-\epsilon t}) \|\tau\|^2,
\end{equation}
\begin{equation}
\label{eq;23.5.12.31}
  2\sqrt{-1}\int_X
  \Bigl(
  h_t\bigl(
  \ttV(\tau,t)^{1,0},
  \ttH(\nu,t)^{1,0}
  \bigr)
- h_t\bigl(
  \ttV(\tau,t)^{0,1},
  \ttH(\nu,t)^{0,1}
  \bigr)
  \Bigr)
=O(e^{-\epsilon t}) \|\tau\|\cdot\|\nu\|.
\end{equation}

Let us explain
an outline of the proof of the estimate (\ref{eq;23.5.12.30}).
Because
$(\delbar_E+\ad t\theta)^{\ast}_{h_t,g_X}\ttV(\tau,t)=0$,
we obtain
\begin{multline}
 \sqrt{-1}\int_X
  \Bigl(
  h_t\bigl(
  \ttV(\tau,t)^{1,0},
  \ttV(\tau,t)^{1,0}
  \bigr)
- h_t\bigl(
  \ttV(\tau,t)^{0,1},
  \ttV(\tau,t)^{0,1}
  \bigr)
  \Bigr)
=-\sqrt{-1}\int_X
  \Bigl(
  h_t\bigl(
  F_{\tau'},
  \ttV(\tau,t)^{0,1}
  \bigr)
  \Bigr)
 \\
=-\sqrt{-1}\int_X
 \tr\Bigl(
 F_{\tau'}\cdot
 \bigl(\ttV(\tau,t)^{0,1}\bigr)^{\dagger}_{h_t}
 \Bigr)
=-\sqrt{-1}\int_{X\setminus\bigcup X_P(1/2)}
 \tr\Bigl(
 F_{\tau'}\cdot
 \bigl(\ttV(\tau,t)^{0,1}\bigr)^{\dagger}_{h_t}
 \Bigr)
\end{multline}
We also have
\[
 \sqrt{-1}\int_{\Sigma_{\theta}}
 \tau\taubar
=\sqrt{-1}\int_{\Sigma_{\theta}}
 \tau'\taubar
=\sqrt{-1}\int_{X}
 \tr\bigl(F_{\tau'}F_{\taubar}\bigr)
=\sqrt{-1}\int_{X\setminus\bigcup X_P(1/2)}
 \tr\bigl(F_{\tau'}F_{\taubar}\bigr).
\]
Note that the support of $F_{\tau'}$ is contained in
$X\setminus\bigcup X_P(1/2)$.
By the estimate (\ref{eq;23.5.12.13})
and the asymptotic orthogonality in \cite{Decouple},
the following holds on $X\setminus\bigcup X_P(1/2)$
for some $B_3,B_4>0$:
\[
\bigl|
 \bigl(\ttV(\tau,t)^{0,1}\bigr)^{\dagger}_{h_t}
 -F_{\taubar}
 \bigr|_{h_t,g_X}
 \leq B_3\exp(-B_4t)\|\tau\|^2.
\]
Then, we obtain the estimate (\ref{eq;23.5.12.30}).
We can obtain the estimate (\ref{eq;23.5.12.31}) similarly.

\paragraph{Acknowledgements}

I thank the reviewer for his/her thorough and careful reading
and helpful comments to improve this manuscript.
This study grew out of my efforts
to understand the attractive works
\cite{Dumas-Neitzke, Fredrickson2, MSWW2}.
I hope that this would be useful
in the further study of this subject.
I thank Qiongling Li and Szilard Szabo
for discussions and collaborations.
This study is partially based on the joint works
\cite{Li-Mochizuki3, Mochizuki-Szabo}.
I am grateful to
Gao Chen,
Laura Fredrickson,
Nianzi Li,
Rafe Mazzeo and Siqi He for some discussions.
In particular,
I appreciate Fredrickson and Mazzeo for sincerely replying
to my inquiries regarding Remark \ref{rem;24.10.23.20}.
I am grateful to Fredrickson, Mazzeo and Szabo for helpful comments
to this manuscript.
I thank Hitoshi Fujioka for discussions.
I thank Yoshifumi Tsuchimoto and Akira Ishii
for their constant encouragements.

I gave talks on the related topics in the conferences 
``The Hitchin system, Langlands duality and mirror symmetry''
in April 2023 in ICMAT,
and ``Minimal Surfaces in Symmetric Spaces'' in May 2023 in IMAG.
The preparation of the talks at these conferences were useful
in this research.
I thank the organizers for giving me the opportunities of the talks.
I gave a series of lectures on this topic
in University of Science and Technology of China
in October 2024.
I thank Bin Xu for the opportunity of lectures,
and all the participants for many stimulating questions.

I am partially supported by
the Grant-in-Aid for Scientific Research (A) (No. 21H04429),
the Grant-in-Aid for Scientific Research (A) (No. 22H00094),
the Grant-in-Aid for Scientific Research (A) (No. 23H00083),
the Grant-in-Aid for Scientific Research (C) (No. 20K03609),
the Grant-in-Aid for Scientific Research (C) (No. 25K06973),
Japan Society for the Promotion of Science.
I am also partially supported by the Research Institute for Mathematical
Sciences, an International Joint Usage/Research Center located in Kyoto
University.

\section{Preliminaries}

\subsection{Some notation}

Let $X$ be any Riemann surface.
We use the notation $K_X$
to denote the sheaf of holomorphic $1$-forms in this section.
The cotangent bundle is denoted by $T^{\ast}X$.
Let $\pi_0:T^{\ast}X\to X$ denote the projection.
Let $\eta$ denote the tautological section
of $\pi_0^{\ast}K_X$,
i.e.,
$\eta$ is the section induced by
the diagonal morphism
$T^{\ast}X\to T^{\ast}X\times_XT^{\ast}X$.

Let $(E,\theta)$ be a Higgs bundle on $X$,
i.e.,
$E$ denotes a holomorphic vector bundle on $X$,
and $\theta$ denotes a holomorphic section of
$\End(E)\otimes K_X$.
The spectral curve of the Higgs bundle $(E,\theta)$
is denoted by $\Sigma_{\theta}$,
i.e.,
$\Sigma_{\theta}$ denotes the $0$-scheme
of the section $\det(\eta\id_E-\theta)$
of $\pi_0^{\ast}\bigl(K_X^{\otimes \rank(E)}\bigr)$.
For any subset $U\subset X$,
we set $\Sigma_{\theta|U}=\Sigma_{\theta}\times_XU$.
We set
\[
 D(\theta):=\bigl\{
 P\in X\,\big|\,|T_P^{\ast}X\cap \Sigma_{\theta}|<\rank(E)
\bigr\}.
\]
The Higgs bundle is called regular semisimple
if $D(\theta)=\emptyset$.
The Higgs bundle is called generically regular semisimple
if $D(\theta)\neq X$.
In this case, $D(\theta)$ is discrete in $X$.
We set $X^{\circ}:=X\setminus D(\theta)$.
The induced Higgs field $\ad(\theta)$ on $\End(E)$
is often denoted by $\theta$
if there is no risk of confusion.

\subsubsection{Local decomposition}
\label{subsection;23.4.15.1}

Suppose that $(E,\theta)$ is generically regular semisimple.
For any $P\in X$,
let $(X_P,z_P)$ denote a holomorphic coordinate neighbourhood of $P$ in $X$
such that $z_P(P)=0$.
We set $X_P^{\circ}=X_P\setminus\{P\}$.
We assume that $X_P$ is simply connected
and that $X_P^{\circ}\cap D(\theta)=\emptyset$.
There exists the decomposition into the connected components
\begin{equation}
\label{eq;23.4.11.10}
 \Sigma_{\theta|X_P}=
 \coprod_{Q\in \Sigma_{\theta|P}}
 (\Sigma_{\theta})_Q,
\end{equation}
where
$(\Sigma_{\theta})_Q\cap T_P^{\ast}X=\{Q\}$.
There exists the decomposition
\begin{equation}
\label{eq;23.4.11.11}
(E,\theta)_{|X_P}
=\bigoplus_{Q\in T_P^{\ast}X\cap \Sigma_{\theta}}
(E_{Q},\theta_Q)
\end{equation}
such that
$\Sigma_{\theta_Q}=(\Sigma_{\theta})_Q$.
If $P\in X^{\circ}$,
we have $\rank E_Q=1$
for any $Q\in \Sigma_{\theta|P}$.

\subsubsection{Filtrations}
\label{subsection;26.1.27.1}

When $\Sigma_{\theta}$ is smooth and there exists a line bundle $L$
on $\Sigma_{\theta}$ such that $E=\pi_{\ast}(L)$,
there exists a decreasing filtration
$\nbigf^j(E_{Q|P})$ $(j=0,\ldots,r_Q-1)$
obtained as the images of
\[
\pi_{\ast}\bigl(L_{|(\Sigma_{\theta})_Q}(-jQ)\bigr)_{|P}
\to
\pi_{\ast}(L_{|(\Sigma_{\theta})_Q})_{|P}.
\]
We also set
$\nbigf^{r_Q}(E_{Q|P})=0$.
We set
$\Gr_{\nbigf}^j(E_{Q|P})=
\nbigf^j(E_{Q|P})\big/
\nbigf^{j+1}(E_{Q|P})$
for $j=0,\ldots,r_Q-1$.
We have
$\Gr_{\nbigf}^j(E_{Q|P})=1$.

\subsection{Deformation complexes of Higgs bundles}
\label{subsection;24.10.22.20}

\subsubsection{Deformation complexes of Higgs bundles}

The deformation complex
$\Def(E,\theta)$
of $(E,\theta)$
is the following complex of sheaves of holomorphic sections:
\[
\End(E)\stackrel{\ad\theta}{\lrarr} \End(E)\otimes K_X.
\]
Here, the first term sits in the degree $0$.

When the Higgs bundle $(E,\theta)$ is equipped with
a non-degenerate symmetric pairing $C$,
let $\End^{\sym}(E;C)$ (resp. $\End^{\asym}(E;C)$)
denote the sheaf of sections $f$ of $\End(E)$
such that $f$ is self-adjoint (resp. anti-self-adjoint)
with respect to $C$.
We obtain the subcomplex
$\Def(E,\theta;C)$
of $\Def(E,\theta)$
given by
\[
 \End^{\asym}(E;C)
 \stackrel{\ad\theta}{\lrarr}
 \End^{\sym}(E;C)\otimes K_X.
\]

\subsubsection{The morphism $\ad(\theta)$
in the case where the spectral curve is smooth}

Suppose that $(E,\theta)$ is generically regular semisimple
and that $\Sigma_{\theta}$ is smooth.
\begin{lem}
\label{lem;23.3.23.2}
The following holds.
\begin{itemize}
 \item The kernel and the cokernel of the morphism
$\ad(\theta):\End(E)\to \End(E)\otimes K_X$ 
are locally free $\nbigo_X$-modules of rank $\rank(E)$.
 \item Let $U$ be any open subset of $X$
       with a holomorphic coordinate $z$.
       Let $f$ be the endomorphism of
       $E_{|U}$ determined by
       $\theta=f\,dz$.
       Then,
       for any $P\in U$,
       $f_{|P}\in\End(E_{|P})$ is regular
       in the following sense:
\begin{equation}
\label{eq;23.3.23.1}
\bigl\{
 g\in\End(E_{|P})\,\big|\,
 [f_{|P},g]=0
\bigr\}
=\bigoplus_{j=0}^{\rank(E)-1}\cnum f_{|P}^j.
\end{equation}
In particular, we obtain
$\dim\bigl\{g\in\End(E_{|P})\,\big|\,
 [f_{|P},g]=0\bigr\}
=\rank(E)$.
 As a result, 
\[
\Ker(\ad\theta)_{|U}
=\bigoplus_{j=0}^{\rank(E)-1} \nbigo_U f^j.
\]
\end{itemize}
\end{lem}
\pf
Let us prove (\ref{eq;23.3.23.1}).
Let $P\in X$
with a holomorphic coordinate neighbourhood
$(X_P,z_P)$ as in \S\ref{subsection;23.4.15.1}.
We obtain the coordinate
$(z,\xi)$ on $T^{\ast}X_P$.
There exists the decomposition (\ref{eq;23.4.11.11})
of $(E,\theta)_{|X_P}$.
Each connected component $(\Sigma_{\theta})_Q$ is smooth.
The eigenvalues of
$f_{Q|P}$ are $\xi(Q)$.
Hence, it is enough to consider the case
where $|\Sigma_{\theta|P}|=1$ and $X=X_P$.
Moreover, we may assume that
$\Sigma_{\theta|P}=\{0\}$ in $T_P^{\ast}X$
by replacing
$\theta$ with $\theta-\frac{1}{r}\tr(\theta)\id_E$.
We set $\Xtilde=\Sigma_{\theta}$.
Let $\pi:\Xtilde\to X$ denote the projection.
We identify
$T^{\ast}X=X\times\cnum$ by using the coordinate $z_P$.
It induces a holomorphic coordinate of $\Xtilde$.
There exists an isomorphism
$E\simeq\pi_{\ast}(\nbigo_{\Xtilde})$
under which the endomorphism $f$ is induced by
the multiplication of $\xi$ on $\nbigo_{\Xtilde}$.
We obtain the frame
$v_j=\pi_{\ast}(\xi^j)$ $(j=0,\ldots,r-1)$
of $E$.
Because
$f_{P}(v_{j|P})=v_{j+1|P}$ $(j<r-1)$
and $f_P(v_{r-1|P})=0$,
we obtain (\ref{eq;23.3.23.1}).
The other claims follow from (\ref{eq;23.3.23.1}).
\hfill\qed

\begin{cor}
Let $\nbigf$ be the filtration in {\rm\S\ref{subsection;26.1.27.1}.}
Let $N_Q$ denote the nilpotent part of $f_{Q|P}$.
Then, we have $N_Q(\nbigf^j(E_{P|Q}))\subset \nbigf^{j+1}(E_{P|Q})$
for $j=0,\ldots,r_Q-1$.
Moreover,
$N_Q$ induces an isomorphism
$\Gr^{\nbigf}_j(E_{Q|P})\simeq
\Gr^{\nbigf}_{j+1}(E_{Q|P})$
for $j=0,\ldots,r_Q-1$.
\hfill\qed
\end{cor}

There exists the natural inclusion
$\Ker(\ad\theta)\otimes K_X\subset
\End(E)\otimes K_X$.
It induces a morphism
\begin{equation}
\label{eq;23.4.23.1}
 \Ker(\ad\theta)\otimes K_X
 \lrarr
 \Cok(\ad\theta).
\end{equation}
The following lemma is obvious.
\begin{lem}
The morphism {\rm(\ref{eq;23.4.23.1})}
induces an isomorphism  
\begin{equation}
 \label{eq;23.4.23.2}
 (\Ker(\ad\theta)\otimes K_X)_{|X\setminus D(\theta)}
\simeq
 \Cok(\ad\theta)_{|X\setminus D(\theta)}.
\end{equation}
As a result, we obtain a splitting
\begin{equation}
\label{eq;23.4.23.3}
 \Cok(\ad\theta)_{|X\setminus D(\theta)}
 \to
 \bigl(
 \End(E)\otimes K_X
 \bigr)
 _{|X\setminus D(\theta)}.
\end{equation}
\hfill\qed
\end{lem}

\subsubsection{Deformation complex in the case
where $\Sigma_{\theta}$ is smooth}

We continue to assume that
 $(E,\theta)$ is generically regular semisimple
 and that $\Sigma_{\theta}$ is smooth.
Let $\pi:\Sigma_{\theta}\to X$ denote the projection.
Let $\eta\in H^0(T^{\ast}X,\pi_0^{\ast}(K_X))$
denote the tautological $1$-form.
The restriction to $\Sigma_{\theta}$
is also denoted by $\eta$. 
There exists a holomorphic line bundle $L$ on $\Sigma_{\theta}$
with $(\pi_{\ast}(L),\pi_{\ast}(\eta))\simeq(E,\theta)$,
where $\pi_{\ast}(\eta)$ is the Higgs field of
$\pi_{\ast}(L)$
induced by $\eta:L\to L\otimes\pi^{\ast}(K_X)$.
We obtain the following refinement of Lemma \ref{lem;23.3.23.2}.

\begin{prop}
\label{prop;23.3.23.10}
There are natural isomorphisms:
\[
 \nbigh^0(\Def(E,\theta))
=\Ker\ad(\theta)\simeq
 \pi_{\ast}\nbigo_{\Sigma_{\theta}},
 \quad\quad
\nbigh^1(\Def(E,\theta))=
 \Cok(\ad(\theta))
\simeq
\pi_{\ast}K_{\Sigma_{\theta}}.
\]
\end{prop}
\pf
Because $L$ is an $\nbigo_{\Sigma_{\theta}}$-module,
$E=\pi_{\ast}L$ is naturally
a $\pi_{\ast}\nbigo_{\Sigma_{\theta}}$-module.
Hence, we obtain the induced morphism
\[
\pi_{\ast}(\nbigo_{\Sigma_{\theta}})
\to
\Ker\ad\theta\subset\End(E).
\]
We can easily check that it is an isomorphism
by using Lemma \ref{lem;23.3.23.2}.
Thus, we obtain
$\pi_{\ast}\nbigo_{\Sigma_{\theta}}\simeq
\nbigh^0(\Def(E,\theta))$.

Let $P\in X^{\circ}$.
There exist a coordinate neighbourhood $(X_P,z_P)$
and a decomposition of the Higgs bundle (\ref{eq;23.4.11.11})
as in \S\ref{subsection;23.4.15.1}.
Let $\psi$ be any holomorphic section of $\End(E)\otimes K_X$
on $X_P$.
We have the decomposition
$\psi=\sum \psi_{Q,Q'}$,
where $\psi_{Q,Q'}\in\Hom(E_{Q'},E_{Q})\otimes K_{X}$.
We may naturally regard
$\sum_Q\psi_{Q,Q}$
as a section of
$\pi_{\ast}(K_{\Sigma_{\theta}})_{|X_P}$.
This induces an isomorphism
$\rho_P:\Cok(\ad\theta)_{|X_P}
\simeq
 \pi_{\ast}(K_{\Sigma_{\theta}})_{|X_P}$.
By considering the isomorphisms $\rho_P$
for any $P\in X^{\circ}$,
we obtain 
the isomorphism
$\rho_{X^{\circ}}:
\Cok\bigl(\ad\theta\bigr)_{|X^{\circ}}
\simeq
 \pi_{\ast}(K_{\Sigma_{\theta}})_{|X^{\circ}}$.
Let us prove that
$\rho_{X^{\circ}}$
extends to an isomorphism
$\rho_X:\Cok(\ad\theta)\simeq \pi_{\ast}(K_{\Sigma_{\theta}})$.
It is enough to consider the case
where
$X$ is a neighbourhood $U$ of $0$ in $\cnum$
and $\Sigma_{\theta|0}=\{0\}$.
Let $f$ be the endomorphism of $E$ determined by
$\theta=f\,dz$.
We set $U^{\circ}=U\setminus\{0\}$
and $\Utilde=\Sigma_{\theta}$.

\begin{lem}
\label{lem;23.3.23.3}
The isomorphism $\rho_{U^{\circ}}$
extends to a morphism
$\rho_U:\Cok(\ad\theta)\to \pi_{\ast}(K_{\Utilde})$. 
\end{lem}
\pf
We set $r:=\rank(E)$.
Let $\varphi_r:\cnum\to\cnum$ be given by
$\varphi_r(\zeta)=\zeta^r$.
We set $U^{(r)}:=\varphi_r^{-1}(U)$.
The induced morphism
$U^{(r)}\to U$ is also denoted by $\varphi_r$.
There exists the isomorphism
$\gamma:U^{(r)}\simeq \Utilde$
such that
$\varphi_r=\pi\circ\gamma$.

There exists the line bundle $L$ on $\Utilde$
with an isomorphism
$(\pi_{\ast}(L),\pi_{\ast}(\eta))\simeq (E,\theta)$.
We obtain the line bundle $L'=\gamma^{\ast}(L)$ on $U^{(r)}$.
Let $e$ be a frame of $L'$ on $U^{(r)}$.
By using the isomorphism $E\simeq \varphi_{r\ast}(L')$,
we obtain the frame
$u_j=\varphi_{r\ast}(\zeta^je)$ $(j=0,\ldots,r-1)$ of $E$.
Let $\omega=\exp(2\pi\sqrt{-1}/r)$.
We set
\[
 v_k=\sum_{j=0}^{r-1}
 \omega^{-jk}\zeta^{-j}\varphi_r^{\ast}(u_j).
\]
They induce a frame of
$\varphi_r^{\ast}(E)(\ast 0)$.

Let $\alpha$ be a holomorphic function $\alpha$ on $U^{(r)}$.
The multiplication of $\alpha$ to $L'$
induces an endomorphism $F_{\alpha}$ of $E$.
We have
\[
 \varphi_r^{\ast}(F_{\alpha})(v_k)=\alpha(\omega^k\zeta)v_k.
\]

By using the identification
$T^{\ast}U=U\times\cnum=\{(z,\xi)\in U\times \cnum\}$,
the endomorphism $f$ is identified with $F_{\gamma^{\ast}\xi}$.
Hence,
the decomposition
\[
 \varphi_r^{\ast}E_{|U^{(r)}\setminus\{0\}}
=\bigoplus_{i=1}^r
 \nbigo_{U^{(r)}\setminus\{0\}}
 v_i
\]
is the eigen decomposition of
$\varphi_r^{\ast}(f)$.

Let $\vecu=(u_0,\ldots,u_{r-1})$
and $\vecv=(v_0,\ldots,v_{r-1})$.
For any holomorphic section $g$ of $\End(E)\otimes K_U$
on $U$,
we obtain the matrix valued holomorphic $1$-form $G$ on $U$
determined by $g(\vecu)=\vecu G$.
Let $B$ be the matrix valued meromorphic function
on $U^{(r)}$
determined by $\vecv=\varphi_r^{\ast}(\vecu)B$.
We have
\[
 \varphi_r^{\ast}(g)(\vecv)=\vecv (B^{-1}\varphi_r^{\ast}(G)B).
\]
We can check that 
the entries of $B^{-1}\varphi_r^{\ast}(G)B$
are holomorphic $1$-forms on $U^{(r)}$
by using
\[
 B_{j,k}=
 \omega^{-jk}\zeta^{-j},
 \quad
 \quad
 (B^{-1})_{k,j}=\frac{1}{r}\omega^{jk}\zeta^j,
 \quad\quad
 \varphi_r^{\ast}(dz)
 =r\zeta^{r-1}d\zeta,
 \quad\quad
 (0\leq j,k\leq r-1).
\]
Then, the claim of Lemma \ref{lem;23.3.23.3} follows.
\hfill\qed

\begin{lem}
\label{lem;23.5.17.1}
Let $V$ be an $r$-dimensional complex vector space.
Let $N$ be a nilpotent endomorphism 
such that $N^{r-1}\neq 0$.
For any $H\in \sum_{j=1}^{r-1}\cnum N^j$,
there exists $G\in \End(V)$
such that $[N,G]=H$.
\end{lem}
\pf
Let $\minisl(2,\cnum)$ denote the Lie algebra of
trace-free $(2\times 2)$-matrices.
Let
\[
 h=\left(
 \begin{array}{cc}
  1 & 0 \\ 0 & -1
 \end{array}
 \right),
 \quad
 e=\left(
  \begin{array}{cc}
   0 & 1 \\ 0 & 0
  \end{array}
  \right),
  \quad
 f=\left(
 \begin{array}{cc}
  0 & 0 \\ 1 & 0
 \end{array}
 \right).
\]
There exists an irreducible representation
$\kappa$ of $\minisl(2,\cnum)$ on $V$
such that $\kappa(f)=N$.
There exists the decomposition
$\End(V)=\cnum\cdot \id_V\oplus \End_0(V)$,
where $\End_0(V)$ denotes
the space of the trace-free endomorphisms of $V$.
It is compatible with the induced $\minisl(2,\cnum)$-action.
It is well known that
$\End_0(V)$ does not contain a trivial representation of
$\minisl(2,\cnum)$,
and $N^j$ $(j=1,\ldots,r-1)$ are contained in
$\End_0(V)\cap\Ker\ad(N)$.
Then, the claim is clear.
\hfill\qed

\vspace{.1in}

We obtain Proposition \ref{prop;23.3.23.10}
from the following lemma.

\begin{lem}
\label{lem;23.3.23.11}
$\pi_{\ast}(K_{\Utilde})$
equals the image of $\rho_U$.
\end{lem}
\pf
Let $\nu$ be a holomorphic $1$-form on $\Utilde$.
The multiplication of $\nu$
induces a section
$F_{\nu}$ of
$\End(E)\otimes K_U(\{0\})$.
We have $[f,F_{\nu}]=0$.
It implies
$[f_{|0},\Res_0(F_{\nu})]=0$.
It is easy to see that
$\Res_0(F_{\nu})\in\sum_{j=1}^{r-1}\cnum f_{|0}^j$.
By Lemma \ref{lem;23.5.17.1},
there exists
$g_0\in\End(E_{|0})$
such that
$[f,g_0]=\Res_0(F_{\nu})$.
Let $g$ be a holomorphic section of $\End(E)$
such that $g_{|0}=g_0$.
Then,
$\Ftilde=F_{\nu}-[\theta,z^{-1}g]$
is a holomorphic section of
$\End(E)\otimes K_U$,
and $\rho_U(\Ftilde)=\pi_{\ast}(\nu)$.
Thus, we obtain
Lemma \ref{lem;23.3.23.11}
and Proposition \ref{prop;23.3.23.10}.
\hfill\qed

\begin{cor}
\label{cor;23.4.11.20}
If moreover $(E,\theta)$ is equipped with a non-degenerate
symmetric pairing $C$,
then we obtain
\[
 \nbigh^0\bigl(\Def(E,\theta;C)\bigr)=0,
 \quad\quad
 \nbigh^1\bigl(\Def(E,\theta;C)\bigr)
 \simeq
 \pi_{\ast}K_{\Sigma_{\theta}}.
\]
\end{cor}
\pf
Let $U$ be any open subset of $X$
with a holomorphic coordinate $z$.
Let $f$ be the endomorphism of $E$
determined by $\theta=f\,dz$.
Note that $f$ is self-adjoint with respect to $C$.
Because
$\Ker(\ad\theta)=\bigoplus_{j=0}^{r-1}\nbigo_Uf^j$,
we obtain
$\nbigh^0\bigl(\Def(E,\theta;C)\bigr)_{|U}=0$.

Let $\tau$ be a holomorphic section of
$\End(E)\otimes K_X$
such that $\tau$ is mapped to $0$
in $\Cok(\ad(\theta))\simeq \pi_{\ast}(K_{\Sigma_{\theta}}^1)$.
There exists a holomorphic section $s$ of $\End(E)$
such that $\tau=[\theta,s]$.
Let $\tau^{\lor}$ (resp. $s^{\lor}$)
denote the adjoint of $\tau$ (resp. $s$) with respect to $C$.
Because $\theta^{\lor}=\theta$,
we have
$\tau^{\lor}=-[\theta,s^{\lor}]$.
If $\tau=\pm\tau^{\lor}$,
we have
$\tau=\frac{1}{2}[\theta,s\mp s^{\lor}]$.
We also note that 
any holomorphic section of $\End^{\asym}(E;C)\otimes K_X$
is mapped to $0$ in
$\Cok(\ad(\theta))$.
Then, we obtain
$\nbigh^1\bigl(\Def(E,\theta;C)\bigr)
 \simeq
 \pi_{\ast}K_{\Sigma_{\theta}}$.
\hfill\qed

\subsubsection{Trace morphisms}
\label{subsection;24.10.22.1}

As the composition of
$\pi_{\ast}(\nbigo_{\Sigma_{\theta}})
=\nbigh^0(\Def(E,\theta))
\to
\End(E)$
and 
$\tr:\End(E)\to\nbigo_X$,
we obtain the morphism
\[
\tr:\pi_{\ast}(\nbigo_{\Sigma_{\theta}})
\to \nbigo_X.
\]
It equals the standard trace map.
Similarly,
for any discrete subset $S\subset X$,
we obtain
\[
 \tr:\pi_{\ast}
 \Bigl(\nbigo_{\Sigma_{\theta}}\bigl(\ast\pi^{-1}(S)\bigr)\Bigr)
\to \nbigo_X(\ast S).
\]

We obtain the morphism
$\End(E)\otimes K_X\to K_X$ 
by taking the trace.
It induces a morphism of complexes
$\Def(E,\theta)\to K_X[1]$.
The induced morphism
$\pi_{\ast}(K_{\Sigma_{\theta}})
=\nbigh^1(\Def(E,\theta))
\lrarr
K_X$
equals the trace morphism
\[
\tr:\pi_{\ast}(K_{\Sigma_{\theta}})
\to K_X
\]
defined as follows.
We obtain
$\pi_{\ast}(K_{\Sigma_{\theta}})_{|X^{\circ}}
 \to
 \bigl(
\End(\pi_{\ast}(\nbigo_{\Sigma_{\theta}}))
 \otimes K_{X}
 \bigr)_{|X^{\circ}}$
as in {\rm\S\ref{subsection;23.4.11.3}} below.
It extends to
\[
\pi_{\ast}(K_{\Sigma_{\theta}})
 \to
 \End(\pi_{\ast}(\nbigo_{\Sigma_{\theta}}))
 \otimes K_{X}(D(\theta)).
\]
Hence, we obtain
$\tr:\pi_{\ast}(K_{\Sigma_{\theta}})
\to K_X(D(\theta))$.
It is easy to check that
the image is contained in $K_X$,
i.e,
$\tr$ induces a morphism
$\tr:\pi_{\ast}(K_{\Sigma_{\theta}})
\to K_X$.

\subsubsection{Locally defined non-degenerate symmetric pairings}
\label{subsection;23.5.3.3}

Let $P\in D(\theta)$.
There exists a decomposition of the Higgs bundle
(\ref{eq;23.4.11.11}).
Let $r_Q:=\rank(E_Q)$,
which equals the ramification index of $\pi$ at $Q$.
Let $C_{L,Q}:
L_{|(\Sigma_{\theta})_Q}\otimes
L_{|(\Sigma_{\theta})_Q}\simeq
\nbigo_{(\Sigma_{\theta})_Q}\bigl(-(r_Q-1)Q\bigr)$
be a symmetric pairing.
Let us recall that
it induces a non-degenerate symmetric pairing
of the Higgs bundle $(E_Q,\theta_Q)$
as studied in \cite{Mochizuki-Szabo}.
Indeed,
we obtain a symmetric pairing
$C_Q:E_Q\otimes E_Q\to \nbigo_X$
as the composition of the natural morphisms
\[
\begin{CD}
E_Q\otimes E_Q
@>{\pi_{\ast}(C_{L,Q})}>>
\pi_{\ast}\Bigl(\nbigo_{\Sigma_Q}\bigl(-(r_Q-1)Q\bigr)\Bigr)
 @>{a}>>
 \nbigo_X.
\end{CD}
\]
Here, the morphism $a$ is induced by
the trace map $\tr$.

\begin{lem}
The restriction of
$C_{Q|P}$ to
 $\nbigf^{j}(E_{Q|P})
 \otimes
 \nbigf^{k}(E_{Q|P})$
is $0$ 
if $k+j>r_Q-1$.
The induced pairings
$\Gr_{\nbigf}^j(E_{Q|P})
 \otimes
 \Gr_{\nbigf}^{r_Q-1-j}(E_{Q|P})
 \lrarr
 \cnum$
are perfect. 
In particular, $C_{Q|P}$ is non-degenerate.
\hfill\qed
\end{lem}

We set $C_P=\bigoplus C_Q$
and $C_{P|P}=\bigoplus C_{Q|P}$.
Let $\nu$ be a holomorphic $1$-form on $\Sigma_{\theta}$.
The restriction
$\nu^{\circ}=\nu_{|X^{\circ}}$
induces a holomorphic section
$F_{\nu^{\circ}}$
of $(\End(E)\otimes K_X)_{|X^{\circ}}$
as in \S\ref{subsection;23.4.11.3} below.
It induces a holomorphic section
$F_{\nu}$ of $\End(E)\otimes K_X(D)$.
Because $F_{\nu}$ is self-adjoint with respect to $C_P$,
$\Res_P(F_{\nu})$ is self-adjoint with respect to $C_{P|P}$.
We also note that $\Res_P(F_{\nu})$ is nilpotent,
and that there exists the decomposition
$\Res_P(F_{\nu})=\sum\Res_P(F_{\nu})_Q$,
where $\Res_P(F_{\nu})_Q$
are endomorphisms of $E_{Q|P}$.
Let $f_P$ be the endomorphism of $E_{|X_P}$
determined by $\theta_{|X_P}=f_P\,dz_P$.
By using Corollary \ref{cor;23.4.11.20},
we obtain the following lemma.
\begin{lem}
\label{lem;23.5.3.2}
There uniquely exists $g_{P}\in \End(E_{|P})$
such that
(i) $[f_{P|P},g_P]=\Res_P(F_{\nu})$,
(ii) $g_P$ is anti-self-adjoint with respect to $C_{P|P}$.
\hfill\qed
\end{lem}

\subsubsection{Dolbeault resolutions}
\label{subsection;23.5.12.2}

We continue to assume that
$(E,\theta)$ is generically regular semisimple
and that $\Sigma_{\theta}$ is smooth.
Let $\nbiga^{\bullet}(\Def(E,\theta))$
denote the Dolbeault resolution of $\Def(E,\theta)$.
The space of global sections
of $\nbiga^j(\Def(E,\theta))$
is denoted by $A^j(\End(E))$
which is naturally identified with
the space of $C^{\infty}$ $\End(E)$-valued $j$-forms.
If moreover
$(E,\theta)$ is equipped with
a non-degenerate symmetric pairing $C$,
let $\nbiga^{\bullet}(\Def(E,\theta;C))$
denote the Dolbeault resolution of $\Def(E,\theta;C)$,
and the space of global sections of
$\nbiga^j(\Def(E,\theta;C))$
is denoted by
$A^j(\End(E),\theta;C)$.

The Cauchy operator of $E$ is denoted by $\delbar_E$.
The induced Cauchy operators of $\End(E)$ is also denoted
by $\delbar_{E}$.

\begin{lem}
\label{lem;23.5.3.4}
Suppose that $\tau\in A^1(\End(E))$ is a $1$-cocycle,
i.e., $(\delbar_{E}+\ad\theta)\tau=0$.
Then, $\tau$ induces a holomorphic section of
$\nbigh^1(\Def(E,\theta))\simeq\pi_{\ast}K_{\Sigma_{\theta}}$.
 The induced holomorphic section is denoted by
 $\rho(\tau)$.
\end{lem}
\pf
There exists the natural morphism
$\Def(E,\theta)\to \nbigh^1(\Def(E,\theta))[-1]$.
By considering the induced morphism of
the Dolbeault resolutions,
we obtain the claim of the lemma.
\hfill\qed

\vspace{.1in}
For any section $\tau$ of $A^{j}(\Def(E,\theta))$,
let $\Supp(\tau)$ denote the support of $\tau$.

\begin{prop}
\label{prop;23.4.15.20}
Let $\tau\in A^1(\Def(E,\theta;C))$ such that
$(\delbar_{E}+\ad\theta)\tau=0$
and that $\rho(\tau)=0$.
Then, there exists a unique section $\sigma$
of $A^0(\Def(E,\theta;C))$ 
such that $(\delbar_{E}+\ad\theta)\sigma=\tau$
and that $\Supp(\sigma)=\Supp(\tau)$.
\end{prop}
\pf
Let us begin with a special case.
\begin{lem}
\label{lem;23.4.15.11}
Let $\tau\in A^1(\Def(E,\theta;C))$ such that
(i) $(\delbar_{E}+\ad\theta)\tau=0$,
(ii) $\Supp\tau\subset X^{\circ}$,
(iii) $\rho(\tau)=0$.
Then, there exists a unique section $\sigma$
of $A^0(\Def(E,\theta;C))$ 
such that $(\delbar_{E}+\ad\theta)\sigma=\tau$
and that $\Supp(\sigma)=\Supp(\tau)$.
\end{lem}
\pf
By the conditions,
there exists a unique $C^{\infty}$-section $\sigma$
of $\End^{\asym}(E;C)$
such that
(i) $\ad(\theta)(\sigma)=\tau^{1,0}$,
(ii) $\Supp(\sigma)=\Supp(\tau^{1,0})$.
Because
$\ad(\theta):\End^{\asym}(E;C)\to \End^{\sym}(E;C)\otimes K_X$
is a monomorphism on $X^{\circ}$,
it is easy to check that
$\tau^{0,1}-\delbar_{E}\sigma=0$.
\hfill\qed

\vspace{.1in}

Let $U$ be any non-compact connected open subset of $X$.
Let $\tau\in A^1(\Def(E,\theta;C)_{|U})$ such that
$(\delbar_{E}+\ad\theta)\tau=0$
and that $\rho(\tau)=0$.

\begin{lem}
\label{lem;23.4.15.10}
There exists a unique section $\sigma$
of $A^0(\Def(E,\theta;C)_{|U})$ 
such that $(\delbar_{E}+\ad\theta)\sigma=\tau$
and that $\Supp\sigma=\Supp\tau$.
\end{lem}
\pf
Because $U$ is non-compact,
there exists a $C^{\infty}$-section $\sigma_1$
of $\End^{\asym}(E;C)$
such that $\delbar_E\sigma_1=\tau^{0,1}$.
We set
$\tau_1=\tau-(\delbar_E+\ad\theta)\sigma_1$.
Then, $\tau_1$ is a holomorphic section of
$\End^{\sym}(E;C)\otimes K_X$
which is mapped to $0$ in $\Cok(\ad(\theta))$.
Because $U$ is non-compact,
there exists a unique holomorphic section
$\sigma_2$ of $\End^{\asym}(E;C)$
such that $\ad(\theta)(\sigma_2)=\tau_1$.
We set $\sigma=\sigma_1+\sigma_2$,
and then $(\delbar_E+\ad\theta)\sigma=\tau$ holds.
By using $\ad(\theta)(\sigma)=\tau^{1,0}$
we obtain the uniqueness of $\sigma$
and $\Supp\sigma=\Supp\tau$.
\hfill\qed

\vspace{.1in}

Let $N$ be an open neighbourhood of $D$
such that each connected component is non-compact.
By applying Lemma \ref{lem;23.4.15.10} to $\tau_{|N}$,
there exists $\sigma_N\in A^0(\Def(E,\theta;C)_{|N})$
such that
$(\delbar_E+\ad\theta)\sigma_N=\tau_{|N}$.
Let $N'$ be a neighbourhood of $D$ in $N$
such that the closure of $N'$ in $X$ is contained in $N$.
Let $\chi:X\to [0,1]$ be a $C^{\infty}$-function
such that $\chi=1$ on $N'$ and $\chi=0$ on $X\setminus N$.
We set
$\tau_1=\tau-(\delbar_E+\ad\theta)(\chi\sigma_N)$.
Because $\Supp(\tau_1)\subset X^{\circ}$,
there exists $\sigma_1\in A^0(\Def(E,\theta;C))$
such that $(\delbar_E+\ad\theta)\sigma_1=\tau_1$
and that $\Supp(\sigma_1)\subset\Supp(\tau_1)$.
We set
$\sigma=\chi_N\sigma_N+\sigma_1$.
Then, we have
$(\delbar_E+\ad\theta)\sigma=\tau$.
By using $\ad(\theta)(\sigma)=\tau^{1,0}$,
we obtain the uniqueness of $\sigma$
and $\Supp\sigma=\Supp\tau$.
\hfill\qed

\subsubsection{Appendix}
\label{subsection;23.4.11.3}

Assume that $(E,\theta)$ is regular semisimple.
There exists a line bundle $L$ on $\Sigma_{\theta}$
with an isomorphism
$(\pi_{\ast}L,\pi_{\ast}\eta)\simeq(E,\theta)$.
For any simply connected open subset $U\subset X$,
there exists the decomposition
into Higgs bundles of rank $1$:
\begin{equation}
\label{eq;23.3.21.1}
 (E,\theta)_{|U}=\bigoplus_{i=1}^{\rank E} (E_{U,i},\theta_{U,i}).
\end{equation}

Let $\nu$ be a $(p,q)$-form on $\Sigma_{\theta}$.
The multiplication of $\nu$ on $L$
induces a $C^{\infty}$-section $F_{\nu}$ of $A^{p,q}(\End(E))$.
\begin{lem}
\mbox{{}}
\begin{itemize}
 \item
 $F_{\nu}$ is diagonal
with respect to the decomposition {\rm(\ref{eq;23.3.21.1})}.      
 \item
 If $\nu$ is a holomorphic function,
then $F_{\nu}$ satisfies
$(\delbar_E+\ad\theta)F_{\nu}=0$.
In particular,
$(\del_{E,h}+\ad\theta_h^{\dagger})F_{\nu}$ 
is a harmonic $1$-form
of a harmonic bundle $(E,\theta,h)$.
  (See Proposition {\rm\ref{prop;23.4.11.1}} below.)
 \item
If $\nu$ is a holomorphic $(1,0)$-form,
then $\delbar_EF_{\nu}=0$.
If $\nu$ is a $(0,1)$-form, then 
$\ad\theta(F_{\nu})=0$.
\hfill\qed      
\end{itemize}
\end{lem}

\subsection{Infinitesimal deformations and relative Higgs bundles}

\subsubsection{Infinitesimal deformations of a Higgs bundle}

For any open subset $U\subset X$,
we set 
$\nbigo_{U^{[1]}}=\nbigo_U[\epsilon]/\epsilon^2$.
We obtain the ringed space
$U^{[1]}=(U,\nbigo_{U^{[1]}})$.
Let $\iota_U:U\to U^{[1]}$ be the morphism of the ringed spaces
defined by the identity map on $U$
and $\iota_U^{\ast}(f_0+f_1\epsilon)=f_0$.
Let $q_U:U^{[1]}\to U$ be the morphism of the ringed spaces
defined by the identity map on $U$
and $q_U^{\ast}(f)=f$.

An infinitesimal deformation of a Higgs bundle $(E,\theta)$ on $X$
is a locally free
$\nbigo_{X^{[1]}}$-module $\Etilde$
with a morphism
$\thetatilde:\Etilde\to\Etilde\otimes q^{\ast}K_{X}$
such  that $\iota^{\ast}(\Etilde,\thetatilde)=(E,\theta)$.
The following lemma is well known.
(See \cite{Nitsure, Yokogawa-infinitesimal}.)

\begin{lem}
$H^1(X,\Def(E,\theta))$
equals the space of
the isomorphism classes of infinitesimal deformations of $(E,\theta)$.
\end{lem}
\pf
We just recall
how an infinitesimal deformation of $(E,\theta)$
induces an element of $H^1(X,\Def(E,\theta))$.
We set $n=\rank E$.
Let $X=\bigcup_{i=0}^m U_i$ be an open covering of $X$
such that $U_i$ are Stein.
Recall that any non-compact Riemann surface is Stein.

Let $(\Etilde,\thetatilde)$ be an infinitesimal deformation
of $(E,\theta)$.
There exist isomorphisms
$a_i:q^{\ast}(E_{|U_i})\simeq \Etilde_{|U_i}$.
Let $\psi_{i,j}$ denote the composition of the following isomorphisms:
\[
 q_{U_i\cap U_j}^{\ast}(E_{|U_i\cap U_j})
 \simeq
 q_{U_j}^{\ast}(E_{|U_j})_{|(U_i\cap U_j)^{[1]}}
 \stackrel{a_j}{\simeq}
 \Etilde_{|(U_i\cap U_j)^{[1]}}
 \stackrel{a_i^{-1}}{\simeq}
 q_{U_i}^{\ast}(E_{|U_i})_{|(U_i\cap U_j)^{[1]}}
 \simeq
 q_{U_i\cap U_j}^{\ast}(E_{|U_i\cap U_j}).
\]
There exist $\rho_{i,j}\in H^0(U_i\cap U_j,\End(E))$
such that
$\psi_{i,j}=\id+\epsilon \rho_{i,j}$.
We have
$\rho_{j,k}-\rho_{i,k}+\rho_{i,j}=0$
on $U_i\cap U_j\cap U_k$.
We also obtain
$\kappa_i\in H^0(U_i,\End(E)\otimes K_{U_i})$
such that
$a_i^{-1}\circ\thetatilde_{|U_i^{[1]}}\circ a_i=
q_{U_i}^{\ast}(\theta)+\epsilon \kappa_i$.
We have
$\kappa_j-\kappa_i-\ad(\theta)\rho_{i,j}=0$
on $U_i\cap U_j$.
Then, the tuple
\[
 (\{\rho_{i,j}\},\{\kappa_i\})
 \in \bigoplus_{i,j} H^0(U_i\cap U_j,\End(E))
 \oplus
 \bigoplus_{i}H^0(U_i,\End(E)\otimes K_{U_i})
\]
gives a $1$-cocycle of
the \v{C}ech cohomology of $\Def(E,\theta)$.

Let $(\{\rhotilde_{i,j}\},\{\kappatilde_i\})$
is a $1$-cocycle induced by
different trivializations $\atilde_i$.
There exist $b_i\in H^0(U_i,\End(E))$
such that $\atilde_i=a_i\circ(1+\epsilon b_i)$.
Then, we have
$\rhotilde_{i,j}=\rho_{i,j}+b_j-b_i$ on $U_i\cap U_j$
and $\kappatilde_i=\kappa_i+\ad(\theta)b_i$ on $U_i$.
Hence,
$(\{\rhotilde_{i,j}\},\{\kappatilde_i\})$
is cohomologous to
$(\{\rho_{i,j}\},\{\kappa_i\})$.
\hfill\qed

\begin{rem}
\label{rem;23.4.19.1}
Let $(\nbigf^{\bullet},d_{\nbigf})$ be a complex of sheaves.
For each $\nbigf^j$,
we obtain the \v{C}ech complex
$\check{C}^{\bullet}(\nbigf^j)$ whose differential is denoted by
$d_{\check{C}}$.
(See {\rm\cite[\S B]{Huybrechts-text-book}}
for the signature of \v{C}ech complex.)
We obtain a morphism of complexes
$d_{\nbigf}:\check{C}^{\bullet}(\nbigf^j)
\to \check{C}^{\bullet}(\nbigf^{j+1})$
satisfying $d_{\nbigf}\circ d_{\nbigf}=0$.
Thus, we obtain a double complex
$\check{C}^{\bullet}(\nbigf^{\bullet})$.
The differential of the total complex
$\Tot\check{C}^{\bullet}(\nbigf^{\bullet})$
is induced by
$d_{\check{C}}:
 \check{C}^{i}(\nbigf^j)
 \lrarr
 \check{C}^{i+1}(\nbigf^j)$
and
$(-1)^id_{\nbigf}:
 \check{C}^{i}(\nbigf^j)
 \lrarr
 \check{C}^{i}(\nbigf^{j+1})$.
\hfill\qed
\end{rem}

\subsubsection{Relative Higgs bundles and
the induced infinitesimal deformations}
\label{subsection;23.5.10.1}

Let $S$ be a complex manifold.
Let $\Omega^1_{S\times X/S}$
denote the pull back of $K_X$ to $S\times X$
by the projection.
Let $\nbige$ be a holomorphic vector bundle on $S\times X$.
A morphism
$\theta_{\nbige}:\nbige\to\nbige\otimes\Omega^1_{S\times X/S}$
is called a relative Higgs field of $\nbige$.
Such a pair $(\nbige,\theta_{\nbige})$ is called
a relative Higgs bundle.

Let $s_0\in S$.
Let $\iota_{s_0}:X\to S\times X$
denote the inclusion $X=\{s_0\}\times X\subset S\times X$.
Let $\nbigo_{S,s_0}$ denote the stalk of $\nbigo_S$ at $s_0$
which is a local ring.
A morphism of local rings
$\nbigo_{S,s_0}\to \cnum[\epsilon]/\epsilon^2$
is equivalent to a tangent vector $v\in T_{s_0}S$.
It induces a morphism of ringed spaces
$\iota_{v}:X^{[1]}\to S\times X$
such that
$\iota_{s_0}=\iota_{v}\circ\iota$.
We obtain an infinitesimal deformation
$\iota_{v}^{\ast}(\nbige,\theta_{\nbige})$
of $\iota_{s_0}^{\ast}(\nbige,\theta_{\nbige})$.

\subsubsection{Appendix: $C^{\infty}$-description of
the induced infinitesimal deformation}

For any complex submanifold $S'$ of $S$,
let $(\nbige_{S'},\theta_{\nbige,S'})$
denote the restriction of $(\nbige,\theta_{\nbige})$
to $S'\times X$.
Let $q_{S'}:S'\times X\to X$ denote the projection.
Let $\nbigc^{\infty}_{X}\nbigo_{S'}$ denote
the sheaf of $C^{\infty}$-functions $f$ on $S'\times X$
such that $\delbar_{S'}f=0$.

\begin{lem}
\label{lem;24.10.20.1}
There exist a neighbourhood $S_0$ of $s_0$ in $S$
and a $C^{\infty}$-isomorphism
$\varphi:
q_{S_0}^{\ast}(\nbige_{s_0})
 \simeq
 \nbige_{S_0}$
such that
$\varphi$ is holomorphic in the $S_0$-direction,
i.e.,
$\varphi$ induces a morphism of
$\nbigc^{\infty}_X\nbigo_{S'}$-modules.
\end{lem}
\pf
We set $n=\rank \nbige$.
Let $X=\bigcup_{i=0}^m U_i$ be an open covering of $X$
such that $U_i\cap U_j=\emptyset$ $(1\leq i\neq j\leq m)$,
i.e.,
$U_i$ $(1\leq i\leq m)$ has the non-empty intersection
with only $U_0$.
There exist a neighbourhood $S_1$ of $s_0$ in $S$
and $\vecv_i$ be frames of $\nbige_{|S_1\times U_i}$.
We obtain holomorphic maps
$A_{i}(s,x):S_1\times(U_0\cap U_i)\to \GL(n,\cnum)$
determined by $\vecv_{0|S_1\times(U_0\cap U_i)}
=\vecv_{i|S_1\times(U_0\cap U_i)}A_i$.
We may assume that $A_i$ extends to
a holomorphic map on a neighbourhood of
$S_1\times\overline{U_0\cap U_i}$.
Let $I_n\in\GL(n,\cnum)$ denote the identity matrix.
Let $\nbigu$ be a convex open neighbourhood of $I_n$
in $\GL(n,\cnum)$.
There exists a neighbourhood $S_2$ of $s_0$ in $S_1$
such that
$A_i(s_0,x)^{-1}A_i(s,x)\in \nbigu$
for any $x\in \overline{U_0\cap U_i}$.

Let $U_i'$ be a relatively compact neighbourhood of
$(X\setminus U_0)\cap U_i$ in $U_i$.
Let $U_i''$ be a relatively compact neighbourhood of
$\overline{U_i}'$ in $U_i$.
We assume that
$\overline{U_i}''\cap \overline{U}_j=\emptyset$
for any $j\neq 0,i$.
Let $\chi_i:X\to[0,1]$ be a $C^{\infty}$-function
such that
$\chi_i=0$ on $U_i'$ and $\chi_i=1$ on $X\setminus U_i''$.

Let $\varphi_{U_0}$ be the isomorphism
$q_{S_2}^{\ast}(\nbige)_{|S_2\times U_0}
\simeq
\nbige_{|S_2\times U_0}$
defined by
$\varphi_{U_0}(q_{S_2}^{\ast}(\vecv_{0|\{s_0\}\times U_0}))
=\vecv_{0|S_2\times U_0}$.
Let $\varphi_{U_i}$ be the $C^{\infty}$-isomorphism
$q_{S_2}^{\ast}(\nbige)_{|S_2\times U_i}
\simeq
\nbige_{|S_2\times U_i}$
defined by
\[
\varphi_{U_i}(q_{S_2}^{\ast}(\vecv_{i|\{s_0\}\times U_0}))
=\vecv_{i|S_2\times U_0}
\Bigl(
I_n+\chi_i\cdot\bigl(
A_i(s,x)A_i(s_0,x)^{-1}
-I_n
\bigr)
\Bigr).
\]
We obtain 
a $C^{\infty}$-isomorphism
$\varphi:q_{S_2}^{\ast}(\nbige_{s_0})\simeq \nbige_{S_2}$
from $\varphi_{U_i}$ $(1\leq i\leq m)$
and the restriction of $\varphi_0$
to $X\setminus\bigcup_{i=1}^m \overline{U}_i''$.
It is holomorphic in the $S_2$-direction.
\hfill\qed

\vspace{.1in}
We set
$(E,\theta):=(\nbige_{s_0},\theta_{\nbige,s_0})$.
We obtain a $C^{\infty}$-section $B^{0,1}$
of $q_{S_0}^{-1}(\End(E)\otimes\Omega^{0,1}_X)$
by
\[
\varphi^{\ast}(\delbar_{\nbige_{S_0}})
=\delbar_{q_{S_0}^{\ast}(E)}+B^{0,1}.
\]
We obtain a $C^{\infty}$-section $B^{1,0}$
of $q_{S_0}^{-1}(\End(E)\otimes\Omega^{1,0}_X)$
by 
\[
\varphi^{\ast}(\theta_{\nbige,S_0})
=q_{S_0}^{\ast}(\theta)
+B^{1,0}.
\]
Because $\varphi$ is holomorphic in the $S_0$-direction,
$B^{1,0}$ and $B^{0,1}$ are holomorphic in the $S_0$-direction.

For any tangent vector $v\in T_{s_0}S$,
we obtain the section $\vtilde$ of
$T_{s_0}S\times X\subset
\iota_0^{\ast}(T(S\times X))$.
By using a $C^{\infty}$-isomorphism
$\varphi:q_{S_0}^{\ast}(E)\simeq \nbige$,
we define the derivative
$\vtilde(B^{1,0})$ and $\vtilde(B^{0,1})$.
We obtain a section
$\tau=\vtilde(B^{1,0})+\vtilde(B^{0,1})$
of $A^1(\Def(E,\theta))$.

\begin{lem}
\label{lem;24.10.20.5}
$\tau$ is a closed $1$-form in $A^1(\Def(E,\theta))$,
and its cohomology class is
$\iota_{v}(\nbige,\theta_{\nbige})$.
\end{lem}
\pf
Because
$\delbar_{q_{S_0}^{\ast}(E)}(B^{1,0})
+\bigl[
q_{S_0}^{\ast}(\theta),B^{0,1}
\bigr]=0$,
we obtain
$(\delbar+\ad\theta)(\tau)=
\delbar_{E}(\vtilde(B^{1,0}))
+[\theta,\vtilde(B^{0,1})]=0$.

There exist the following naturally defined quasi-isomorphisms:
\[
 A^{\bullet}(\Def(E,\theta))
 \lrarr
 \Tot^{\bullet}
 \check{C}^{\bullet}\bigl(
 \nbiga^{\bullet}(\Def(E,\theta))
 \bigr)
 \llarr
 \Tot^{\bullet}
 \check{C}^{\bullet}\bigl(
 \Def(E,\theta)
 \bigr).
\]
Here, $\nbiga^{\bullet}(\Def(E,\theta))$
is the Dolbeault resolution of $\Def(E,\theta)$.
See Remark \ref{rem;23.4.19.1}
for the construction of
$\Tot^{\bullet}\check{C}^{\bullet}(\nbigf^{\bullet})$
from a complex of sheaves $\nbigf^{\bullet}$.

Let us compare
the cocycles
induced by $\tau$ and
$\iota_v(\nbige,\theta_{\nbige})$
in the complex
$\Tot^{\bullet}
 \check{C}^{\bullet}\bigl(
 \nbiga^{\bullet}(\Def(E,\theta))
 \bigr)$.
In the following argument,
we shall identify $\nbige$
and $q_{S_0}^{-1}(E)$ by $\varphi$,
and we omit to denote $\varphi$.
Let $X=\bigcup_{i=0}^m U_i$
be an open covering.
We may assume that there exist holomorphic frames
$\vecv_i$ of $\nbige_{|S_0\times U_i}$.
We obtain the matrix valued holomorphic function $D_{i,j}$
on $S_0\times(U_{i}\cap U_j)$
determined by
$\vecv_{j|S_0\times(U_i\cap U_j)}
=\vecv_{i|S_0\times(U_i\cap U_j)}D_{ij}$.
We also obtain the matrix valued holomorphic $1$-form
$\Theta_i$ on $S_0\times U_i$
determined by
$\theta_{\nbige}\vecv_i=\vecv_i\Theta_i$.
The restriction of $\vecv_i$
to $\{s_0\}\times U_i$
are denoted by $\vecv_{i|s_0}$.
Let $\kappa_i$ be the holomorphic section of
$\End(E)\otimes\Omega_{U_i}^1$
determined by
$\kappa_i\vecv_{i|s_0}=\vecv_{i|s_0}\vtilde(\Theta_i)$.
Let $\rho_{i,j}$ be the holomorphic section of $\End(E)_{U_i\cap U_j}$
determined by
$\rho_{i,j}\vecv_j=\vecv_i \vtilde(D_{i,j})$.
By definition,
$(\{\rho_{i,j}\},\{\kappa_i\})$
is the $1$-cocycle of 
$\Tot^{\bullet}
\check{C}^{\bullet}\bigl(
\Def(E,\theta)
\bigr)$
induced by $\iota_v(\nbige)$.
Let $M_i$ be the matrix valued $C^{\infty}$-function
determined by
$\vecv_i=q^{-1}_{S_0}(\vecv_{i|s_0})M_i$.
The entries of $M_i$ are contained in
$\nbigc^{\infty}_X\nbigo_{S_0}$.
Because
$\delbar_{\nbige}(q_{S_0}^{-1}(\vecv_{i|s_0}))
=q_{S_0}^{-1}(\vecv_{i|s_0})(-\delbar(M_i)M_i^{-1})$,
we obtain
\[
 B^{0,1}q^{-1}(\vecv_{i|s_0})
 = q^{-1}(\vecv_{i|s_0})
 (-\delbar(M_i)M_i^{-1}).
\]
Hence, we obtain
\begin{equation}
\label{eq;24.10.20.2}
 \vtilde(B^{0,1})q^{-1}(\vecv_{i|s_0})
 =q^{-1}(\vecv_{i|s_0})
 \bigl(
 -\delbar \vtilde(M_i)
 \bigr).
\end{equation}
Because
$q_{S_0}^{-1}(\theta)q_{S_0}^{\ast}(\vecv_{i|s_0})
=q_{S_0}^{\ast}(\vecv_{i|s_0}) M_i\Theta_i M_i^{-1}$,
we obtain
\[
 B^{1,0} q_{S_0}^{\ast}(\vecv_{i|s_0})
 =q_{S_0}^{\ast}(\vecv_{i|s_0})
 \bigl(
 M_i\Theta_iM_i^{-1}-q_{S_0}^{\ast}(\Theta_{i|s_0})
 \bigr).
\]
Here, $\Theta_{i|s_0}$ denotes the restriction of
$\Theta_i$ to $\{s_0\}\times U_i$.
We obtain
\begin{equation}
\label{eq;24.10.20.3}
 \vtilde(B^{1,0})
 q_{S_0}^{\ast}(\vecv_{i|s_0})
 =q_{S_0}^{\ast}(\vecv_{i|s_0})
 \bigl(
 \vtilde(\Theta_i)
-[\Theta_{i|s_0},\vtilde(M_i)] 
 \bigr).
\end{equation}
Let $D_{ij|s_0}$ denote the restriction of $D_{ij}$
to $\{s_0\}\times (U_i\cap U_j)$.
Because $M_iD_{ij}=q_{S_0}^{\ast}(D_{ij|s_0})M_j$,
we obtain
\begin{equation}
\label{eq;24.10.20.4}
 \vtilde(D_{ij})=
 D_{ij|s_0}\vtilde(M_j)
 -\vtilde(M_i)D_{ij|s_0}.
\end{equation}
Let $\sigma_i$ be the $C^{\infty}$-section of $\End(E)$
on $U_i$
determined by
$\sigma_i\vecv_i=\vecv\bigl(\vtilde(M_i)\bigr)$.
By (\ref{eq;24.10.20.2}), (\ref{eq;24.10.20.3})
and (\ref{eq;24.10.20.4}),
we obtain
\[
\kappa_i-
 (B^{1,0}+B^{0,1})_{|U_i}
=(\delbar+\ad\theta)\sigma_i,
 \quad
 \rho_{i,j}
=\bigl(
 \sigma_{j|U_i\cap U_j}-\sigma_{i|U_i\cap U_j}
 \bigr).
\]
Hence, the $1$-cocycles induced by $\tau$ and $\iota_v(\nbige)$
are cohomologous in 
the complex
$\Tot^{\bullet}
\check{C}^{\bullet}\bigl(
\nbiga^{\bullet}(\Def(E,\theta))
\bigr)$.
\hfill\qed

\begin{rem}
For a different $C^{\infty}$-isomorphism
in Lemma {\rm\ref{lem;24.10.20.1}},
we obtain another $1$-cocycle
$\tau'$ of $A^1(\Def(E,\theta))$.
The difference $\tau-\tau'$ is a $1$-coboundary.
It follows from Lemma {\rm\ref{lem;24.10.20.5}}.
We can also check it directly.
\hfill\qed
\end{rem}

\subsection{Harmonic $1$-forms}

\subsubsection{Harmonic $1$-forms of $E$}

Let $h$ be a harmonic metric of $(E,\theta)$.
Let $\nabla=\delbar_E+\del_{E,h}$
be the Chern connection of $(E,\theta)$.
Let $\theta^{\dagger}_h$
denote the adjoint of $\theta$ with respect to $h$.
We fix a conformal metric $g_X$ of $X$.
Let $(\delbar_E+\theta)^{\ast}_{h,g_X}:
A^{i+1}(E)\to A^{i}(E)$
denote the formal adjoint of $\delbar_E+\theta$
with respect to $h$ and $g_X$.

\begin{df}
$\tau\in A^1(E)$
is called a harmonic $1$-form
of $(E,\theta,h)$
if $(\delbar_E+\theta)\tau=0$
and 
$(\delbar_E+\theta)^{\ast}_{h,g_X}\tau=0$.
\hfill\qed
\end{df}

\begin{rem}
Let $\Delta_{h,g_X}=
(\delbar_E+\theta)_{h,g_X}^{\ast}(\delbar_E+\theta)
+(\delbar_E+\theta)(\delbar_E+\theta)_{h,g_X}^{\ast}$.
When $X$ is compact,
we have $\Delta_{h,g_X}\tau=0$
if and only if
$(\delbar_E+\theta)\tau
=(\delbar_E+\theta)_{h,g_X}^{\ast}\tau=0$.
\hfill\qed
\end{rem}

\begin{lem}
\label{lem;23.5.5.130}
$\tau$ is a harmonic $1$-form if and only if
$(\delbar_E+\theta)\tau=
(\del_{E,h}+\theta_h^{\dagger})\tau=0$. 
\end{lem}
\pf
Let $\Lambda_X$ denote the adjoint of
the multiplication of the K\"ahler form associated with $g_X$.
According to \cite[Lemma 3.1]{s1},
we have
\[
(\delbar_E+\theta)^{\ast}_{h,g_X}\tau
=-\sqrt{-1}[\Lambda_X,\del_{E,h}+\theta_h^{\dagger}]\tau
=-\sqrt{-1}\Lambda_X((\del_{E,h}+\theta_h^{\dagger})\tau).
\]
It is $0$ if and only if
$(\del_{E,h}+\theta_h^{\dagger})\tau=0$.
\hfill\qed

\vspace{.1in}
For any complex number $\lambda$,
we obtain
$\DDlambda=\delbar_E+\theta+\lambda(\del_{E,h}+\theta^{\dagger}_h):
A^i(E)\to A^{i+1}(E)$.
Let $\DD^{\lambda\ast}_{h,g_X}:A^{i+1}(E)\to A^i(E)$
denote the formal adjoint of $\DDlambda$
with respect to $h$ and $g_X$.

\begin{cor}
\label{cor;23.5.4.11}
For a harmonic $1$-form $\tau$,
we have $\DDlambda\tau=0$
and  $\DD^{\lambda\ast}_{h,g_X}\tau=0$
for any $\lambda$.
\hfill\qed
\end{cor}

Let $X^{\dagger}$ denote the conjugate of $X$.
We obtain the Higgs bundle
$(E,\del_{E,h},\theta^{\dagger}_h)$ on $X^{\dagger}$,
and $h$ is a harmonic metric of
$(E,\del_{E,h},\theta^{\dagger}_h)$.
\begin{cor}
\label{cor;23.5.5.131}
$\tau\in A^1(E)$ is a harmonic $1$-form for
$(E,\delbar_E,\theta,h)$ 
if and only if
it is a harmonic $1$-form for 
$(E,\del_{E,h},\theta^{\dagger}_h,h)$.
\hfill\qed
\end{cor}

The following proposition
is useful to construct harmonic $1$-forms.
\begin{prop}
\label{prop;23.4.11.1}
Let $s$ be a section of $E$
such that $(\delbar_E+\theta)s=0$.
Then, $\tau=(\del_{E,h}+\theta_h^{\dagger})s$
is a harmonic $1$-form of
$(E,\theta,h)$.
Similarly,
if $s'$ is a section of $E$
such that $(\del_{E,h}+\theta^{\dagger}_h)s'=0$,
then $(\delbar_E+\theta)s'$ is a harmonic $1$-form of 
$(E,\theta,h)$. 
\end{prop}
\pf
We have
$(\del_{E,h}+\theta_h^{\dagger})\tau
=(\del_{E,h}+\theta_h^{\dagger})^2s=0$.
We also have
\[
(\delbar_E+\theta)\tau
=(\delbar_E+\theta)(\del_{E,h}+\theta_h^{\dagger})s
=-(\del_{E,h}+\theta_h^{\dagger})(\delbar_E+\theta)s
=0.
\]
Thus, we are done.
\hfill\qed

\vspace{.1in}

The following proposition is useful
in the study of pairings of harmonic $1$-forms.
\begin{prop}
\label{prop;23.5.4.20}
Let $\tau_i$ $(i=1,2)$ be harmonic $1$-forms
of $(E,\theta,h)$ on $X$.
Let $U$ be a simply connected relatively compact open subset of $X$
such that $\del U$ is smooth.
Note that,
because $\DD^1\tau_1=0$,
there exists a section $\sigma_1$ of $E$ on $U$
such that $\tau_{1|U}=\DD^1\sigma_1$. 
Then, we obtain
\begin{equation}
\label{eq;23.5.17.10}
\int_U
 \Bigl(
 \sqrt{-1}
 h(\tau_1^{1,0},\tau_2^{1,0})
-\sqrt{-1}
 h(\tau_1^{0,1},\tau_2^{0,1})
 \Bigr)
 =\sqrt{-1}\int_{\del U}
 \bigl(
 h(\sigma_1,\tau_2^{1,0})
-h(\sigma_1,\tau_2^{0,1}) 
 \bigr).
\end{equation}
\end{prop}
\pf
It follows from $\DD_{h,g_X}^{1\ast}\tau_2=0$.
\hfill\qed

\vspace{.1in}
Note that the left hand side of (\ref{eq;23.5.17.10}) is
the Hermitian product of
the space of $L^2$-sections of
$\nbiga^1(\Def(E,\theta))$ on $U$.

\begin{cor}
\label{cor;23.5.5.30}
In Proposition {\rm\ref{prop;23.5.4.20}},
we moreover assume that
$\tau_1=\tau_2$
and that there exist positive constants
$B_i$ $(i=1,2)$
such that
$|\tau_1|_{h,g_X}\leq B_1$
and
$|\sigma_1|_{h}\leq B_2$
on a neighbourhood of $\del U$.
Then,
there exists $B_3>0$,
depending only on $B_i$ $(i=1,2)$
such that
\[
 \|\tau_{1|U}\|_{L^2,h}
 \leq
 B_3.
\]
Here, $\|\tau_{1|U}\|_{L^2,h}$
denotes the $L^2$-norm of 
$\tau_{1|U}$ with respect to $h$ and $g_X$.
Note that it is independent of the choice of
a conformal metric $g_X$.
\hfill\qed
\end{cor}

\begin{rem}
Proposition {\rm\ref{prop;23.5.4.20}}
may be related with
the exactness in 
{\rm\cite[Lemma 12]{Dumas-Neitzke}}.
\hfill\qed
\end{rem}

\subsubsection{Harmonic $1$-forms of $\End(E)$}

Suppose that $(E,\theta)$
is generically regular semisimple
and that $\Sigma_{\theta}$ is smooth.
Let $C$ be a non-degenerate symmetric pairing of $(E,\theta)$.
Let $h$ be a harmonic metric of $(E,\theta)$
compatible with $C$.
The Chern connection of $\End(E)$
is also denoted by
$\nabla_h=\del_{E,h}+\delbar_E$.

Let $s$ be a holomorphic section of $\nbigh^0(\Def(E,\theta))$.
We obtain the harmonic $1$-form
$\tau=(\del_{E,h}+\ad\theta_h^{\dagger})s$ of $(\End(E),\ad\theta,h)$
by Proposition \ref{prop;23.4.11.1}.
There exists a holomorphic function $\alpha$
on $\Sigma_{\theta}$
such that $s=F_{\alpha}$
by Proposition \ref{prop;23.3.23.10}.

\begin{prop}\mbox{{}}
\label{prop;23.5.4.10}
\begin{itemize}
\item $\tau\in A^{1}(\Def(E,\theta;C))$.
      Let $\rho(\tau)$ denote the induced section of
      $\nbigh^1(\Def(E,\theta;C))\simeq\pi_{\ast}\Omega^1_{\Sigma_{\theta}}$.
 \item $\rho(\tau)=\pi_{\ast}d\alpha$.
 \end{itemize}
\end{prop}
\pf
Let $P\in X^{\circ}$.
Let $(X_P,z_P)$ be
a simply connected coordinate neighbourhood
around $P$
such that $X_P\cap D(\theta)=\emptyset$.
There exists a decomposition
$(E,\theta)_{|X_P}
=\bigoplus_{i=1}^{\rank(E)}
(E_{P,i},\theta_{P,i})$
into Higgs bundles of rank $1$.
For each $i$, let $v_{P,i}$ be a holomorphic frame of $E_{P,i}$
such that $C(v_{P,i},v_{P,i})=1$.
We obtain a frame $\vecv_P=(v_{P,i})$.
Let $B$ be the matrix 
determined by
$\del_{E,h}\vecv=\vecv\cdot B\,dz$.
Because $\nabla_h=\delbar_E+\del_{E,h}$
is compatible with $C$,
we obtain that the matrix $B$ is anti-symmetric.
Then, it is easy to check the claims.
\hfill\qed

\begin{rem}
Let $X^{\dagger}$ denote the conjugate of $X$,
and let $E^{\dagger}_h$ denote the holomorphic
vector bundle $(E,\del_{E,h})$ on $X^{\dagger}$.
Then, $\theta^{\dagger}_h$
is a Higgs field of $E^{\dagger}_h$,
and $h$ is a harmonic metric of 
$(E^{\dagger}_h,\theta^{\dagger}_h)$.
Because
$\del_{E,h}(s^{\dagger}_h)=(\delbar_Es)^{\dagger}_h=0$
and 
$[\theta^{\dagger}_h,s^{\dagger}_h]=-[\theta,s]^{\dagger}_h=0$, 
$s^{\dagger}_h$ is a holomorphic section of
$\nbigh^0(E^{\dagger},\theta^{\dagger}_h)$.
Hence,
$(\delbar_E+\ad\theta^{\dagger}_h)s^{\dagger}_h$
is also a harmonic $1$-form
of $(\End(E),\ad\theta,h)$.
\hfill\qed
\end{rem}

\subsection{Some estimates for harmonic bundles}

\subsubsection{Uhlenbeck theorem}

For $r>0$,
we set $\nbigb(r)=\{|z|<r\}$.
For a matrix valued $C^{\infty}$-function $s$ on $\nbigb(r)$,
and for $0<r_1\leq r$,
let $\|s_{|\nbigb(r_1)}\|_{L_{\ell}^p,g_0}$
denote the $L_{\ell}^p$-norm of $s_{|\nbigb(r_1)}$
with respect to the Euclidean metric $g_0$
and an appropriate norm of the space of matrices.
If $r_1=r$,
$\|s_{|\nbigb(r)}\|_{L_{\ell}^p,g_0}$ is denoted by
$\|s\|_{L_{\ell}^p,g_0}$.

Let $V$ be a holomorphic vector bundle on $\nbigb(1)$
equipped with a Hermitian metric $h$.
Let $\nabla_h$ denote the Chern connection,
and let $R(h)=\gbigr\,dz\,d\zbar$ denote the curvature.
For a $C^{\infty}$-section $s$ of $V$ or $\End(V)$
and for $0<r_1<1$,
let $\|s_{|\nbigb(r_1)}\|_{L_{\ell}^p,h,g_0}$
denote the $L_{\ell}^p$-norm of $s_{|\nbigb(r_1)}$
with respect to
$\nabla_h$, $h$ and $g_0$.
For a frame $\vecv$ of $V$,
let $\nbiga^{\vecv}$ denote
the connection form of $\nabla_h$ with respect to $\vecv$.
If $r_1=1$,
$\|s_{|\nbigb(1)}\|_{L_{\ell}^p,h,g_0}$
is denoted by
$\|s\|_{L_{\ell}^p,h,g_0}$.
We recall the following theorem of Uhlenbeck.
(See \cite[Theorem 1.3]{u1} for more precise statement.)

\begin{thm}[Uhlenbeck]
\label{thm;23.5.5.20}
Let $p>2$.
There exist positive constants $\epsilon_0$ and $C_{0}$,
which are independent of $(V,h)$,
such that the following holds:
\begin{itemize}
 \item Suppose that $|\gbigr|_{h}\leq \epsilon_0$ on $\nbigb(1)$.
       Then, there exists an orthonormal frame $\vecv$
       of $V$ such that
       $\|\nbiga^{\vecv}\|_{L_1^p,g_0}
       \leq
       C_{0}\|\gbigr\|_{L^p,h,g_0}$
       and
       $d^{\ast}\nbiga^{\vecv}=0$.
       Here, $d^{\ast}$
       denotes the formal adjoint of
       the exterior derivative $d$.
\hfill\qed
\end{itemize} 
\end{thm}

\subsubsection{Universal boundedness of higher derivatives of
Higgs fields}

Let $(E,\theta,h)$ be
a harmonic bundle on $\nbigb(r_0)$ for some $r_0>1$.
Let $f$ be the endomorphism of $E$
determined by $\theta=f\,dz$.
Suppose that there exists $B_1>0$ such that the following holds.
\begin{itemize}
 \item For any $z\in \nbigb(r_0)$,
       any eigenvalue $\alpha$ of
       $f_{|z}$
       satisfies
       $|\alpha|\leq B_1$.
\end{itemize}

Let $\nabla_h$ denote the Chern connection of $(E,h)$.
Let $R(h)=\gbigr\,dz\,d\zbar$ denote the curvature of $\nabla_h$.
Let $g_0=dz\,d\zbar$ be the Euclidean metric.
For a section $\alpha$ of $\End(E)$,
and for $\ell\in\seisuu_{\geq 0}$ and $p\geq 1$,
let $\|\alpha_{|\nbigb(r)}\|_{L_{\ell}^p,h,g_0}$
denote the $L_{\ell}^p$-norm
of $\alpha$ on $\nbigb(r)$
with respect to $h$ and $g_0$
and the derivatives induced by $\nabla_{h}$.

Let $0<r_1<r_0$.
We set $r_2=\frac{r_0+r_1}{2}$.
By a variant of Simpson's main estimate
(see \cite[Proposition 2.1]{Decouple}),
we obtain the following estimate
in \cite[Proposition 4.1]{Mochizuki-Szabo}.
\begin{prop}
\label{prop;23.5.5.60}
There exists $B_2>0$ such that
\[
 \sup|f_{|\nbigb(r_2)}|_{h}
+\sup|f^{\dagger}_{h|\nbigb(r_2)}|_{h}
+\sup|\gbigr_{|\nbigb(r_2)}|_{h}\leq B_{2}.
\] 
Moreover, for any $\ell\in\seisuu_{\geq 0}$
and $p\geq 1$,
there exists $B_{3}(\ell,p)>0$,
depending only on $r_0-r_1$, $B_1$ and $(\ell,p)$
such that
\[
 \|f_{|\nbigb(r_2)}\|_{L_{\ell}^p,h,g_0}
+\|f^{\dagger}_{h|\nbigb(r_2)}\|_{L_{\ell}^p,h,g_0}
+\|\gbigr_{|\nbigb(r_2)}\|_{L_{\ell}^p,h,g_0}\leq B_{3}(\ell,p).
\]
\hfill\qed
\end{prop}

\subsubsection{The sup norm and the $L^2$-norm of harmonic $1$-forms}

Let $\tau\in A^1(E)$ be a harmonic $1$-form of
$(E,\theta,h)$ on $\nbigb(r_0)$.
Let $\|\tau\|_{L^2,h}$ be the $L^2$-norm of
$\tau$ with respect to $h$ and $g_0$.
Note that $\|\tau\|_{L^2,h}$ is independent of the choice of
a conformal metric of $\nbigb(r_0)$.
\begin{prop}
\label{prop;23.5.5.31}
There exists $B_4>0$,
depending only on $r_0-r_1$ and $B_1$,
 such that
the following holds: 
\[
 \sup_{\nbigb(r_1)}|\tau|_{h,g_0}
\leq B_4\|\tau\|_{L^2,h}. 
\] 
\end{prop}
\pf
Let $\epsilon_0$ be as in Theorem \ref{thm;23.5.5.20}.
Let $B_2$ be as in Proposition \ref{prop;23.5.5.60}.
We set
$\epsilon_1:=\min\{2^{-1}(r_2-r_1),B_2^{-1/2}\epsilon_0^{1/2}\}$,
which satisfies
$0<\epsilon_1<r_2-r_1$
and $B_2\epsilon_1^2\leq\epsilon_0$.
For any $z_0\in\nbigb(r_1)$,
let $\psi_{z_0}:\nbigb(1)\lrarr\nbigb(r_2)$
be defined by $\psi_{z_0}(\zeta)=z_0+\epsilon_1\zeta$.
We obtain
$(\Ehat_{z_0},\thetahat_{z_0},\hhat_{z_0})
=\psi_{z_0}^{\ast}(E,\theta,h)$
and
$\tauhat_{z_0}=\psi_{z_0}^{\ast}(\tau)$ on $\nbigb(1)$.
Let $\nabla_{\hhat_{z_0}}$ denote the Chern connection of
$(\Ehat_{z_0},\hhat_{z_0})$,
and let $R(\hhat_{z_0})$ denote the curvature of
$\nabla_{\hhat_{z_0}}$.
Let 
$\|\tauhat_{z_0}\|_{L^2,\hhat_{z_0}}$
denote the $L^2$-norm of $\tauhat_{z_0}$
with respect to $\hhat_{z_0}$ and $g_0$.
We obtain the following on $\nbigb(1)$:
\[
 |\thetahat_{z_0}|_{\hhat_{z_0}}
 \leq
 B_2\epsilon_1,
 \quad
 |R(\hhat_{z_0})|_{\hhat_{z_0}}
 \leq\epsilon_0,
 \quad
 \|\tauhat_{z_0}\|_{L^2,\hhat_{z_0}}
 \leq \|\tau\|_{L^2,h}.
\]
Let $\vecv$ be an orthonormal frame of $(\Ehat_{z_0},\hhat_{z_0})$
as in Theorem \ref{thm;23.5.5.20}.
Let $\nbiga$ be the connection form of
$\nabla_{\hhat_{z_0}}$
with respect to $\vecv$.
The $(p,q)$-part is denoted by $\nbiga^{p,q}$.
Let $\Theta$
denote the matrix valued $(1,0)$-form
determined by
$\thetahat_{z_0}\vecv=\vecv\Theta$.
We have the expression $\tau=\sum \tau_iv_i$,
and we obtain the vector-valued
$1$-form $\vectau_{z_0}=(\tau_i)$.
The $(p,q)$-part of $\vectau_{z_0}$
is denoted by $\vectau_{z_0}^{p,q}$.
Because $\tau$ is harmonic,
we have
\begin{equation}
\label{eq;26.1.27.10}
 \delbar\vectau_{z_0}^{1,0}
 +\nbiga^{0,1}(\vectau_{z_0}^{1,0})
 +\Theta(\vectau_{z_0}^{0,1})=0,
 \quad
  \del\vectau_{z_0}^{0,1}
 +\nbiga^{1,0}(\vectau_{z_0}^{0,1})
 +\lefttop{t}\overline{\Theta}(\vectau_{z_0}^{1,0})=0.
\end{equation}
By (\ref{eq;26.1.27.10}),
the $L_1^2$-norm of
$\vectau_{z_0}$ on $\nbigb(2/3)$
is dominated by
$\|\tau\|_{L^2,h}$.
It implies that for $p>2$,
the $L^p$-norm of $\vectau_{z_0}$ on $\nbigb(1/2)$
is dominated by
$\|\tau\|_{L^2,h}$.
By (\ref{eq;26.1.27.10}) again,
the $L_1^p$-norm of 
$\vectau_{z_0}$
on $\nbigb(1/3)$
is dominated by
$\|\tau\|_{L^2,h}$.
Hence, we can show that there exists $B_{5}>0$,
depending only on $B_2$,
such that 
$|\tauhat_{z_0}|_{\hhat_{z_0},g_0}\leq B_5\|\tau\|_{L^2,h}$
on $\{|\zeta|<1/4\}$.
It implies that
$|\tau|_{h,g_0}\leq \epsilon_1^{-1}B_5\|\tau\|_{L^2,h}$
on $\nbigb(r_1)$.
\hfill\qed

\subsubsection{Poisson equation}

Let $s$ be a $C^{\infty}$-section of $E$ on $\nbigb(r_0)$.

\begin{prop}
\label{prop;23.5.5.61}
There exists $B_{10}>0$
depending only on $B_1$ and $r_0-r_1$
such that
the following holds:
\begin{equation}
\label{eq;23.5.5.65}
 \sup_{\nbigb(r_1)}\bigl|s\bigr|_{h}
 +\sup_{\nbigb(r_1)}
 \bigl|
 (\delbar_E+\theta)s
 \bigr|_{h,g_0}
 \leq
 B_{10}
 \|s\|_{L^2,h,g_0}
+B_{10}
 \sup_{\nbigb(r_2)}\bigl|
 (\delbar_E+\theta)^{\ast}_{h,g_0}
 (\delbar_E+\theta)s
 \bigr|_{h}.
\end{equation}
\end{prop}
\pf
We set $\sigma=(\delbar_E+\theta)s$
and $\rho=(\del_{E,h}+\theta^{\dagger})(\delbar_E+\theta)s$.
We use the notation in the proof of Proposition \ref{prop;23.5.5.31}.
We obtain the section
$s_{z_0}=\psi_{z_0}^{-1}(s)$,
the $1$-form 
$\sigma_{z_0}=\psi_{z_0}^{-1}(\sigma)$
and the $2$-form
$\rho_{z_0}=\psi_{z_0}^{-1}(\sigma)$
on $\nbigb(1)$.
Let $\kappa:\real\to \{0\leq a\leq 1\}$
be a $C^{\infty}$-function
such that
(i) $\kappa(t)=1$ $(t\leq 3/4)$ and $\kappa(t)=0$ $(t\geq 4/5)$,
(ii) $\kappa'/\kappa^{1/2}$ is bounded on $\{t\,|\,\kappa(t)>0\}$.
Let $\chi(\zeta)=\kappa(|\zeta|)$.
We obtain
\begin{multline}
 \int_{\nbigb(1)} \chi|\sigma_{z_0}|^2_{\hhat_{z_0},g_0}
 \leq
 \int_{\nbigb(1)} |\delbar\chi|_{g_0}\cdot
 |s_{z_0}|_{\hhat_{z_0}}
 \cdot
 |\sigma_{z_0}|_{\hhat_{z_0},g_0}
 +
 \int_{\nbigb(1)}
 \chi |s_{z_0}|_{\hhat_{z_0}}
 \cdot |\rho_{z_0}|_{\hhat_{z_0},g_0}
\\
 \leq
 \left(
 \int_{\nbigb(1)} |\chi^{-1/2}\delbar\chi|^2|s_{z_0}|^2_{\hhat_{z_0}}
 \right)^{1/2}
 \left(
 \int_{\nbigb(1)} \chi|\sigma_{z_0}|^2_{\hhat_{z_0},g_0}
 \right)^{1/2}
 +
 \left(
 \int_{\nbigb(1)}
 |s_{z_0}|_{\hhat_{z_0}}^2
 \right)^{1/2}
 \cdot
 \left(
 \int_{\nbigb(1)}
 |\rho_{z_0}|_{\hhat_{z_0},g_0}^2
 \right)^{1/2}.
\end{multline}
There exist $C_1,C_2>0$, depending only on $\kappa$,
such that
\begin{equation}
\label{eq;26.1.27.4}
 \left(
 \int_{\nbigb(3/4)}|\sigma_{z_0}|^2_{\hhat_{z_0},g_0}
 \right)^{1/2}
\leq
 \left(
 \int_{\nbigb(1)} \chi|\sigma_{z_0}|^2_{\hhat_{z_0},g_0}
 \right)^{1/2}
 \leq
 C_1
 \left(
 \int_{\nbigb(1)} |s_{z_0}|^2_{\hhat_{z_0}}
 \right)^{1/2}
 + C_2
  \left(
 \int_{\nbigb(1)}|\rho_{z_0}|^2_{\hhat_{z_0},g_0}
 \right)^{1/2}. 
\end{equation}

We have the expressions
$s_{z_0}=\sum s_{z_0,i}v_i$,
$\sigma_{z_0}=\sum \sigma_{z_0,i}v_i$
and
$\rho_{z_0}=\sum \rho_{z_0,i}v_i$.
We set
$\vecs=(s_i)$,
$\vecsigma_{z_0}=(\sigma_{z_0,i})$,
and
$\vecrho_{z_0}=(\rho_{z_0,i})$.
We have
\begin{equation}
\label{eq;26.1.27.2}
 (\delbar+\nbiga^{0,1})\vecs_{z_0}
 +\Theta\vecs=\vecsigma_{z_0},
\end{equation}
\begin{equation}
\label{eq;26.1.27.3}
 (\del+\nbiga^{1,0})\vecsigma_{z_0}^{0,1}
 +\lefttop{t}\overline{\Theta}(\vecsigma^{1,0}_{z_0})
 =\vecrho_{z_0},
 \quad
(\delbar+\nbiga^{0,1})\vecsigma_{z_0}^{1,0}
+\Theta(\vecsigma_{z_0}^{0,1})=0.
\end{equation}
By (\ref{eq;26.1.27.4}), (\ref{eq;26.1.27.3})
and the elliptic regularity,
the $L_1^2$-norm of $\vecrho_{z_0}$ on $\nbigb(2/3)$
is dominated by
$\|s\|_{L^2,h,g_0}$
and
$\sup|\rho_{z_0}|_{\hhat_0,g_0}$.
It implies that for $p>2$,
the $L^p$-norm of
$\vecrho_{z_0}$ on $\nbigb(1/2)$
is dominated by
$\|s\|_{L^2,h,g_0}$
and
$\sup|\rho_{z_0}|_{\hhat_0,g_0}$.
By (\ref{eq;26.1.27.4}), (\ref{eq;26.1.27.3})
and the elliptic regularity again,
the $L_1^p$-norm of $\vecrho_{z_0}$ on $\nbigb(1/3)$
is dominated by
$\|s\|_{L^2,h,g_0}$
and
$\sup|\rho_{z_0}|_{\hhat_0,g_0}$.
Hence,
$|\sigma_{z_0}|_{\hhat_0,g_0}$ on $\nbigb(1/4)$
is dominated by
$\|s\|_{L^2,h,g_0}$
and
$\sup|\rho_{z_0}|_{\hhat_0,g_0}$.
By (\ref{eq;26.1.27.2}) and the elliptic regularity,
$|s_{z_0}|$ on $\nbigb(1/5)$
is dominated by
$\|s\|_{L^2,h,g_0}$
and
$\sup|\rho_{z_0}|_{\hhat_0,g_0}$.
Then, we obtain (\ref{eq;23.5.5.65}).
\hfill\qed

\subsection{Real bilinear forms on complex vector spaces}

This is a preliminary for \S\ref{subsection;23.4.24.10}.

\subsubsection{Bilinear forms on real vector spaces}
\label{subsection;23.5.16.10}

Let $V_{\real}$ be a finite dimensional $\real$-vector space.
We set $V_{\real}^{\lor}=\Hom_{\real}(V_{\real},\real)$.
Let $g_{\real}$ be
a non-degenerate symmetric bilinear form on $V_{\real}$.
It induces the isomorphism $\Psi_{g_{\real}}:V_{\real}\to V_{\real}^{\lor}$
by $\Psi_{g_{\real}}(u)(v)=g_{\real}(u,v)$.
We obtain the induced symmetric form $g_{\real}^{\lor}$
on $V_{\real}^{\lor}$ by
\[
 g_{\real}^{\lor}(a,b)
=g_{\real}\bigl(
 \Psi_{g_{\real}}^{-1}(a),
 \Psi_{g_{\real}}^{-1}(b)
 \bigr).
\]
\begin{lem}
We have $\Psi_{g_{\real}}^{\lor}\circ\Psi_{g_{\real}}=\id_{V_{\real}}$
and $(g_{\real}^{\lor})^{\lor}=g_{\real}$.
If $g_{\real}$ is positive definite,
$g_{\real}^{\lor}$ is also positive definite.
\hfill\qed
\end{lem}

Let $\omega$ be a symplectic form on $V_{\real}$.
It induces the isomorphism
$\Psi_{\omega}:V_{\real}\to V_{\real}^{\lor}$
by $\Psi_{\omega}(u)(v)=\omega(u,v)$.
We obtain the induced symplectic form $\omega^{\lor}$
on $V_{\real}^{\lor}$ by 
\[
 \omega^{\lor}(a,b)
 =\omega\bigl(
\Psi_{\omega}^{-1}(a),
\Psi_{\omega}^{-1}(b)
 \bigr).
\]
\begin{lem}
We have $\Psi_{\omega^{\lor}}\circ\Psi_{\omega}=-\id_V$
and $(\omega^{\lor})^{\lor}=\omega$.
\hfill\qed
\end{lem}

\subsubsection{Real bilinear forms on complex vector spaces}
\label{subsection;23.5.16.11}

Let $V$ be a finite dimensional $\cnum$-vector space.
Let $V_{\real}$ denote the underlying $\real$-vector space.
We set $V^{\lor}=\Hom_{\cnum}(V,\cnum)$
and $V_{\real}^{\lor}=\Hom_{\real}(V_{\real},\real)$.
The multiplication of $\sqrt{-1}$ on $V$
induces the $\real$-automorphism $I$ on $V_{\real}$
such that $I^2=-\id_{V_{\real}}$.
It induces an $\real$-linear automorphism $I^{\lor}$
on $V_{\real}^{\lor}$ by $I^{\lor}(f)(v)=f(Iv)$.

Let us recall the identification of
$V_{\real}^{\lor}$ and $V^{\lor}$.
There exists the decomposition
\[
 V_{\real}^{\lor}
 \otimes\cnum
 =\Hom_{\real}(V_{\real},\cnum)
 =V^{\lor}
 \oplus
 \overline{V^{\lor}},
\]
where $V^{\lor}$ is the $\sqrt{-1}$-eigen space of $I^{\lor}$,
and $\overline{V^{\lor}}$ is the $-\sqrt{-1}$-eigen space of
$I^{\lor}$.
We obtain the isomorphism of $\real$-vector spaces
$\Phi:V_{\real}^{\lor}\to V^{\lor}$
as the composition of
the inclusion
$V_{\real}^{\lor}\to V_{\real}^{\lor}\otimes\cnum$
and the projection
$V_{\real}^{\lor}\otimes\cnum
=V^{\lor}\oplus \overline{V^{\lor}}\to V^{\lor}$.
It is explicitly described as
\[
 \Phi(f)=\frac{1}{2}\bigl(f-\sqrt{-1}I^{\lor}(f)\bigr).
\]
\begin{rem}
If we regard $V_{\real}^{\lor}$ with $I^{\lor}$
as a $\cnum$-vector space,
then $\Phi$ is an isomorphism of $\cnum$-vector spaces.
\hfill\qed
\end{rem}

We can check the following lemmas by direct computations.
\begin{lem}
\mbox{{}}\label{lem;23.5.16.13}
Let $g_{\real}$ be a non-degenerate symmetric bilinear form on $V_{\real}$.
Let $\omega$ be a symplectic form on $V_{\real}$.
\begin{itemize}
 \item 
If $I$ is an isometry with respect to $g_{\real}$,
we have $\Psi_{g_{\real}}\circ I=-I^{\lor}\circ\Psi_{g_{\real}}$.
If $I$ is an isometry with respect to $\omega$,
we have $\Psi_{\omega}\circ I=-I^{\lor}\circ\Psi_{\omega}$.
\item
     Suppose that $I$ is an isometry with respect to $g_{\real}$,
     and that $\omega(u,v)=g_{\real}(u,Iv)$ for any $u,v$.
Then, we obtain
$\omega^{\lor}(u^{\lor},v^{\lor})
=-g_{\real}^{\lor}(u^{\lor},I^{\lor}v^{\lor})$.
In particular,
if $\omega(I\cdot,\cdot)$ is positive definite,
$-\omega^{\lor}(I^{\lor}\cdot,\cdot)$ is positive definite.
     \hfill\qed
 \end{itemize}
\end{lem}

\subsection{Bilinear forms on the first cohomology group
of a compact Riemann surface}
\label{subsection;23.4.24.10}

Let $\Sigma$ be a compact Riemann surface.
There exists a natural symplectic form $\omega_{\Sigma}$
on $H^1(\Sigma,\real)$
given by
\[
 \omega_{\Sigma}(a,b)
 =\int_{\Sigma}a\wedge b.
\]
There exists the Hodge decomposition:
\[
 H^1(\Sigma,\real)
 \otimes\cnum
 =H^1(\Sigma,\nbigo_{\Sigma})
 \oplus
 H^0(\Sigma,K_{\Sigma}).
\]
The projections induce isomorphisms of $\real$-vector spaces
\[
 \ttF_1:
 H^1(\Sigma,\real)
 \simeq
 H^1(\Sigma,\nbigo_{\Sigma}),
 \quad\quad
 \ttF_2:
 H^1(\Sigma,\real)
 \simeq
 H^0(\Sigma,K_{\Sigma}). 
\]

\subsubsection{Bilinear forms on $H^1(\Sigma,\nbigo_{\Sigma})$}
Let $\omega_1$ be the real symplectic form
on $H^1(\Sigma,\nbigo_{\Sigma})$
given by
\[
 \omega_1(a_1,a_2)
=-\omega_{\Sigma}(\ttF_1^{-1}(a_1),\ttF_1^{-1}(a_2)).
\]
There exist the harmonic representatives $\eta_i$ of $a_i$,
i.e.,
$\eta_i$ are the $(0,1)$-forms such that
$\del_{\Sigma}\eta_i=0$,
and $[\eta_i]=a_i$ in $H^1(\Sigma,\nbigo_{\Sigma})$.
We have
\[
 \ttF_1^{-1}([\eta_i])
=[\eta_i+\overline{\eta_i}].
\]
Therefore,
\[
 \omega_1([\eta_1],[\eta_2])
 =-\int_{\Sigma}
 (\eta_1+\overline{\eta_1})
 \wedge
 (\eta_2+\overline{\eta_2})
 =-\int_{\Sigma}
 \bigl(
 \eta_1\overline{\eta_2}
+\overline{\eta_1}\eta_2
 \bigr).
\]
The associated positive definite symmetric bilinear form
$g_{1\real}$ is given as follows:
\[
 g_{1\real}([\eta_1],[\eta_2])
=\omega_1([\sqrt{-1}\eta_1],[\eta_2])
=-\sqrt{-1}\int_{\Sigma}
(\eta_1\etabar_2-\overline{\eta_1}\eta_2)
=-\sqrt{-1}\int_{\Sigma}
(\eta_1\etabar_2+\eta_2\overline{\eta_1}).
\]
The associated Hermitian metric
$g_1=g_{1\real}+\sqrt{-1}\omega_1$ is given as follows:
\[
 g_1([\eta_1],[\eta_2])
=-2\sqrt{-1}\int_{\Sigma}
 \eta_1\etabar_2.
\]

\subsubsection{Bilinear forms on $H^0(\Sigma,K_{\Sigma})$}

Let $\omega_2$ be the real symplectic form
on $H^0(\Sigma,K_{\Sigma})$
given by
\[
 \omega_2(a_1,a_2)
=\omega_{\Sigma}(\ttF_2^{-1}(a_1),\ttF_2^{-1}(a_2)).
\]
We may naturally regard $a_i$
as holomorphic $1$-forms on $\Sigma$.
We have
\[
 \ttF_2^{-1}(a_i)
=a_i+\overline{a_i}.
\]
Therefore,
\[
 \omega_2(a_1,a_2)
 =\int_{\Sigma}
 (a_1+\overline{a_1})
 \wedge
 (a_2+\overline{a_2})
 =\int_{\Sigma}
 \bigl(
 a_1\overline{a_2}
+\overline{a_1}a_2
 \bigr).
\]
The associated positive definite symmetric bilinear form
$g_{2\real}$ is given as follows.
\[
 g_{2\real}(a_1,a_2)
=\omega_2(\sqrt{-1}a_1,a_2)
=\sqrt{-1}\int_{\Sigma}
 (a_1 \overline{a_2}-\overline{a_1}a_2)
=\sqrt{-1}\int_{\Sigma}
(a_1\abar_2+a_2\overline{a_1}).
\]
The associated hermitian metric
$g_{2}$ is given as follows:
\[
 g_2(a_1,a_2)
=2\sqrt{-1}\int_{\Sigma}a_1\abar_2.
\]

\subsubsection{Duality}

By the Serre duality,
the following pairing is perfect.
\begin{equation}
\label{eq;23.5.3.20}
 H^1(\Sigma,\nbigo_{\Sigma})
 \otimes
 H^0(\Sigma,K_{\Sigma})
\lrarr \cnum,
\quad\quad
 (a,b)
 \longmapsto
 -\int_{\Sigma} a\wedge b.
\end{equation}
By this perfect pairing,
we may regard
$H^0(\Sigma,K_{\Sigma})
\simeq
H^1(\Sigma,\nbigo_{\Sigma})^{\lor}$.
We also have the isomorphism
$H^1(\Sigma,\nbigo_{\Sigma})^{\lor}
\simeq
H^1(\Sigma,\nbigo_{\Sigma})_{\real}^{\lor}$
as explained in \S\ref{subsection;23.5.16.11}.
By the composition,
we obtain the following isomorphism of $\real$-vector spaces:
\begin{equation}
\label{eq;23.5.16.12}
 H^0(\Sigma,K_{\Sigma})_{\real}
 \simeq
 H^1(\Sigma,\nbigo_{\Sigma})_{\real}^{\lor}.
\end{equation}

We have the bilinear forms $g_{1\real}^{\lor}$
and $\omega_1^{\lor}$
on $H^1(\Sigma,\nbigo_{\Sigma})_{\real}^{\lor}$
as explained in \S\ref{subsection;23.5.16.10}.

\begin{lem}
\label{lem;23.5.3.21}
Under the isomorphism {\rm(\ref{eq;23.5.16.12})},
we have $\omega_2=-\omega_1^{\lor}$
and $g_{2\real}=g_{1\real}^{\lor}$.
\end{lem}
\pf
Let $\Psitilde_{\omega_1}:H^1(\Sigma,\nbigo_{\Sigma})_{\real}
\simeq
H^0(\Sigma,K_{\Sigma})_{\real}$
denote the composite of $\Psi_{\omega_1}$
and the inverse of (\ref{eq;23.5.16.12}).
(See \S\ref{subsection;23.5.16.10} for $\Psi_{\omega_1}$.)
For harmonic $(0,1)$-forms $\eta_i$ $(i=1,2)$,
the following holds:
\[
 \frac{1}{2}\Bigl(
 \Psi_{\omega_1}(\eta_1)
 -\sqrt{-1}I^{\lor}\Psi_{\omega_1}(\eta_1)
 \Bigr)(\eta_2)
=-\int_{\Sigma}\etabar_1\eta_2.
\] 
We obtain
$\Psitilde_{\omega_1}(\eta_1)=-\etabar_1$.

Let $a_i\in H^0(\Sigma,K_{\Sigma})$ $(i=1,2)$.
We obtain
\[
-\omega_1^{\lor}(a_1,a_2)
=-\omega_1\bigl(\Psitilde_{\omega_1}^{-1}(a_1),
\Psitilde_{\omega_1}^{-1}(a_2)\bigr) 
=-\omega_1(\abar_1,\abar_2)
=\int_{\Sigma}
\bigl(
 \abar_1a_2+a_1\abar_2
 \bigr)
=\omega_2(a_1,a_2).
\]
We obtain $g_{2\real}=g_{1\real}^{\lor}$
from Lemma \ref{lem;23.5.16.13}.
\hfill\qed

\subsection{Estimate for the family of the associated Laplacian operators}
\label{subsection;23.5.3.30}

Let $X$ be a compact Riemann surface
with a conformal metric $g_X$.
Let $(E,\theta)$ be a generically regular semisimple
Higgs bundle of degree $0$
whose spectral curve is smooth.
We fix a flat metric $h_{\det(E)}$ of $\det(E)$.
For each $t>0$,
there exists a harmonic metric $h_t$
of $(E,t\theta)$.
We recall that
the sequence $h_t$ is convergent to
the limiting configuration $h_{\infty}$
in the $C^{\infty}$-sense locally on $X^{\circ}$.

Let $\End_0(E)$ denote the trace free part of
$\End(E)$.
Let $(\delbar_E+t\ad\theta)^{\ast}_{h_t}$
denote the formal adjoint of
$\delbar_E+t\ad\theta$ with respect to $h_t$ and $g_X$.
We consider the operators
$\Delta_{t}=
(\delbar_E+t\ad\theta)^{\ast}_{h_t}\circ(\delbar_E+t\ad\theta)$
on $A^0(\End_0(E))$.
For any sections $s$ of $A^0(\End_0(E))$,
let $\|s\|_{L^2,h_t}$
denote the $L^2$-norm of $s$
with respect to $h_t$ and $g_X$.

\begin{prop}
\label{prop;23.4.12.1}
There exists $C>0$
such that
 $\|\Delta_t(s)\|_{L^2,h_t}\geq C\|s\|_{L^2,h_t}$
for any $t\geq 1$ and $s\in A^0(\End_0(E))$. 
\end{prop}
\pf
If $s\in A^0(\End_0(E))$ satisfies $\Delta_ts=0$,
it is standard to obtain
$(\delbar_E+t\ad\theta)s=0$.
Because $(E,t\theta)$ is stable,
we obtain $s=0$.

Let $a_t>0$ denote the first eigenvalue of $\Delta_t$.
It is enough to prove that $\inf_{t\geq 1} a_t>0$.
Suppose that there exists a subsequence $t(i)\to\infty$
such that $a_{t(i)}\to 0$,
and we shall deduce a contradiction.
There exists a sequence $s_i\in A^0(\End_0(E))$
such that
(i) $\|s_i\|_{L^2,h_{t(i)}}=1$,
(ii) $\Delta_{t(i)}s_i=a_{t(i)}\cdot s_i$.
We obtain
\[
 \|t[\theta,s_i]\|_{L^2,h_{t(i)}}
+\|t[\theta^{\dagger}_{h_{t(i)}},s_i]\|_{L^2,h_{t(i)}}
+\|\delbar_E(s_i)\|_{L^2,h_{t(i)}}
+\|\del_{E,h_{t(i)}}(s_i)\|_{L^2,h_{t(i)}}
\to 0.
\]
By taking a subsequence,
we may assume that
the sequence $s_i$ is weakly convergent in $L_1^2$
locally on $X^{\circ}$.
Let $s_{\infty}$ denote the limit.
The sequence $s_i$ $(i=1,2,\ldots)$ is convergent to $s_{\infty}$
in $L^p$ for any $p>1$ locally on $X^{\circ}$.
We obtain $\theta(s_{\infty})=0$
and $\delbar_E(s_{\infty})=0$.
We also have $\|s_{\infty}\|_{L^2,h_{\infty}}\leq 1$.
Because $\tr s_i=0$,
we obtain $\tr s_{\infty}=0$.

We set
$b_i:=\max_X|s_{i}|_{h_{t(i)}}$.
Because
$\|s_i\|_{L^2,h_{t(i)}}=1$,
there exists $C_1>0$
such that $b_i\geq C_1$.

\begin{lem}
\label{lem;23.5.4.2}
The sequence $b_i$ is bounded.
\end{lem}
\pf
We set $\stilde_{i}=b_i^{-1}s_{i}$
and $h_i=h_{t(i)}$.
We take $P_i\in X$
such that
$|\stilde_{i|P_i}|_{h_i}=1$.
By taking a subsequence,
we may assume that the sequence $P_i$ $(i=1,2,\ldots)$
is convergent.
Let $P_{\infty}$ denote the limit.

Let $U$ be a relatively compact
neighbourhood of $P_{\infty}$ in $X$
with a holomorphic coordinate $z$.
Let $(\delbar_E+t\ad\theta)^{\ast}_{h_i}$
denote the formal adjoint
of $\delbar_E+t\ad\theta$
with respect to
$h_i$ and the Euclidean metric $dz\,d\zbar$.
By Lemma \ref{lem;23.5.4.10} below,
we obtain 
\[
 -\del_{z}\del_{\zbar}|\stilde_{i}|_{h_{i}}^2
 \leq
 2
 \bigl|
 (\delbar_E+t\ad\theta)^{\ast}_{h_i}
 (\delbar_E+t\ad\theta)(\stilde_i)
 \bigr|_{h_{i}}
 \cdot
 |\stilde_i|_{h_{i}}
 =2a_{t(i)}g_U
 |\stilde_i|_{h_{i}}^2
\leq 2a_{t(i)}g_U.
\]
Here, $g_U$ is an $\real_{>0}$-valued
$C^{\infty}$-function on $U$
such that $g_U$ and $g_U^{-1}$ are bounded.
There exists a bounded $C^{\infty}$-function $\beta_U$ on $U$
such that 
$\del_z\del_{\zbar}\beta_U=g_U$.
We obtain
\[
  -\del_{z}\del_{\zbar}
\bigl(
  |\stilde_{i}|_{h_{i}}^2
  -2a_{t(i)}\beta_U
\bigr)
 \leq 0.
\]

Let us consider the case where
$P_{\infty}\in X^{\circ}$.
We may assume that $U$ is relatively compact in $X^{\circ}$.
There exists $\epsilon>0$
such that
$B(P_{\infty},\epsilon)=
\bigl\{P\in U\,\big|\,|z(P)-z(P_{\infty})|<\epsilon\bigr\}$
is relatively compact subset in $U$.
By the mean value property
of the subharmonic function
$|\stilde_i|^2_{h_i}-2a_{t(i)}\beta_U$,
we obtain
\[
 (|\stilde_{i}|_{h_i}^2-2a_{t(i)}\beta_U)(P)
 \leq
 \frac{1}{\pi\epsilon^2}
 \int_{B(P_{\infty},\epsilon)}
 \bigl(
 |\stilde_i|^2_{h_i}-2a_{t(i)}\beta_U
 \bigr)
 \frac{|dz\,d\zbar|}{2}.
\]
Because $a_{t(i)}\to 0$,
there exists $C_4>0$ such that the following holds:
\[
 C_4\leq
 \int_{B(P_{\infty},\epsilon)}|\stilde_i|_{h_{i}}^2
=b_i^{-2} \int_{B(P_{\infty},\epsilon)}|s_{i}|_{h_{i}}^2.
\]
Because $s_i$ is convergent to $s_{\infty}$ in $L^2$ on $U$,
there exists $C_5>0$ such that $b_i\leq C_5$.

Let us consider the case $P_{\infty}\in D(\theta)$.
Let $U'$ be a relatively compact neighbourhood of $P_{\infty}$ in $U$
such that $U'\cap D(\theta)=\{P_{\infty}\}$.
By applying the maximum principle
to the subharmonic function
$|\stilde_{i}|_{h_i}^2-2a_{t(i)}\beta_U$,
there exist $Q_i\in \del U'$
such that
\[
\bigl(
 |\stilde_{i}|_{h_i}^2
 -2a_{t(i)}\beta_U
\bigr)(P_i)
 \leq
 \bigl(
  |\stilde_{i}|_{h_i}^2
  -2a_{t(i)}\beta_U
  \bigr)(Q_i).
\]
Because $a_{t(i)}\to 0$,
there exists $C_6>0$
such that
\[
 |\stilde_{i|Q_i}|_{h_i}^2\geq C_6.
\]
By taking a subsequence,
we may assume that the sequence $Q_i$ is convergent.
Let $Q_{\infty}\in \del U'$ denote the limit.
There exists $\epsilon>0$
such that
$B(Q_{\infty},\epsilon)
=\bigl\{
 P\in U\,\big|\,
 |z(P)-z(Q_{\infty})|<\epsilon
\bigr\}$
is relatively compact in $U\setminus\{P_{\infty}\}$.
By the mean value property
of the subharmonic function
$|\stilde_i|_{h_i}^2-a_{t(i)}\beta_U$,
there exists $C_7>0$
such that
\[
 C_7\leq
 \int_{B(Q_{\infty},\epsilon)}|\stilde_i|_{h_i}^2
 =b_i^{-2}
  \int_{B(Q_{\infty},\epsilon)}|s_i|_{h_i}^2.
\]
Because
$s_i$ is convergent in $L^2$ on $B(Q_{\infty},\epsilon)$,
the sequence $b_i$ is bounded.
Thus, we obtain Lemma \ref{lem;23.5.4.2}.
\hfill\qed

\vspace{.1in}
By the Lebesgue theorem
and the boundedness of the sequence $b_i$,
we obtain the following lemma.

\begin{lem}
$s_{\infty}$
is a non-zero bounded section of $E$
on $X^{\circ}$.
\hfill\qed
\end{lem}

Hence, $s_{\infty}$ is an automorphism of
the filtered Higgs bundle
$(\nbigp^{h_{\infty}}_{\ast}E,\theta)$.
(See \cite{Mochizuki-KH-Higgs} for filtered Higgs bundles,
for example.)
Because $\Sigma_{\theta}$ is connected,
$(\nbigp^{h_{\infty}}_{\ast}E,\theta)$
is stable,
and hence we obtain $s_{\infty}=\alpha\cdot\id_E$
for some $\alpha\in\cnum$.
Because $\tr s_{\infty}=0$,
we obtain $\alpha=0$,
i.e., $s_{\infty}=0$.
It contradicts $s_{\infty}\neq 0$.
Thus, we obtain a contradiction as desired,
and the proof of Proposition \ref{prop;23.4.12.1}
is completed.
\hfill\qed

\vspace{.1in}
We recall the following general lemma.
\begin{lem}
\label{lem;23.5.4.10}
Let $(E,\theta,h)$ be a harmonic bundle on $X$.
Let $U$ be an open subset of $X$
with a holomorphic coordinate $z$.
Let $f$ be the automorphism of $E_{|U}$
determined by $\theta=f\,dz$.
Let $\nabla$ denote the Chern connection of $(E,h)$. 
Then, for any $s\in A^0(U,E)$,
we have
\begin{equation}
\label{eq;23.5.4.1}
 -\del_z\del_{\zbar}|s|_h^2
=2\Re\bigl(
 (-\nabla_z\nabla_{\zbar}+f_h^{\dagger}\circ f)s,s
 \bigr)
-|\nabla_zs|_h^2
-|\nabla_{\zbar}s|_h^2
-|f(s)|_h^2
-|f^{\dagger}_h(s)|_h^2.
\end{equation} 
\end{lem} 
\pf
Note that
$\bigl(
 \nabla_{z}\nabla_{\zbar}
-\nabla_{\zbar}\nabla_z
 \bigr)
+[f,f^{\dagger}_h]=0$.
We obtain
\begin{multline}
 -\del_z\del_{\zbar}|s|_{h}^2
=
-h\bigl(
\nabla_{z}\nabla_{\zbar}s,s
 \bigr)
-h\bigl(s,
\nabla_{\zbar}\nabla_{z}s
 \bigr)
-\bigl|\nabla_{\zbar}s\bigr|_h^2
-\bigl|\nabla_{z}s\bigr|_h^2
 \\
=-h\bigl(
\nabla_{z}\nabla_{\zbar}s,s
 \bigr)
-h\bigl(s,
\nabla_{z}\nabla_{\zbar}s
 \bigr)
-h\bigl(s,[f,f^{\dagger}_h]s\bigr)
-\bigl|\nabla_{\zbar}s\bigr|_h^2
-\bigl|\nabla_{z}s\bigr|_h^2
 \\
=2\Re h\bigl(
-\nabla_z\nabla_{\zbar}s,s
\bigr)
+2h\bigl(s,f^{\dagger}_h\circ f(s)\bigr)
-\bigl|f(s)\bigr|_h^2
-\bigl|f^{\dagger}_h(s)\bigr|_h^2
-\bigl|\nabla_{\zbar}s\bigr|_h^2
-\bigl|\nabla_{z}s\bigr|_h^2
\\
=2\Re h\bigl(
 (-\nabla_z\nabla_{\zbar}+f^{\dagger}_hf)s,s
 \bigr)
-\bigl|f(s)\bigr|_h^2
-\bigl|f^{\dagger}_h(s)\bigr|_h^2
-\bigl|\nabla_{\zbar}s\bigr|_h^2
-\bigl|\nabla_{z}s\bigr|_h^2.
\end{multline}
Thus, we obtain (\ref{eq;23.5.4.1}).
\hfill\qed

\section{Horizontal deformations}
\label{section;23.5.8.20}

\subsection{Line bundles on spectral curves}

\subsubsection{Spectral curves}
\label{subsection;23.4.17.60}

We recall the spectral curves by following
\cite{Beauville-Narasimhan-Ramanan}.
Let $X$ be a compact Riemann surface
with genus $g_X\geq 2$.
In this section,
we prefer to use the notation $K_X$
to denote the sheaf of holomorphic $1$-forms
because we consider $K_X^{j/2}$.
We set
\[
 \proj=\Proj(\Sym(\nbigo_X\oplus K_X^{-1})),
\]
i.e., $\proj$ denotes the projective completion of $K_X$.
Let $\pi_{\proj}:\proj\to X$ denote the projection.
Let $\nbigo_{\proj}(1)$ denote the tautological bundle.
We have
$\pi_{\proj\ast}(\nbigo_{\proj}(1))
=\nbigo_X\oplus K_X^{-1}$.
Let $y$ denote the section of
$\nbigo_{\proj}(1)$ induced by
the section $1$ of $\nbigo_X\subset\nbigo_X\oplus K_X^{-1}$.
Let $x$ denote the section of
$\pi_{\proj}^{\ast}K_X\otimes \nbigo_{\proj}(1)$
induced by the section $1$ of
$\nbigo_X\subset K_X\otimes(\nbigo_X\oplus K_X^{-1})$.

We set $\AAA_{X,n}=\bigoplus_{j=1}^nH^0(X,K_X^{j})$.
For any $s=(s_1,\ldots,s_n)\in \AAA_{X,n}$,
we obtain the section
\[
 \xi_s=x^n+\sum_{j=1}^n (-1)^js_j y^j x^{n-j}
 \in
 H^0\bigl(
 \proj,
 \pi_{\proj}^{\ast}K_X^n
 \otimes
 \nbigo_{\proj}(n)
 \bigr).
\]
Let $Z_s$ be the subscheme of $\proj$
obtained as the $0$-set of $\xi_s$.
It is contained in $T^{\ast}X\subset\proj$.
Let $\AAA'_{X,n}$ denote the subset of $s\in\AAA_{X,n}$
such that $Z_{s}$ are smooth subschemes.
Because
$\pi_{\proj}^{\ast}K_X^n\otimes\nbigo_{\proj}(n)$
is an ample line bundle on $\proj$,
the theorem of Bertini implies that
$\AAA'_{X,n}$ is a non-empty Zariski open subset of
$\AAA_{X,n}$,
and that $Z_s$ is irreducible for any $s\in \AAA'_{X,n}$.
Let $\pi_s:Z_s\to X$ denote the projection.
For each $Q\in Z_s$,
let $r_Q$ denote the ramification index of $\pi_s$ at $Q$.
We obtain the effective divisor $\Dtilde_s=\sum (r_Q-1)Q$.
We have $K_{Z_s}=\pi^{\ast}(K_X)(\Dtilde_s)$.
We recall the following 
(see \cite[\S3]{Beauville-Narasimhan-Ramanan}).
\begin{lem}
For any $s\in \AAA'_{X,n}$,
the genus of $Z_s$ is $n^2(g_X-1)+1$.
It implies $\deg(\Dtilde_s)=2n(n-1)(g_X-1)$.
\hfill\qed
\end{lem}

Let $Z\subset \AAA'_{X,n}\times T^{\ast}X$ denote the universal curve.
The natural morphism $Z\to \AAA'_{X,n}$ is smooth.

\subsubsection{Singular flat metrics}
\label{subsection;23.4.17.2}

We consider $s\in \AAA'_{X,n}$
and a holomorphic line bundle $L$ on $Z_s$
of degree $n(n-1)(g_X-1)$.
We set $Z_s^{\circ}=Z_s\setminus|\Dtilde_s|$.
Because
\[
 \deg(L)-\sum_{Q}
 \frac{1}{2}(r_Q-1)
 =\deg(L)-\frac{1}{2}\deg(\Dtilde_s)
  =\deg(L)-n(n-1)(g_X-1)=0,
\]
we obtain the following lemma.
\begin{lem}
\label{lem;23.5.10.2}
There exists a Hermitian metric $h_L$ of
$L_{|Z_s^{\circ}}$
satisfying the following conditions.
\begin{itemize}
 \item The Chern connection $\nabla_{L}$
       of $(L_{|Z_s^{\circ}},h_L)$ is flat.
 \item For each $Q\in Z_s$,
       there exists a frame $v_Q$ of $L$ around $Q$
       such that
       $h_L(v_Q,v_Q)=|\zeta_Q|^{-r_Q+1}$,
       where $\zeta_Q$ is a local coordinate around $Q$
       such that $\zeta_Q(Q)=0$,
       and $r_Q$ is the ramification index of $\pi_s$ at $Q$.
\end{itemize}
Such a Hermitian metric is unique up to the multiplication of
positive constants. 
\hfill\qed
 \end{lem}
 
Note that for each $Q\in Z_s$
there exists a neighbourhood $(Z_s)_Q$
and a pairing
$C_Q:L_{|(Z_s)_Q}\otimes L_{|(Z_s)_Q}\lrarr
\nbigo_{Z_s}((r_Q-1)Q)$
such that
$C_Q$ is compatible with $h_L$,
i.e.,
$|C_Q(v_Q,v_Q)|=|\zeta_Q|^{-r_Q+1}$
for $v_Q$ as above.
We set $U(1)=\bigl\{a\in\cnum\,\big|\,|a|=1\bigr\}$.
The monodromy representation
induces a homomorphism
$H_1(Z_s\setminus|\Dtilde_s|,\seisuu)\to U(1)$
denoted by $\kappa_{L}$.

\begin{lem}
\label{lem;23.4.16.10}
Let $L_i$ $(i=1,2)$ be holomorphic line bundles on $Z_s$
of degree $n(n-1)(g_X-1)$.
Let $U_j$  $(j=1,\ldots,m)$ be simply connected open subsets of $Z_s$
such that $U_{j}\cap U_k=\emptyset$ $(j\neq k)$
and that
$|\Dtilde_s|\subset\bigcup U_j$.
Let 
\[
 \iota:H_1\left(Z_s\setminus\bigcup U_j,\seisuu\right)
 \lrarr
 H_1\left(Z_s\setminus|\Dtilde_s|,\seisuu\right)
\]
denote the morphism induced by the inclusion.
If $\kappa_{L_1}\circ\iota=\kappa_{L_2}\circ\iota$,
there exists an isomorphism $L_1\simeq L_2$.
\end{lem}
\pf
For any $Q\in|\Dtilde_s|$,
let $\gamma_Q$ be a loop around $Q$.
We obtain $\kappa_{L_i}(\gamma_Q)=(-1)^{r_Q-1}$
by the definition of $h_{L_i}$.
Hence, the assumption of Lemma \ref{lem;23.4.16.10}
implies $\kappa_{L_1}=\kappa_{L_2}$.
We obtain the isomorphism of the flat bundles
$(L_{1|Z_s\setminus|\Dtilde_s|},\nabla_{L_1})
\simeq
(L_{2|Z_s\setminus|\Dtilde_s|},\nabla_{L_2})$.
It extends to an isomorphism of
holomorphic line bundles
$L_1\simeq L_2$.
\hfill\qed

\subsubsection{A special case}

We fix a line bundle $K_X^{1/2}$
and an isomorphism $(K_X^{1/2})^2\simeq K_X$.
We have
\begin{equation}
 \deg(\pi_s^{\ast}K_X^{(n-1)/2})=n(n-1)(g_X-1)
=\frac{1}{2}\deg(\Dtilde_s).
\end{equation}

\begin{prop}
\label{prop;23.4.16.11}
The image of
 $\kappa_{\pi_s^{\ast}(K_X^{(n-1)/2})}:
 H_1(Z_s\setminus|\Dtilde_s|,\seisuu)\to U(1)$
is contained in $\{\pm 1\}$.
\end{prop}
\pf
Recall that $Z_s$ is the $0$-set of
$\xi_s\in H^0(\proj,\pi_{\proj}^{\ast}K_X^n\otimes\nbigo_{\proj}(n))$,
and that the line bundle $\nbigo_{\proj}(n)$ is trivialized
by the section $\pi_{\proj}^{\ast}(y)$
on a neighbourhood of $Z_s$.
(See \S\ref{subsection;23.4.17.60} for the section $y$.)
Hence, there exists a natural isomorphism
between $\pi_s^{\ast}K_X^n$
and the normal bundle of $Z_s\subset T^{\ast}X$.
By the symplectic structure of $T^{\ast}X$,
we obtain the isomorphism
$\pi_s^{\ast}K_X^n\simeq K_{Z_s}$.
Hence, there exists a natural isomorphism
\begin{equation}
 \label{eq;23.4.30.1}
 C_{ \pi_s^{\ast}(K_X^{(n-1)/2})}:
 \pi_s^{\ast}(K_X^{(n-1)/2})
 \otimes
 \pi_s^{\ast}(K_X^{(n-1)/2})
 \simeq
 \pi_s^{\ast}(K_X^{n-1})
 \simeq
 K_{Z_s}\otimes \pi_s^{\ast}(K_X^{-1})
 \simeq
 \nbigo_{Z_s}(\Dtilde_s).
\end{equation}
For any local section $u$ of $\pi_s^{\ast}K_X^{(n-1)/2}$,
we define
\[
 h_{\pi_s^{\ast}K_X^{(n-1)/2}}(u,u)
=|C_{ \pi_s^{\ast}(K_X^{(n-1)/2})}(u,u)|.
\]
In this way,
we obtain a singular flat metric $h_{\pi_s^{\ast}K_X^{(n-1)/2}}$
as in Lemma \ref{lem;23.5.10.2}.
Because $h_{\pi_s^{\ast}(K_X^{(n-1)/2})}$
is compatible with the non-degenerate symmetric bilinear product,
we obtain the claim of Proposition \ref{prop;23.4.16.11}.
(See \cite[Lemma 2.7]{Mochizuki-Szabo}.)
\hfill\qed

\subsubsection{Some representation}
\label{subsection;23.5.10.3}

Let $L$ be a holomorphic line bundle on $Z_s$
of degree $n(n-1)(g_X-1)$.
We obtain the holomorphic line bundle
$\Lhat=L\otimes \pi_s^{\ast}(K_X^{(n-1)/2})^{-1}$ of degree $0$.
It is equipped with the unique unitary flat connection
$\nabla_{\Lhat}$
whose $(0,1)$-part equals the holomorphic structure of $\Lhat$.
The monodromy representation of $\nabla_{\Lhat}$
induces a homomorphism
$\kappahat_L:H_1(Z_s,\seisuu)\to U(1)$.

\begin{rem}
Let $(K_X^{1/2})'$ be another square root of $K_X$. 
We obtain the holomorphic line bundle
$M=(K_X^{1/2})^{-1}\otimes (K_X^{1/2})'$ of degree $0$ on $X$.
It is equipped with a unitary flat connection $\nabla_M$
which induces $\kappa_M:H_1(X,\seisuu)\to \{\pm 1\}$.
It induces
$\pi_s^{\ast}(\kappa_M):H_1(Z_s,\seisuu)\to \{\pm 1\}$.
Let $\kappahat'_{L}:H_1(Z_s,\seisuu)\to U(1)$ be
the homomorphism
obtained from $L$ and $(K_X^{1/2})'$.
Then, we have
$\kappahat_{L}=\kappahat'_{L}\otimes \pi_s^{\ast}(\kappa_M)$.
\hfill\qed 
\end{rem}

There exists a natural isomorphism
$L\simeq
 \Lhat\otimes
 \pi_s^{\ast}(K_X^{(n-1)/2})$.
For $\gamma\in H_1(Z_s\setminus|\Dtilde_s|,\seisuu)$,
we have
\[
 \kappa_{L}(\gamma)
 =\kappahat_{L}(\gamma)
 \cdot\kappa_{\pi_s^{\ast}(K_X^{(n-1)/2})}(\gamma)
 \in U(1).
\]

\subsection{Horizontal deformation of line bundles}

\subsubsection{Horizontal family of line bundles}

Let $B$ be a simply connected Hausdorff space
with a base point $b_0$.
Let $\varphi:B\to \AAA'_{X,n}$ be a continuous map.
We set $Z_b:=Z_{\varphi(b)}$ for any $b\in B$.
Let $L_{b}$ $(b\in B)$ be holomorphic line bundles
on $Z_{b}$ of degree $n(n-1)(g_X-1)$. 
For any $b_1,b_2\in B$,
we choose a continuous path $\gamma$ in $B$
such that $\gamma(0)=b_1$ and $\gamma(1)=b_2$,
and then we obtain the isomorphisms
$H_1(Z_{b_1},\seisuu)\simeq
H_1(Z_{b_2},\seisuu)$
along the path $\varphi\circ\gamma$,
which is independent of the choice of $\gamma$.

\begin{df}
We say that the family $L_b$ $(b\in B)$ is horizontal
if $\kappahat_{L_b}:H_1(Z_{b},\seisuu)\to U(1)$
are independent of $b\in B$.
This condition is independent of the choice of $K_X^{1/2}$.
\hfill\qed
\end{df}

The following lemma is obvious.
\begin{lem}
Let $N$ be any holomorphic line bundle on $X$
with an isomorphism $N^2\simeq K_X$.
Then,  
the family $\pi_{\varphi(b)}^{\ast}(N^{n-1})$ $(b\in B)$
is horizontal. 
\hfill\qed
\end{lem}

The following lemma is well known.
\begin{lem}
\label{lem;23.4.17.1}
Let $L$ be any holomorphic line bundle on $Z_{b_0}$
of degree $n(n-1)(g_X-1)$.
There exists
a horizontal family $L_b$ $(b\in B)$ such that $L_{b_0}\simeq L$.
If $L_b'$ $(b\in B)$ is another horizontal family
such that $L'_{b_0}\simeq L$,
then there exists an isomorphism $L_b\simeq L_b'$ for each $b\in B$.
\end{lem}
\pf
We set $\Lhat=L\otimes \pi_s^{\ast}(K_X^{(n-1)/2})^{-1}$,
which is a holomorphic line bundle of degree $0$.
We obtain the homomorphism
$\kappahat_{L}:H_1(Z_{b_0},\seisuu)\to U(1)$.
It induces homomorphisms $H_1(Z_{b},\seisuu)\to U(1)$
and unitary flat line bundles $(\Lhat_b,\nabla_{\Lhat_b})$ on $Z_b$
$(b\in B)$.
Let $\Lhat_b$ denote the underlying holomorphic line bundle.
We set
$L_b=\Lhat_b\otimes \pi_{\varphi(b)}^{\ast}(K_X^{(n-1)/2})$.
Then, $L_b$ $(b\in B)$ is a horizontal family.

Let $L'_b$ $(b\in B)$ be another horizontal family
such that $L'_{b_0}\simeq L$.
Because $\kappahat_{L_b}=\kappahat_{L'_b}$ for any $b\in B$,
there exist isomorphisms
$L_b\otimes \pi_{\varphi(b)}^{\ast}(K_X^{(n-1)/2})^{-1}
\simeq
L'_b\otimes \pi_{\varphi(b)}^{\ast}(K_X^{(n-1)/2})^{-1}$,
and hence
$L_b\simeq L_b'$.
\hfill\qed

\subsubsection{Complex analytic property of horizontal families}
\label{subsection;23.4.17.40}

Let $B$ be a multi-disc
with the base point $b_0=(0,\ldots,0)$.
Let $\varphi:B\to \AAA'_{X,n}$ be a holomorphic map.
We set $Z_B:=Z\times_{\AAA'_{X,n}}B$.
The natural morphism $q_0:Z_B\to B$ is smooth.
We use the notation
$\Dtilde_b$ and $\pi_b$
instead of 
$\Dtilde_{\varphi(b)}$
and $\pi_{\varphi(b)}$,
respectively.
Let $\pi_B:Z_B\to B\times X$ be the induced morphism.
For any holomorphic line bundle $\nbigl$ on $Z_B$,
we set $\nbigl_b=\nbigl_{|Z_b}$ $(b\in B)$.
If the family $\nbigl_{b}$ $(b\in B)$ is horizontal,
we say that $\nbigl$ is horizontal.

\begin{lem}
\label{lem;23.4.17.13}
Let $L$ be any holomorphic line bundle on $Z_{b_0}$
of degree $n(n-1)(g_X-1)$.
Then, there exists a horizontal holomorphic line bundle $\nbigl$ on $Z_B$
such that $\nbigl_{b_0}\simeq L$.
If $\nbigl'$ is another horizontal holomorphic line bundle on $Z_B$
such that $\nbigl'_{b_0}\simeq L$,
then there exists a holomorphic isomorphism
$\nbigl\simeq\nbigl'$. 
\end{lem}
\pf
We obtain the homomorphism
$\kappahat_{L}:H_1(Z_{b_0},\seisuu)\to U(1)$
as in the proof of Lemma \ref{lem;23.4.17.1}.
Because $Z_B$ is homotopy equivalent to $Z_{b_0}$,
we obtain $\kappahat_L':H_1(Z_B,\seisuu)\to U(1)$.
There exists a unitary flat line bundle
$(\nbiglhat,\nabla_{\nbiglhat})$ on $Z_B$
which induces the homomorphism $\kappahat_L'$.
Let $q_2:B\times X\to X$ denote the projection.
We obtain the line bundle
$(q_2\circ\pi_B)^{\ast}(K_X^{(n-1)/2})$ on $Z_B$.
We set
\[
 \nbigl=
 \nbiglhat\otimes
 (q_2\circ\pi_B)^{\ast}(K_X^{(n-1)/2}).
\]
Then, $\nbigl$ is horizontal,
and there exists an isomorphism
$\nbigl_{b_0}\simeq L$.

Let $\nbigl'$ be a horizontal holomorphic line bundle
with an isomorphism
$\nbigl'_{b_0}\simeq L$.
We set $\nbigm=\nbigl'\otimes\nbigl^{-1}$.
By Lemma \ref{lem;23.4.17.1},
$\nbigm_{b}\simeq\nbigo_{Z_b}$ for any $b\in B$.
Hence,
$q_{0\ast}(\nbigm)$ is a holomorphic line bundle on $B$.
Because $B$ is a multi-disc,
there exists an isomorphism
$q_{0\ast}(\nbigm)\simeq\nbigo_B$.
A nowhere vanishing section of 
$q_{0\ast}(\nbigm)$
induces an isomorphism $\nbigl\simeq\nbigl'$.
\hfill\qed

\begin{rem}
We shall also explain another construction of
a horizontal holomorphic line bundle
in terms of locally defined symmetric pairings
in {\rm\S\ref{subsection;23.4.17.12}}.
\hfill\qed
\end{rem}

\subsubsection{Horizontal property and symmetric pairings}

We continue to use the notation in \S\ref{subsection;23.4.17.40}.
Let $\nbigr_{\pi_B}$ denote the ramification divisor 
of $\pi_B$ on $Z_B$.
(See \S\ref{subsection;23.4.30.2}.)
Let $H$ be a closed complex analytic hypersurface of $Z_B$
such that $H\to B$ is proper and finite,
and that $|\nbigr_{\pi_B}|\subset H$.
For any $b\in H$, we set $H_b:=H\cap Z_{b}$.
We assume the following.

\begin{condition}
\label{condition;23.4.17.20}
There exist open subsets
$(Z_{B})_0$ and $(Z_B)_Q$ $(Q\in H_{b_0})$
such that the following holds.
 \begin{itemize}  
  \item $Z_B=(Z_B)_0\cup\bigcup_{Q\in H_{b_0}}(Z_B)_Q$.
  \item $(Z_B)_0$, $(Z_B)_Q$
	and $(Z_B)_0\cap (Z_B)_Q$
	are Stein.
  \item $\overline{(Z_B)_Q}\cap \overline{(Z_B)_{Q'}}=\emptyset$
	if $Q\neq Q'$.
  \item The boundaries $\del (Z_B)_0$
	and $\del (Z_B)_Q$
	are submanifolds of $Z_B$,
	and the projections $\del (Z_B)_0\to B$
	and $\del (Z_B)_Q\to B$ are proper and fiber bundles.
  \item There exists a bi-holomorphic isomorphism
       $(Z_B)_0\simeq B\times ((Z_B)_0\cap Z_{b_0})$
       over $B$.
 \item Each $(Z_B)_Q\cap Z_{b_0}$
       is bi-holomorphic to an open disc.
       Each
       $(Z_B)_0\cap
       (Z_B)_Q\cap Z_{b_0}$
       is homeomorphic to
       an annulus $\{1/2<|z|<1\}$.
\hfill\qed
 \end{itemize} 
\end{condition}

\begin{rem}
We shall explain a way to obtain such an open covering
in {\rm\S\ref{subsection;23.4.17.30}}.
\hfill\qed
\end{rem}

We set
$(Z_b)_0:=(Z_B)_0\cap Z_b$,
$(Z_b)_Q:=(Z_B)_Q\cap Z_b$
and
$(Z_b)_{0,Q}:=(Z_B)_0\cap(Z_B)_Q\cap Z_b$.
We also set
$(Z_B)_{0,Q}=(Z_B)_0\cap (Z_B)_Q$.
We note the following lemma which follows from the conditions.
\begin{lem}
For any $Q\in H_{b_0}$ and $b\in B$,
$(Z_b)_Q$ is homeomorphic to an open disc,
and
$(Z_b)_{0,Q}$ is homeomorphic to
an annulus.
\hfill\qed
\end{lem}

Let $p_{0}:(Z_B)_0\to (Z_{b_0})_0$ denote the morphism
induced by the isomorphism
$(Z_B)_0\simeq B\times (Z_{b_0})_0$
and the projection
$B\times(Z_{b_0})_0\to (Z_{b_0})_0$.
We shall prove the following proposition
in \S\ref{subsection;23.4.17.11}--\ref{subsection;23.5.17.20}.

\begin{prop}
\label{prop;23.4.17.10}
Let $\nbigl$ be a holomorphic line bundle on $Z_B$
such that $\deg\nbigl_b=n(n-1)(g_X-1)$ for any $b\in B$.
Let $h_{\nbigl_{b_0}}$
denote the singular flat metric of 
$\nbigl_{b_0}$ as in {\rm\S\ref{subsection;23.4.17.2}}.
Then, $\nbigl$ is horizontal
if and only if  
there exist an isomorphism
$\varphi_0:p_0^{\ast}(\nbigl_{|(Z_{b_0})_0})
\simeq
\nbigl_{|(Z_B)_0}$
and 
non-degenerate symmetric pairings
$C_Q:\nbigl_{|(Z_B)_Q}\otimes
 \nbigl_{|(Z_B)_Q}
 \simeq
 \nbigo_{Z_B} (\nbigr_{\pi_B})_{|(Z_B)_Q}$
$(Q\in H_{b_0})$ 
such that the following holds.
\begin{itemize}
 \item $C_{Q|(Z_{b_0})_{Q}}$
       is compatible with
       $h_{\nbigl_{b_0}|(Z_{b_0})_Q}$.
       Moreover,
       on $(Z_B)_{0,Q}$,
       we have
       $p_0^{\ast}(C_{Q|(Z_{b_0})_{0,Q}})
       =C_Q$
       under the isomorphism $\varphi_{0}$.       
\end{itemize}
\end{prop}

\begin{rem}
As we shall see in {\rm\S\ref{subsection;23.4.17.11}},
for each $b\in B$,
there exists a singular flat metric
$h_{\nbigl_b}$
of $\nbigl_b$
such that
$h_{\nbigl_b|(Z_b)_Q}$
is compatible with $C_{Q|(S_b)_Q}$.
\hfill\qed 
\end{rem}

\subsubsection{The ``if'' part of Proposition \ref{prop;23.4.17.10}}
\label{subsection;23.4.17.11}

Let us prove the ``if'' part.
Let $b\in B$.
Let $h_{0,b}$ be 
the flat Hermitian metric
of $\nbigl_{b|(Z_b)_0}$
induced by
$h_{\nbigl_{b_0}|(Z_{b_0})_0}$
and the isomorphism $\varphi_0$.
Let $h_{Q,b}$ be the singular flat Hermitian metric
of $\nbigl_{b|(Z_b)_Q}$
compatible with $C_{Q|(Z_b)_Q}$.
By the third condition,
we obtain a singular flat Hermitian metric
$h_{b}$ of $\nbigl_{b}$
from $h_{0,b}$ and $h_{Q,b}$ $(Q\in H_{b_0})$.
It is a singular flat Hermitian metric $h_{\nbigl_b}$
as in Lemma \ref{lem;23.5.10.2}.

We remark the following lemma
as a consequence of Proposition \ref{prop;23.4.16.11}.
\begin{lem}
\label{lem;23.4.16.12}
For $b\in B$,
let $\kappa_{b}:H_1((Z_b)_0,\seisuu)\to U(1)$
be the representation 
induced by the Chern connection of
$\pi_b^{\ast}(K_X^{(n-1)/2})$
with the singular metric $h_{\pi_b^{\ast}(K_X^{(n-1)/2})}$.
Then, $\kappa_b$ are independent of $b\in B$.
\hfill\qed
\end{lem}

Let $L_b$ be a holomorphic line bundle on $Z_{b}$
such that $\kappahat_{L_b}$
is equal to
$\kappahat_{\nbigl_{b_0}}$.
(See \S\ref{subsection;23.5.10.3}
for the representation $\kappahat_{L_b}$.)
Let
\[
 \iota_b:
 H_1\Bigl(
 Z_b\setminus\bigcup (Z_b)_Q,\seisuu
 \Bigr)
 \lrarr
 H_1\bigl(Z_b,\seisuu\bigr)
\]
denote the homomorphism induced by the inclusion.
By Lemma \ref{lem;23.4.16.12}
and the construction of $h_{\nbigl_b}$,
we obtain
$\kappa_{L_b}\circ\iota
=\kappa_{\nbigl_b}\circ\iota$.
Hence, by Lemma \ref{lem;23.4.16.10},
$L_b$ is isomorphic to
$\nbigl_{b}$,
i.e., $\nbigl$ is horizontal.

\subsubsection{Construction of horizontal family via symmetric pairings}
\label{subsection;23.4.17.12}

Let $L$ be a holomorphic line bundle on $Z_{b_0}$
of degree $n(n-1)(g_X-1)$.
We obtain the holomorphic line bundle
$\nbigl_0=p_0^{\ast}(L_{|(Z_{b_0})_0})$
on $(Z_B)_0$.
We set $\nbigl_Q=\nbigo_{(Z_B)_Q}$.
We fix an isomorphism
$\nbigl_{Q|(Z_{b_0})_Q}
 \simeq
 L_{|(Z_{b_0})_Q}$.

\begin{lem}
For each $Q$,
there exists an isomorphism
$\varphi_{1,Q}:
 \nbigl_{0|(Z_B)_{0,Q}}
 \simeq
 \nbigl_{Q|(Z_B)_{0,Q}}$
such that
$\varphi_{1,Q|(Z_{b_0})_{0,Q}}$
equals the identity.
\end{lem}
\pf
We use Oka-Grauert principle
(see \cite[\S1.1.2, Theorem B']{Gromov}).
There exists a topological isomorphism
$\nbigl_{0|(Z_B)_{0,Q}}
 \simeq
 \nbigl_{Q|(Z_B)_{0,Q}}$.
Because $(Z_B)_{0,Q}$ is Stein,
there exists a holomorphic isomorphism
$\varphi'_{1,Q}:
\nbigl_{0|(Z_B)_{0,Q}}
\simeq
\nbigl_{Q|(Z_B)_{0,Q}}$.
The restriction
$\varphi'_{1,Q|(Z_{b_0})_{0,Q}}$ can be regarded
as a holomorphic function
$\alpha:(Z_{b_0})_{0,Q}\to\cnum^{\ast}$.
Both $(Z_{b_0})_{0,Q}$ and $\cnum^{\ast}$
are homotopic to
$S^1=\bigl\{a\in\cnum\,\big|\,|a|=1\bigr\}$.
Because $(Z_B)_{0,Q}$ is Stein,
we can modify $\varphi'_{1,Q}$
such that
the induced map
$H_1((Z_{b_0})_{0,Q},\seisuu)
\to
H_1(\cnum^{\ast},\seisuu)$
is $0$.
Then,
there exists $\beta:(Z_{b_0})_{0,Q}\to \cnum$
such that $\alpha=\exp(\beta)$.
Because $(Z_B)_{0,Q}$ is Stein,
there exists $\betatilde:(Z_B)_{0,Q}\to\cnum$
whose restriction to $(Z_{b_0})_{0,Q}$
equals $\beta$.
Then,
$\varphi_{1,Q}
=\exp(-\betatilde)\varphi'_{1,Q}$
satisfies the desired condition.
\hfill\qed

\vspace{.1in}

Because $(Z_B)_Q$ is Stein,
there exists a non-degenerate symmetric bilinear product
\[
 \Ctilde_Q:
 \nbigl_Q\otimes\nbigl_Q
 \simeq
 \nbigo_{Z_B}(\nbigr_{\pi_B})_{|(Z_B)_Q}.
\]
Let $h_L$ be the singular flat Hermitian metric of $L$
as in \S\ref{subsection;23.4.17.2}.
There exists a non-degenerate symmetric pairing
\[
C_{Q}:
L_{|(Z_{b_0})_Q}
\otimes
L_{|(Z_{b_0})_Q}
\simeq
\nbigo_{Z_{b_0}}(\nbigr_{\pi_{b_0}})
_{|(Z_{b_0})_Q}
\]
compatible with $h_L$.
We obtain the nowhere vanishing holomorphic function $g_Q$
on $(Z_B)_{0,Q}$
determined by
\[
 p_0^{\ast}C_Q
=g_Q\cdot \varphi_{1,Q}^{\ast}\Ctilde_Q.
\]
\begin{lem}
There exists a holomorphic function $g_Q^{1/2}$
on $(Z_B)_{0,Q}$
such that $(g_Q^{1/2})^{2}=g_Q$. 
\end{lem}
\pf
The restriction of $g_Q$
to $(Z_{b_0})_{0,Q}$
has a square root.
Then, it is easy to see that
$g_Q$ also has a square root.
\hfill\qed

\vspace{.1in}

We obtain the following isomorphisms:
\[
 \varphi_{2,Q}=g_Q^{-1/2}\varphi_{1,Q}:
 \nbigl_{0|(Z_B)_{0,Q}}
 \simeq
 \nbigl_{Q|(Z_B)_{0,Q}}.
\]
We obtain a holomorphic line bundle $\Psi(L)$ on $Z_B$
by gluing $\nbigl_0$ and $\nbigl_Q$
$(Q\in H_{b_0})$
via the isomorphisms $\varphi_{2,Q}$.
Because we have already proved
``the if part'' of Proposition \ref{prop;23.4.17.10}
in \S\ref{subsection;23.4.17.11},
we obtain that $\Psi(L)$ is horizontal.
Because $(g^{1/2}_Q)_{|(Z_{b_0})_{0,Q}}$
extends to a nowhere vanishing holomorphic function
on $(Z_{b_0})_Q$,
we can check that $\Psi(L)_{b_0}$
is isomorphic to $L$.

\subsubsection{The ``only if'' part of Proposition \ref{prop;23.4.17.10}}
\label{subsection;23.5.17.20}

From $\nbigl_{b_0}$,
we obtain the horizontal holomorphic line bundle
$\Psi(\nbigl_{b_0})$
such that $\Psi(\nbigl_{b_0})_{b_0}\simeq \nbigl_{b_0}$
by applying the construction in \S\ref{subsection;23.4.17.12}.
By the uniqueness in Lemma \ref{lem;23.4.17.13},
$\Psi(\nbigl_{b_0})\simeq\nbigl$.
By the construction of $\Psi(\nbigl_{b_0})$,
there exist an isomorphism
$\varphi_0:
p_0^{\ast}\bigl(
\Psi(\nbigl_{b_0})_{|(Z_{b_0})_0}
\bigr)
\simeq
\Psi(\nbigl_{b_0})_{|(Z_B)_0}$
and
non-degenerate symmetric pairings
$C_Q:\Psi(\nbigl_{b_0})_{|(Z_B)_Q}\otimes
 \Psi(\nbigl_{b_0})_{|(Z_B)_Q}
 \simeq
 \nbigo_{Z_B}(\nbigr_{\pi_B})_{|(Z_B)_Q}$
$(Q\in H_{b_0})$ 
satisfying the condition
in Proposition \ref{prop;23.4.17.10}.
Thus, the ``only if'' part of Proposition \ref{prop;23.4.17.10}
is proved.
\hfill\qed

\subsubsection{Complement}
\label{subsection;23.4.17.30}
We set
$D_{b}=\pi_b(|\Dtilde_b|)$.
(See \S\ref{subsection;23.4.17.40} for $\Dtilde_b$.)
For each $P\in D_{b_0}$,
let $X_P$ be a simply connected open neighbourhood of $P$ in $X$
such that $X_P\cap D_{b_0}=\{P\}$.
Let $X_P'$ be a simply connected relatively compact
open neighbourhood of $P$ in $X_P$.
\begin{condition}
\label{condition;23.4.16.50}
$D_{b}\subset \bigcup X_P'$
for any $b\in B$.
\hfill\qed
\end{condition}
We set $X_1=X\setminus \bigcup_{P\in D_{b_0}}\overline{X'_{P}}$.
We set
$H=\pi_B^{-1}\bigl(\pi_B(|\nbigr_{\pi_B}|)\bigr)\subset Z_B$.
We set
$(Z_B)_0=\pi_B^{-1}(B\times X_1)$.
For each
$Q\in \pi_{b_0}^{-1}(D_{b_0})$,
let $(Z_B)_Q$
denote the connected component of
$\pi_B^{-1}(B\times X_{\pi_{b_0}(Q)})$
which contains $Q$.

\begin{lem}
\label{lem;23.4.17.41}
Condition {\rm\ref{condition;23.4.17.20}}
is satisfied for
$(Z_B)_0$ and $(Z_B)_Q$ $(Q\in H_{b_0})$.
\end{lem}
\pf
We can check that
$(Z_B)_0$, $(Z_B)_Q$ and $(Z_B)_{0,Q}$ are Stein
by using \cite[\S V.1 Theorem 1]{Grauert-Remmert-Stein}.
The natural morphisms
$a:(Z_B)_0\to B\times X_1$
and
$b:B\times (Z_{b_0})_{0}\to B\times X_1$
are covering maps over $B$.
Because they are the same maps over $\{b_0\}\times X_1$,
there exists the unique homeomorphism
$c:(Z_B)_0\lrarr B\times (Z_{b_0})_0$
such that $a=b\circ c$.
It is easy to see that $c$ and $c^{-1}$
are holomorphic.
The other conditions can also be checked easily.
\hfill\qed

\subsubsection{Appendix: Ramification divisors}
\label{subsection;23.4.30.2}

Let $a_i:\nbigx_i\to B$ be smooth morphisms of complex manifolds
such that $\dim \nbigx_1=\dim\nbigx_2=\dim B+1$.
Let $\varphi:\nbigx_1\to \nbigx_2$ be
a proper morphism of complex manifolds
over $B$,
i.e., $a_2\circ\varphi=a_1$.
Let $K_{\nbigx_i}$ $(i=1,2)$
denote the canonical line bundle of $\nbigx_i$.
There exists the effective divisor
$\nbigr_{\varphi}$ on $\nbigx_{1}$
such that the natural morphism
$\varphi^{\ast}K_{\nbigx_2}\to K_{\nbigx_1}$ induces
an isomorphism
\[
 K_{\nbigx_1}
 \simeq
 \varphi^{\ast}(K_{\nbigx_2})
 \otimes
 \nbigo_{\nbigx_1}(\nbigr_{\varphi}).
\]
The divisor $\nbigr_{\varphi}$ is called
the ramification divisor of $\varphi$.

In general, for a complex manifold $Y$,
let $\Omega^1_Y$ denote the sheaf of holomorphic $1$-forms on $Y$.
Let $\Omega^1_{\nbigx_i/B}$ denote the cokernel of
$a_i^{\ast}\Omega^1_B\to \Omega^1_{\nbigx_i}$,
which is a holomorphic line bundle on $\nbigx_i$.
The natural morphism
$\varphi^{\ast}\Omega^1_{\nbigx_2/B}
 \lrarr
 \Omega^1_{\nbigx_1/B}$
induces an isomorphism
\[
 \Omega^1_{\nbigx_1/B}
 \simeq
 \varphi^{\ast}\Omega^1_{\nbigx_2/B}
 \otimes
 \nbigo_{\nbigx_1}(\nbigr_{\varphi}).
\]

Let $f:B'\to B$ be any morphism of complex manifolds.
We set $\nbigx_i'=\nbigx_i\times_BB'$.
We obtain the induced morphisms
$a_i':\nbigx_i'\to B'$
and $\varphi':\nbigx_1'\to\nbigx_2'$.
We also obtain the induced morphism
$f_1:\nbigx'_1\to\nbigx_1$.
The ramification divisor is compatible with the base change,
i.e.,
$f_1^{\ast}(\nbigr_{\varphi})=\nbigr_{\varphi'}$.

\subsection{Horizontal deformations of Higgs bundles
with smooth spectral curves}

\subsubsection{Horizontal family of Higgs bundles with smooth spectral curves}

Let $B$ be a simply connected Hausdorff space.
Let $(E_b,\theta_b)$ $(b\in B)$ be
a family of Higgs bundles of rank $n$
of degree $0$ on $X$ such that the following holds.
\begin{itemize}
 \item Each $\Sigma_{\theta_b}$ is smooth.
       The induced map
       $\varphi:B\to \AAA'_{X,n}$ is continuous.
\end{itemize}
According to Hitchin and Beauville-Narasimhan-Ramanan,
for each $b\in B$,
there exists a line bundle $L_b$ on $\Sigma_{\theta_b}=Z_b$
of degree $n(n-1)(g_X-1)$
with an isomorphism
$\pi_{\varphi(b)\ast}(L_b)\simeq (E_b,\theta_b)$.
Such an $L_b$ is unique up to an isomorphism.

\begin{df}
The family $(E_b,\theta_b)$ $(b\in B)$
is called horizontal
if the family $L_b$ $(b\in B)$ is horizontal.
\hfill\qed 
\end{df}

The following lemma is well known
and an obvious translation of Lemma \ref{lem;23.4.17.1}.
\begin{lem}
Let $(E,\theta)$ be a Higgs bundle of degree $0$ on $X$
whose spectral curve is smooth.
Let $\varphi:B\to \AAA'_{X,n}$ is a continuous map.
Then, there exists a horizontal family
$(E_b,\theta_b)$ $(b\in B)$
such that $(E_{b_0},\theta_{b_0})\simeq (E,\theta)$
and $\Sigma_{\theta_{b}}=Z_{\varphi(b)}$.
For another such family $(E'_b,\theta'_b)$ $(b\in B)$,
 there exist isomorphisms
$(E_b,\theta_b)\simeq(E'_b,\theta'_b)$.
\hfill\qed 
\end{lem}

\subsubsection{Complex analytic property of horizontal families}

Let $B$ and $\varphi$ be as in \S\ref{subsection;23.4.17.40}.
\begin{df}
Let $(\nbige,\theta)$ be a relative Higgs bundle on $B\times X$
such that $\Sigma_{\nbige,\theta}=Z_B$.
We set
$(\nbige_b,\theta_{b}):=(\nbige,\theta)_{|\{b\}\times X}$.
Assume that $\deg(\nbige_b)=0$.
We say that $(\nbige,\theta)$ is horizontal
if the family $(\nbige_b,\theta_b)$ $(b\in B)$ is horizontal.
\hfill\qed
\end{df}

\begin{lem}
\label{lem;23.5.3.1}
Let $(E,\theta)$ be a Higgs bundle of degree $0$ on $X$
such that $\Sigma_{\theta}=Z_{b_0}$.
There exists
a horizontal relative Higgs bundle $(\nbige,\theta_{\nbige})$
on $B\times X$
such that
$(\nbige_{b_0},\theta_{\nbige,b_0})=(E,\theta)$
and that $\Sigma_{\nbige,\theta}=Z_B$.
If $(\nbige',\theta_{\nbige'})$
be another such relative Higgs bundle,
then there exists an isomorphism
$(\nbige,\theta_{\nbige})\simeq(\nbige',\theta_{\nbige'})$.
\end{lem}
\pf
Let $L$ be a line bundle on $Z_{b_0}$
of degree $n(n-1)(g_X-1)$
corresponding to $(E,\theta)$.
By Lemma \ref{lem;23.4.17.1},
there exists a horizontal line bundle $\nbigl$ on $Z_B$
such that $\nbigl_{b_0}\simeq L$.
According to \cite[Chapter 2, \S4.2]{Grauert-Remmert},
$\nbige=\pi_{B\ast}\nbigl$
is a locally free $\nbigo_{B\times X}$-module.
It is naturally equipped with the relative Higgs field
$\theta_{\nbige}$,
and $(\nbige_{b_0},\theta_{\nbige,b_0})\simeq(E,\theta)$.
The second claim also follows from Lemma \ref{lem;23.4.17.1}.
\hfill\qed

\subsubsection{Horizontal property and symmetric pairings}

Let $X_P$, $X_{P'}$ and $X_1$
be as in \S\ref{subsection;23.4.17.30}.
Let $(\nbige,\theta_{\nbige})$ be a relative Higgs bundle on $B\times X$
such that
(i) $\Sigma_{\nbige,\theta_{\nbige}}=Z_B$,
(ii) $\deg(\nbige_b,\theta_{\nbige,b})=0$ for any $b\in B$.

Let $q_2:B\times X_1\to X_1$ denote the projection.
Let $\nbigl_1$ be the line bundle
on $\pi_B^{-1}(B\times X_1)$
which induces $(\nbige,\theta_{\nbige})_{|B\times X_1}$.
There exists the holomorphic isomorphism
$\rho_1:B\times \pi_{b_0}^{-1}(X_1)\simeq
\pi_B^{-1}(B\times X_1)$
over $B\times X_1$.

We obtain the following proposition
as a consequence of Proposition \ref{prop;23.4.17.10}
and Lemma \ref{lem;23.4.17.41}.
\begin{prop}
\label{prop;23.4.17.50}
$(\nbige,\theta_{\nbige})$ is horizontal
if and only if
there exist a holomorphic isomorphism
 $\varphi:q_2^{\ast}
 (\nbige_{b_0|X_1})
 \simeq
  \nbige_{|B\times X_1}$
and non-degenerate symmetric pairings
$C_P$ of $(\nbige,\theta_{\nbige})_{|B\times X_P}$
$(P\in D_{b_0})$
such that the following holds.
\begin{itemize}
 \item $\varphi$ is induced by a holomorphic isomorphism
       $\rho_1^{\ast}(\nbigl_{1|\pi_{b_0}^{-1}(X_1)})
       \simeq
       \nbigl_{1|\pi_{B}^{-1}(B\times X_1)}$.
 \item Each $C_{P|\{b_0\}\times X_P}$
       is compatible with the limiting configuration
       of $(\nbige_{b_0},\theta_{\nbige,b_0})$.
       Moreover, on $B\times(X_1\cap X_P)$,
       we have
       $q_{2}^{\ast}C_{P|\{b_0\}\times (X_1\cap X_P)}
       =C_{P|B\times (X_1\cap X_P)}$
       under the isomorphism $\varphi$.
\hfill\qed
\end{itemize} 
\end{prop}

\subsection{Horizontal infinitesimal deformations}

Suppose that $X$ is compact.
Let $(E,\theta)$ be a Higgs bundle of degree $0$ on $X$
which is generically regular semisimple with smooth spectral curve.

\begin{df}\mbox{{}}
\begin{itemize}
 \item 
An infinitesimal deformation $(E^{[1]},\theta^{[1]})$
of $(E,\theta)$
is called horizontal
if there exists a relative Higgs bundle
$(\nbige,\theta_{\nbige})$ on $S\times X$
such that
(i) $(\nbige,\theta_{\nbige})$ is horizontal,
(ii) there exist $s_0\in S$ and $v\in T_{s_0}S$
such that
 $\iota_{v}^{\ast}(\nbige,\theta_{\nbige})
       \simeq (E^{[1]},\theta^{[1]})$.
       (See {\rm\S\ref{subsection;23.5.10.1}}
       for the notation.)
 \item An element of $H^1(X,\Def(E,\theta))$ is called
       horizontal
       if it corresponds to a horizontal infinitesimal deformation.
       Let $H^1(X,\Def(E,\theta))^{\hor}\subset
       H^1(X,\Def(E,\theta))$
       denote the subset of horizontal elements.
\hfill\qed
\end{itemize}
\end{df}

Let $\nbigm_H$ denote the moduli space of
stable Higgs bundles of degree $0$ and rank $n$ on $X$.
There exist the Hitchin fibration
$\Phi_H:\nbigm_H\to \AAA_{X,n}$
and its restriction
$\Phi_H:\nbigm'_H\to \AAA'_{X,n}$.
We obtain the point $(E,\theta)\in\nbigm_H'$
corresponding to $(E,\theta)$.
Let $\nbigu$ be any simply connected open neighbourhood of
$\Phi_H(E,\theta)$ in $\AAA'_{X,n}$.
According to Lemma \ref{lem;23.5.3.1},
there exists a morphism of complex manifolds
$\Psi:\nbigu\to \nbigm'_H$
such that
(i) $\Phi_H\circ\Psi=\id_{\nbigu}$,
(ii) the family of Higgs bundles
$\Psi(x)$ $(x\in\nbigu)$ is horizontal,
(iii) $\Psi(\Phi_H(E,\theta))=(E,\theta)$.

The tangent space 
$T_{(E,\theta)}\nbigm_H$
is identified with
$H^1(X,\Def(E,\theta))$.
We obtain the following lemma by definition.
\begin{lem}
 $T_{(E,\theta)}\Psi(\nbigu)
 \subset
 T_{(E,\theta)}\nbigm_H$
 is identified with
 $H^1(X,\Def(E,\theta))^{\hor}$.
In particular,
$H^1(X,\Def(E,\theta))^{\hor}$
is a $\cnum$-vector subspace of
$H^1(X,\Def(E,\theta))$
such that 
$2\dim H^1(X,\Def(E,\theta))^{\hor}=\dim H^1(X,\Def(E,\theta))$.
\hfill\qed
\end{lem}

\subsubsection{\v{C}ech description of
horizontal infinitesimal deformations}
\label{subsection;23.5.10.4}

For each $P\in D(\theta)$,
let $(X_P,z_P)$ denote a holomorphic coordinate
neighbourhood around $P$ as in \S\ref{subsection;23.4.15.1}.
We set $(E_P,\theta_P):=(E,\theta)_{|X_P}$.
Let $C_P$ be a non-degenerate symmetric pairing of
$(E_P,\theta_P)$ compatible with $h_{\infty}$
on $X_P\setminus\{P\}$.
Let $X_P'$ be a relatively compact open neighbourhood of $P$
in $X_P$.
We set $X_0=X\setminus \bigcup_{P\in D(\theta)}\overline{X_P'}$.
We obtain an open covering
$X=X_0\cup \bigcup_{P\in D(\theta)} X_P$.
We assume $X_P\cap X_{P'}=\emptyset$ $(P\neq P')$.

We may regard $\Cok(\ad\theta)_{|X_0}$
as a subsheaf of
$(\End(E)\otimes K_X)_{|X_0}$
by the splitting (\ref{eq;23.4.23.3}).
We may also regard
$\Cok(\ad\theta)_{|X_0\cap X_P}$
as a subsheaf of
$\bigl(
\Hom^{\sym}(E_P;C_P)\otimes K_X
\bigr)_{|X_0\cap X_P}$.

Then, we define a subcomplex $\gbiga^{\bullet}$ of
$\Tot \check{C}^{\bullet}(\Def(E,\theta))$
as follows (see Remark \ref{rem;23.4.19.1} for
$\Tot \check{C}^{\bullet}(\Def(E,\theta))$):
\[
 \gbiga^0:=
 \bigoplus_{P\in D(\theta)}
 H^0(X_P,\End^{\asym}(E_P;C_P)),
\]
\[
 \gbiga^1:=
 H^0(X_0,\Cok\ad\theta)
 \oplus
 \bigoplus_{P\in D(\theta)}
 H^0(X_P,\End^{\sym}(E_P;C_P)\otimes K_X)
\oplus
 \bigoplus_{P\in D(\theta)}
 H^0(X_0\cap X_P,\End^{\asym}(E_P;C_P)),
\]
\[
 \gbiga^2:=
 \bigoplus_{P\in D(\theta)}
 H^0(X_0\cap X_P,\End^{\sym}(E_P;C_P)\otimes K_X).
\]
We obtain the morphism
\[
\Upsilon_0:H^1\bigl(\gbiga^{\bullet}\bigr)
 \lrarr
 H^1(X,\Def(E,\theta)). 
\]

By the isomorphism
$\Cok(\theta)\simeq \pi_{\ast}K_{\Sigma_{\theta}}$
in Proposition \ref{prop;23.3.23.10},
we obtain the following natural morphism:
\begin{equation}
\label{eq;23.3.26.1}
 \Upsilon_1:
 H^1(X,\Def(E,\theta))
 \to
 H^0(\Sigma_{\theta},K_{\Sigma_{\theta}}).
\end{equation}

\begin{thm}
\label{thm;23.4.24.1}
The composite $\Upsilon_1\circ\Upsilon_0$ is an isomorphism.
The image of $\Upsilon_0$
equals
the space of isomorphism classes of horizontal infinitesimal deformations
$H^1(X,\Def(E,\theta))^{\hor}$.
\end{thm}
\pf
Let us prove that $\Upsilon_1\circ\Upsilon_0$ is injective.
Let $\Upsilon_1\circ\Upsilon_0(a)=0$
for $a\in H^1(\gbiga^{\bullet})$.
Let $(\varphi_0,\varphi_P,\psi_{P})\in\gbiga^1$
be a representative of $a$.
Because $\Upsilon_1\circ\Upsilon_0(a)=0$,
we obtain that
$\varphi_0=0$
and $[\varphi_P]=0$ in $\Cok(\ad\theta_P)$ $(P\in D(\theta))$.
There exists $\nu_P\in H^0(X_P,\End^{\asym}(E_P;C_P))$
such that $\ad(\theta_P)\nu_P=\varphi_P$.
Because we have
$\varphi_{P|X_P\cap X_0}=\ad\theta_P(\psi_P)$,
we obtain
$\nu_{P|X_P\cap X_0}=\psi_P$.
Hence, we obtain $a=0$ in $H^1(\gbiga^{\bullet})$.

We obtain
$\dim H^1(\gbiga^{\bullet})
\leq
 \dim H^0(\Sigma,K_{\Sigma})
=\frac{1}{2}\dim H^1(X,\Def(E,\theta))$.
By Proposition \ref{prop;23.4.17.50},
any horizontal infinitesimal deformation
is contained in the image of $\Upsilon_0$.
The dimension
of the space of the isomorphism classes of horizontal deformations
equals $\dim \AAA=\frac{1}{2}\dim H^1(X,\Def(E,\theta))$.
Hence, we obtain that
the image of $\Upsilon_0$
equals the space of the isomorphism classes of 
horizontal infinitesimal deformations.
Because
$\dim H^1(\gbiga^{\bullet})
=\dim H^0(\Sigma,K_{\Sigma})$,
$\Upsilon_1\circ\Upsilon_0$ is an isomorphism.
\hfill\qed

\begin{cor}
$\Upsilon_0$ is injective,
and 
$\Upsilon_1$ is surjective.
\hfill\qed 
\end{cor}

We define the isomorphism
$\iota^{\hor}:
H^0(\Sigma_{\theta},K_{\Sigma_{\theta}})
\to
H^1(X,\Def(E,\theta))^{\hor}$
by
\[
 \iota^{\hor}(\nu)
 =\Upsilon_0\circ
 \bigl(
\Upsilon_1\circ\Upsilon_0
 \bigr)^{-1}(\nu).
\]
We also regard $\iota^{\hor}$
as a morphism
$H^0(\Sigma_{\theta},K_{\Sigma_{\theta}})
\to
H^1(X,\Def(E,\theta))$
in a natural way.

\vspace{.1in}

Let us give an explicit description
of $(\Upsilon_1\circ\Upsilon_0)^{-1}(\nu)\in H^1(\gbiga^{\bullet})$
for $\nu\in H^0(\Sigma_{\theta},K_{\Sigma_{\theta}})$.
As in \S\ref{subsection;23.5.3.3},
we have
\[
 F_{\nu}\in H^0(X,\End(E)\otimes K_{\Sigma_{\theta}}(D)).
\]
By Lemma \ref{lem;23.5.3.2},
there exists a section
$g_P\in H^0(X_P,\End^{\asym}(E_P;C_P)(D))$
on $X_P$
such that
\[
\varphi_P:=F_{\nu|X_P}-\ad(\theta)(g_P)
\in H^0(X_P,\End^{\sym}(E_P;C_P)\otimes K_{X_P}). 
\]
We set
$\psi_P:=g_{P|X_P\cap X_0}$
and
$\varphi_0:=F_{\nu|X_0}$.
Then, 
we obtain a $1$-cocycle
$a=(\varphi_0,\varphi_P,\psi_P)$
in $\gbiga^1$.
The following lemma is clear by the construction.
\begin{lem}
The cohomology class
$[a]\in H^1(\gbiga^{\bullet})$ is
$(\Upsilon_1\circ\Upsilon_0)^{-1}(\nu)$.
Its image in $H^1(X,\Def(E,\theta))$
is $\iota^{\hor}(\nu)$.
\hfill\qed
\end{lem}

\begin{rem}
We shall study the harmonic representative of
$\iota^{\hor}(\nu)$ later.
\hfill\qed 
\end{rem}

\subsubsection{$C^{\infty}$-description of
horizontal infinitesimal deformations}

Let $\nu$ be a holomorphic $1$-form on $\Sigma_{\theta}$.
We obtain the following proposition
from Theorem \ref{thm;23.4.24.1}.
\begin{cor}
\label{cor;23.5.12.1}
Let $\tau\in A^1(\Def(E,\theta))$
be a $1$-cocycle such that
$\rho(\tau)=\nu$.
(See Lemma {\rm\ref{lem;23.5.3.4}} for
the map $\rho$.) 
We also assume the following conditions.
\begin{itemize}
 \item $\tau_{|X\setminus\bigcup X_P'}
       =F_{\nu|X\setminus\bigcup X_P'}$.
 \item $\tau_{|X_P}\in A^1(\Def(E_P,\theta_P;C_P))$
      for each $P\in D(\theta)$.       
\end{itemize}
Then, $[\tau]=\iota^{\hor}(\nu)$ holds
in $H^1(X,\Def(E,\theta))$.
In particular,
$[\tau]\in H^1(X,\Def(E,\theta))^{\hor}$.
Conversely,
any element of $H^1(X,\Def(E,\theta))^{\hor}$
is described in this way.
\hfill\qed
\end{cor}

\subsubsection{Lagrangian property of the horizontal subspaces}
\label{subsection;24.10.23.10}

For any Higgs bundle $(E,\theta)$ of rank $n$ of degree $0$ on $X$,
let $\omega_{\nbigm_{H},(E,\theta)}$
denote the alternative bilinear form on $A^1(\Def(E,\theta))$
defined as follows \cite[(6.5)]{Hitchin-self-duality}:
\begin{equation}
 \omega_{\nbigm_H,(E,\theta)}(\tau_1,\tau_2)
=-\int_X\Tr(\tau_1\wedge\tau_2).
\end{equation}
It induces the alternative bilinear form on $H^1(X,\Def(E,\theta))$,
which is also denoted by $\omega_{\nbigm_H,(E,\theta)}$.
According to Hitchin,
the family of the alternative bilinear forms
$\bigcup_{(E,\theta)\in\nbigm_H}\omega_{\nbigm_H,(E,\theta)}$
induces a holomorphic symplectic form $\omega_{\nbigm_H}$
of $\nbigm_H$.

\begin{cor}
Let $(E,\theta)\in\nbigm_H'$.
Then, the horizontal subspace
$H^1(X,\Def(E,\theta))^{\hor}$
of $H^1(X,\Def(E,\theta))$
is Lagrangian with respect to
$\omega_{\nbigm_H,(E,\theta)}$. 
\end{cor}
\pf
For any
$\iota^{\hor}(\nu_i)\in H^1(X,\Def(E,\theta))^{\hor}$ $(i=1,2)$,
there exist representatives
$\tau_i\in A^1(\Def(E,\theta))$
satisfying the conditions in Corollary \ref{cor;23.5.12.1}.
On $X\setminus\bigcup X_P'$,
we have $\tau_i=\tau_i^{1,0}$,
and hence
$\tau_1\wedge \tau_2=0$.
On $X_P$,
$\tau_{i}^{1,0}$
(resp. $\tau_i^{0,1}$)
are symmetric
(anti-symmetric)
with respect to $C_P$.
We recall that
for a symmetric $(n\times n)$-matrix $A$
and an anti-symmetric $(n\times n)$-matrix $B$,
we obtain $\tr(AB)=0$
because
\[
 \tr(AB)
 =\tr(BA)=\tr\bigl(
 \lefttop{t}(BA)
 \bigr)
 =\tr\bigl(
 \lefttop{t}A\lefttop{t}B
 \bigr)
=-\tr(AB).
\]
Hence, we obtain
$\tr(\tau_1\wedge\tau_2)
 =
\tr\bigl(
 \tau_{1}^{1,0}\wedge\tau_2^{0,1}
 \bigr)
+\tr\bigl(
 \tau_{1}^{0,1}\wedge\tau_2^{1,0}
 \bigr)
=0$
on $X_P$.
Hence, we obtain
$\omega_{\nbigm_H,(E,\theta)}(\iota^{\hor}(\nu_1),\iota^{\hor}(\nu_2))
=\omega_{\nbigm_H,(E,\theta)}(\tau_1,\tau_2)=0$.
\hfill\qed

\section{Comparison of the Hitchin metric and the semi-flat metric}

\subsection{The semi-flat metric}

\subsubsection{Vertical infinitesimal deformations}

We recall the following.
(See \S\ref{subsection;23.5.10.1}
for relative Higgs bundles and the induced infinitesimal deformation.)
\begin{df}
\mbox{{}}
\begin{itemize}
 \item 
A family of Higgs bundles
$(E_s,\theta_s)$ $(s\in S)$ is called vertical
if the spectral curves
$\Sigma_{\theta_s}\subset T^{\ast}X$ are independent of $s$.
\item
       A relative Higgs bundle $(\nbige,\theta_{\nbige})$
on $S\times X$ is called vertical
if the associated family
$(\nbige,\theta_{\nbige})_{|\{s\}\times X}$ $(s\in S)$
is vertical.
\item
     An infinitesimal deformation of $(E,\theta)$
is called vertical
if it is induced by a vertical relative Higgs bundle.
Let $H^1(X,\Def(E,\theta))^{\ver}\subset H^1(X,\Def(E,\theta))$
denote the subset of vertical infinitesimal deformations.
\hfill\qed
\end{itemize}
\end{df}

Suppose that $(E,\theta)$ is generically regular semisimple
and that the spectral curve $\Sigma_{\theta}$ is smooth.
The natural morphism
$\pi_{\ast}\nbigo_{\Sigma_{\theta}}
=\nbigh^0\bigl(\Def(E,\theta)\bigr)
 \lrarr
 \Def(E,\theta)$
induces a morphism
\begin{equation}
\label{eq;23.4.24.2}
 \iota^{\ver}:
 H^1(\Sigma_{\theta},\nbigo_{\Sigma_{\theta}})
 \lrarr
 H^1(X,\Def(E,\theta)).
\end{equation}

\begin{lem}
The morphism {\rm(\ref{eq;23.4.24.2})}
is injective. 
The image equals
$H^1(X,\Def(E,\theta))^{\ver}$.
\end{lem}
\pf
Because
$\pi_{\ast}\nbigo_{\Sigma_{\theta}}
=\nbigh^0(\Def(E,\theta))$,
we obtain the first claim.
Let $U\subset\cnum$ be a neighbourhood of $0$.
For any $a\in H^1(\Sigma_{\theta},\nbigo_{\Sigma_{\theta}})$,
there exists a holomorphic line bundle
$\nbigl$ on $U\times \Sigma_{\theta}$
which induces $a$ as the infinitesimal part.
Then, the second claim is clear.
\hfill\qed

\vspace{.1in}

Any $a\in H^1(\Sigma_{\theta},\nbigo_{\Sigma_{\theta}})$
is represented by a $(0,1)$-form $\nu$
such that the support of $\nu$ is compact
in $\Sigma_{\theta}\setminus \Sigma_{\theta|D(\theta)}$.
We obtain the $C^{\infty}$-section $F'_{\nu}$
of $\End(E)\otimes\Omega^{0,1}$ on $X^{\circ}$
as in \S\ref{subsection;23.4.11.3}
whose support is compact.
It naturally extends to
a $C^{\infty}$-section $F_{\nu}$
of $\End(E)\otimes\Omega^{0,1}$ on $X$.
Because $[\theta,F_{\nu}]=0$,
it is a $1$-cocycle in $A^1(\Def(E,\theta))$.
The image $\iota^{\ver}(a)\in H^1(X,\Def(E,\theta))$
is represented by $F_{\nu}$.
\begin{rem}
We shall study the harmonic representative of
$\iota^{\ver}(a)$ later.
\hfill\qed
\end{rem}

As we recalled in \S\ref{subsection;24.10.23.10},
there exists the complex symplectic structure
$\omega_{\nbigm_H}$ of $\nbigm_H$.
The following is well known.
\begin{lem}
$H^1(X,\Def(E,\theta))^{\ver}$
is a Lagrangian subspace
with respect to $\omega_{\nbigm_H,(E,\theta)}$. 
\end{lem}
\pf
For $\iota^{\ver}(\nu_i)\in H^1(X,\Def(E,\theta))^{\ver}$,
let $F_{\nu_i}$ be representatives as above.
Then, $F_{\nu_1}\wedge F_{\nu_2}=0$.
Hence,
$\omega_{\nbigm_H,(E,\theta)}(\iota^{\ver}(\nu_1),\iota^{\ver}(\nu_2))=0$.
\hfill\qed

\subsubsection{The induced perfect pairings}
\label{subsection;23.5.10.20}

There exists the following decomposition:
\[
 H^1(X,\Def(E,\theta))
=H^1(X,\Def(E,\theta))^{\ver}
\oplus
H^1(X,\Def(E,\theta))^{\hor}.
\]
Because 
$H^1(X,\Def(E,\theta))^{\ver}$ is Lagrangian,
we obtain the following perfect pairings
induced by $\omega_{\nbigm_H,(E,\theta)}$.
\begin{equation}
 \label{eq;23.5.3.10}
 H^1(X,\Def(E,\theta))^{\ver}
 \otimes
 H^1(X,\Def(E,\theta))^{\hor}
 \lrarr\cnum
\end{equation}

\begin{lem}
\label{lem;23.4.24.11}
Under the isomorphisms
 $\iota^{\ver}:
 H^1(\Sigma_{\theta},\nbigo_{\Sigma_{\theta}})
  \simeq
 H^1(X,\Def(E,\theta))^{\ver}$
and 
 $\Upsilon_1:
 H^1(X,\Def(E,\theta))^{\hor}\simeq
H^0(\Sigma_{\theta},K_{\Sigma_{\theta}})$,
the perfect pairing {\rm(\ref{eq;23.5.3.10})}
is identified with the perfect pairing {\rm(\ref{eq;23.5.3.20})}.
\end{lem}
\pf
There exists a $(0,1)$-form $\nu$ on $\Sigma_{\theta}$
such that
(i) it represents $a$,
(ii) the support of $\nu$ is compact in
$\Sigma_{\theta}\setminus \Sigma_{\theta|D(\theta)}$.
Then, $\iota^{\ver}(a)$ is represented by $F_{\nu}$.
For $b\in H^0(\Sigma_{\theta},K_{\Sigma_{\theta}})$,
there exists a closed $1$-form $\tau$
in $A^1(\Def(E,\theta))$
such that $\rho(\tau)=b$.
(See Lemma \ref{lem;23.5.3.4} for $\rho$.)
We have
\[
\omega_{\nbigm_H,(E,\theta)}(\iota^{\ver}(a),[\tau])
=-\int_X\Tr(F_{\nu}\cdot \tau)
=-\int_{\Sigma_{\theta}}a\wedge b.
\]
Thus, we obtain Lemma \ref{lem;23.4.24.11}.
\hfill\qed

\subsubsection{Hermitian metrics on $H^1(X,\Def(E,\theta))$}

Let $g_1$ and $g_2$ be the Hermitian metrics of
$H^1(\Sigma_{\theta},\nbigo_{\Sigma_{\theta}})$
and $H^0(\Sigma_{\theta},K_{\Sigma_{\theta}})$
as in \S\ref{subsection;23.4.24.10}.
There exists the decomposition
\[
 H^1(X,\Def(E,\theta))
=H^1(X,\Def(E,\theta))^{\ver}
\oplus
H^1(X,\Def(E,\theta))^{\hor}
\simeq
H^1(\Sigma_{\theta},\nbigo_{\Sigma_{\theta}})
\oplus
H^0(\Sigma_{\theta},K_{\Sigma_{\theta}}).
\]
Let $g_1\oplus g_2$ denote
the Hermitian metric of $H^1(X,\Def(E,\theta))$
induced by $g_i$ $(i=1,2)$.

\subsubsection{The semi-flat metric}

Note that $\Phi_H:\nbigm_H'\to \AAA_{X,n}'$
is an algebraic integrable system.
For each $s\in\AAA_{X,n}'$,
we obtain the spectral curve $\Sigma_s\subset T^{\ast}X$,
and $\Phi_H^{-1}(s)$ is identified
with the moduli space
$\Pic_d(\Sigma_s)$ of holomorphic line bundles
on $\Sigma_s$ of degree $d=n(n-1)(g_X-1)$.
By choosing $K_X^{1/2}$,
we obtain the Hitchin section $\Psi$ of $\Phi_H$.
By using the Hitchin section,
we may regard
$\Pic_d(\Sigma_s)$
as the quotient of
$H^1(\Sigma_s,\nbigo_{\Sigma_s})$
by a lattice.
The space
$H^1(\Sigma_s,\nbigo_{\Sigma_s})$
is equipped with the real symplectic form $\omega_s$
which induces a positive symmetric bilinear form
on $H^1(\Sigma_s,\nbigo_{\Sigma_s})$.
The cohomology class $[\omega_s]\in H^2(\Pic_d(\Sigma_s),\real)$
is a positive polarization in the sense of \cite[\S3]{Freed}.
According to \cite[Theorem 3.8]{Freed},
we obtain a hyperk\"ahler metric $g_{\semiflat,\real}$,
which is called the semi-flat metric.
Let $g_{\semiflat}$ denote the associated K\"ahler
metric with respect to the natural complex structure
of $\nbigm'_H$
as an open subset of the moduli space of Higgs bundles.

\begin{prop}
\label{prop;23.5.16.1}  
Under the identification
$T_{(E,\theta)}\nbigm'_H\simeq H^1(X,\End(E,\theta))$,
the metric $(g_{\semiflat})_{(E,\theta)}$
equals $g_1\oplus g_2$.
 \end{prop}
\pf
Let $\omega_1$ and $\omega_2$
be the real symplectic forms on
$H^1(\Sigma_{\theta},\nbigo_{\Sigma_{\theta}})$
and
$H^0(\Sigma_{\theta},K_{\Sigma_{\theta}})$
as in \S\ref{subsection;23.4.24.10}.
Let us recall the construction of
the semi-flat metric at $T_{(E,\theta)}\nbigm_H$
by following \cite[\S2 and \S3]{Freed}.

Let $\omega_{\nbigm_H}$ denote the complex symplectic
structure of $\nbigm_H$
in \cite[(6.5)]{Hitchin-self-duality}.
Under the isomorphism
$T_{(E,\theta)}\nbigm_H\simeq
H^1(X,\Def(E,\theta))$,
$(\omega_{\nbigm_H})_{(E,\theta)}$
equals the symplectic form $\omega_{\nbigm_H,(E,\theta)}$
in \S\ref{subsection;23.5.10.20}.
The tangent space
$T_{\Phi_H(E,\theta)}\AAA_{X,n}$
is identified with
the dual space of
the space of the isomorphism classes of
vertical infinitesimal deformations:
\begin{equation}
\label{eq;23.4.24.11}
 T_{\Phi_H(E,\theta)}\AAA_{X,n}
 \simeq
 H^1(X,\Def(E,\theta))^{\hor}
\simeq
 \bigl(H^1(X,\Def(E,\theta))^{\ver}\bigr)^{\lor}
=H^1(\Sigma_{\theta},\nbigo_{\Sigma_{\theta}})^{\lor}.
\end{equation}
We obtain the positive symplectic form
$\omega_{\Phi_H(E,\theta)}=-\omega_1^{\lor}$
on $T_{\Phi_H(E,\theta)}\AAA_{X,n}$.
In this way, we obtain the K\"ahler form
$\omega_{\AAA'_{X,n}}=\{\omega_s\}_{s\in \AAA'_{X,n}}$ on $\AAA'_{X,n}$.
The decomposition
$H^1(X,\Def(E,\theta))
=H^1(X,\Def(E,\theta))^{\ver}
\oplus
H^1(X,\Def(E,\theta))^{\hor}$
induces the flat connection on the tangent bundle
of $\AAA'_{X,n}$.
They give the special K\"ahler structure on $\AAA'_{X,n}$.
(See \cite[\S3]{Freed}.)
We obtain the special K\"ahler metric $g_{3}$ of $\AAA'_{X,n}$
induced by $\omega_{\AAA'_{X,n}}$.
Let $g_{3,\real}$ denote the real part of $g_3$.

Let $g_{3,\real|(E,\theta)}$ denote
the induced positive definite symmetric bilinear form of 
$H^1(X,\Def(E,\theta))^{\hor}
\simeq
T_{\Phi_H(E,\theta)}\AAA'_{X,n}$.
We obtain the positive definite symmetric bilinear form
$g_{3,\real|(E,\theta)}^{\lor}$
of $H^1(X,\Def(E,\theta))^{\ver}$
by the duality.
By definition of $g_{\semiflat}$,
we have
\[
(g_{\semiflat,\real})_{(E,\theta)}
=g_{3,\real|(E,\theta)}^{\lor}\oplus g_{3,\real|(E,\theta)}.
\]

Under the identification
$H^1(X,\Def(E,\theta))^{\hor}
\simeq
H^0(\Sigma_{\theta},K_{\Sigma_{\theta}})$,
$g_{3,\real|(E,\theta)}$
is the positive definite symmetric bilinear form associated with
$\omega_2=-\omega_1^{\lor}$
by Lemma \ref{lem;23.5.3.21},
we obtain
$g_{3,\real|(E,\theta)}=g_{2,\real}$.
Because $g_{1,\real}=g_{2,\real}^{\lor}$ by Lemma \ref{lem;23.5.3.21},
we obtain
$g_{3,\real|(E,\theta)}^{\lor}=g_{1,\real}$.
\hfill\qed

\subsection{The Hitchin metric and the comparison with the semi-flat metric}

\subsubsection{The Hitchin metric}
\label{subsection;23.5.16.41}

Let $(E,\theta)\in\nbigm_H$.
Let $h$ be a harmonic metric of $(E,\theta)$.
Let $\Harm^1(\End(E),\ad\theta,h)$ denote
the space of harmonic $1$-forms of
$(\End(E),\ad\theta,h)$.
There exist the following isomorphisms:
\[
 T_{(E,\theta)}\nbigm_H
 \simeq
 H^1(X,\Def(E,\theta))
 \simeq
 \Harm^1(\End(E),\ad\theta,h).
\]
For $\tau\in\Harm^1(\End(E),\ad\theta,h)$,
we set
\[
 g_{H,\real|(E,\theta)}(\tau,\tau)
 =2\sqrt{-1}\int_X
 \bigl(
 h\bigl(
 \tau^{1,0},
 \tau^{1,0}
 \bigr)
-h\bigl(
 \tau^{0,1},
 \tau^{0,1}
 \bigr)
 \bigr)
=2\sqrt{-1}\int_X
  \tr\bigl(
  \tau^{1,0}
  \wedge
  (\tau^{1,0})^{\dagger}_h
  -\tau^{0,1}\wedge
  (\tau^{0,1})^{\dagger}_{h}
  \bigr).
\]
We obtain
the positive definite symmetric product
$g_{H,\real|(E,\theta)}$
of the $\real$-vector space $T_{(E,\theta)}\nbigm_H$,
and the Hermitian metric $g_{H|(E,\theta)}$
of $T_{(E,\theta)}\nbigm_H$.
According to Hitchin,
they induce the hyperk\"ahler metric $g_{H,\real}$
and the K\"ahler metric $g_H$ of $\nbigm_H$,
respectively.
(See \S\ref{subsection;23.5.16.40} below
for a more precise comparison with the Riemannian metric constructed in
\cite{Hitchin-self-duality}.)
We have $g_{H,\real}=\Re g_H$.
The metrics $g_H$ and $g_{H,\real}$
are called the Hitchin metric.

\subsubsection{Main Theorem}

Let $s$ be the $C^{\infty}$-automorphism
of the complex vector bundle $T\nbigm_H'$
determined by $g_{H}=g_{\semiflat}\cdot s$.
Similarly,
let $s_{\real}$ be the $C^{\infty}$-automorphism
of the real vector bundle $T\nbigm_H'$
determined by
$g_{H,\real}=g_{\semiflat,\real}\cdot s_{\real}$.
We obtain the following theorem
from Proposition \ref{prop;23.4.12.3},
Proposition \ref{prop;23.5.5.230}
and Proposition \ref{prop;23.5.5.231}
below.
\begin{thm}
\label{thm;23.5.11.20}
For any compact subset $K\subset\nbigm_H'$,
there exist positive constants $B_i>0$ $(i=1,2)$
such that
$\bigl|(s-\id)_{(E,t\theta)}\bigr|_{g_{\semiflat}}
\leq B_1\exp(-B_2t)$
and
$\bigl|(s_{\real}-\id)_{(E,t\theta)}
\bigr|_{g_{\semiflat,\real}}
\leq B_1\exp(-B_2t)$
as $t\to\infty$ 
for any $(E,\theta)\in K$.
\end{thm}

Though the statements
in \S\ref{subsection;23.5.4.12}--\ref{subsection;23.5.3.31}
are given for a fixed Higgs bundle $(E,\theta)$,
the estimates can be given locally uniformly on $\nbigm'_H$.

\subsubsection{Appendix}
\label{subsection;23.5.16.40}

It is well known that
$g_{H,\real}$ in \S\ref{subsection;23.5.16.41}
equals the Riemannian metric of $\nbigm_H$
introduced in \cite[\S6]{Hitchin-self-duality}.
We include an explanation for the convenience of readers.
Let $(E,\theta,h)$ be a harmonic bundle.
Let $\nabla_h=\delbar_E+\del_{E,h}$ denote
the Chern connection of $(E,h)$.
Let $\theta^{\dagger}_h$ denote the adjoint of
$\theta$ with respect to $h$.
We obtain the following equations
for $\tau^{1,0}+\tau^{0,1}\in A^1(\End(E))$
as the linearization of the Hitchin equation
for a pair of a holomorphic structure and a Higgs field:
\begin{equation}
\label{eq;23.5.16.30}
  \del_{E,h}\tau^{0,1}+[\theta^{\dagger}_h,\tau^{1,0}]
   -\Bigl(
   \del_{E,h}\tau^{0,1}+[\theta^{\dagger}_h,\tau^{1,0}]
   \Bigr)^{\dagger}_h=0,
\end{equation}
\begin{equation}
 \label{eq;23.5.16.31}
  \delbar_E\tau^{1,0}+[\theta,\tau^{0,1}]=0.
\end{equation}
Let $\gminiu(E,h)\subset \End(E)$ denote the subbundle of
endomorphisms of $E$ which are anti-self-adjoint
with respect to $h$.
The linearization $d_1:A^0(\gminiu(E,h))\to A^1(\End(E))$
of the action of the gauge group
is given by
\[
 d_1(a)=[\theta,a]+\delbar_E(a).
\]
We consider the real symmetric bilinear forms
on $A^0(\gminiu(E,h))$ and $A^1(\End(E))$
induced by $h$ and a conformal metric $g_X$ of $X$.
The formal adjoint
$d_1^{\ast}:A^1(\End(E))\to A^0(\gminiu(E,h))$ is given by
\begin{equation}
\label{eq;23.5.16.32}
 d_1^{\ast}(\tau)=
 -\sqrt{-1}\Lambda_{g_X}\Bigl(
  \del_{E,h}\tau^{0,1}+[\theta^{\dagger}_h,\tau^{1,0}]
  +
  \bigl(
  \del_{E,h}\tau^{0,1}+[\theta^{\dagger}_h,\tau^{1,0}]
  \bigr)^{\dagger}_h
  \Bigr).
\end{equation}
The system of equations \cite[(6.1)]{Hitchin-self-duality}
is equivalent to
the system of equations
(\ref{eq;23.5.16.30}), (\ref{eq;23.5.16.31})
and $d_1^{\ast}\tau=0$.
It is equivalent to
the system of equations
(\ref{eq;23.5.16.31})
and $(\del_{E,h}+\ad\theta^{\dagger}_h)\tau=0$,
i.e.,
the condition for $\tau$ to be a harmonic $1$-form.
Because the (real) Hitchin metric at $(E,\theta)\in\nbigm_H$
is defined as the restriction of the $L^2$-metric
to the space of the solutions of
\cite[(6.1)]{Hitchin-self-duality},
it equals $g_{H,\real|(E,\theta)}$.

\subsection{Harmonic representatives of horizontal infinitesimal deformations}
\label{subsection;23.5.4.12}

Let $\iota^{\hor}_t:
H^0(\Sigma_{\theta},K_{\Sigma_{\theta}})
\lrarr
H^1(X,\Def(E,t\theta))$
denote the composite of the following morphisms:
\[
 H^0(\Sigma_{\theta},K_{\Sigma_{\theta}})
 \stackrel{c_1}{\simeq}
 H^0(\Sigma_{t\theta},K_{\Sigma_{t\theta}})
 \stackrel{c_2}{\lrarr}
 H^1(X,\Def(E,t\theta)).
\]
Here, $c_1$ is induced by
$\Sigma_{\theta}\simeq \Sigma_{t\theta}$,
and $c_2$ is the morphism $\iota^{\hor}$
in (\ref{subsection;23.5.10.4})
for $(E,t\theta)$.
Let $\nu\in H^0(\Sigma_{\theta},K_{\Sigma_{\theta}})$.
Let us study the harmonic $1$-form 
$\ttH(\nu,t)\in A^1(\Def(E,t\theta))$
such that
$[\ttH(\nu,t)]=\iota^{\hor}_t(\nu)$.
Let $\|\nu\|$ denote the $L^2$-norm of $\nu$ on $\Sigma_{\theta}$.
Let $g_X$ denote a conformal metric of $X$.

\subsubsection{Asymptotically harmonic representative on $X^{\circ}$}

We recall $X^{\circ}=X\setminus D(\theta)$.
For simplicity of the description,
$D(\theta)$ is also denoted by $D$.
We obtain
$F_{\nu}\in
H^0(X^{\circ},\End(E)\otimes K_X)$
as in \S\ref{subsection;23.4.11.3}.
Let $N(D)$ be a neighbourhood of $D$.

\begin{lem}
\label{lem;23.5.5.1}
There exist positive constants
$B_i$ $(i=0,1,2)$ such that the following holds
on $X\setminus N(D)$ for any $t\geq 1$:
\[
 |F_{\nu}|_{h_t,g_X}\leq B_0\|\nu\|,
 \quad\quad
 \bigl|
 (\delbar_E+\ad t\theta)^{\ast}_{h_t,g_X}
 F_{\nu}
 \bigr|_{h_t}
 \leq
 B_1\exp(-B_2t)\|\nu\|.
\]
\end{lem}
\pf
We obtain the first from \cite[Corollary 2.6]{Decouple},
and the second from \cite[Corollary 2.5]{Decouple}.
\hfill\qed

\begin{lem}
\label{lem;23.5.5.111}
There exist positive constants $B_i$ $(i=1,2)$
such that the following holds:
\begin{equation}
\label{eq;23.5.5.101}
\left|
 \sqrt{-1}
 \int_{X\setminus N(D)}
 h_{t}(F_{\nu},F_{\nu})
-\sqrt{-1}
 \int_{\Sigma_{\theta|X\setminus N(D)}}
 \nu\wedge\nubar
 \right|
 \\
 \leq B_1\exp(-B_2t)
 \|\nu\|^2.
\end{equation}
\end{lem}
\pf
We have
\[
\sqrt{-1}
\int_{\Sigma_{\theta|X\setminus N(D)}}
 \nu\wedge\nubar
 =\sqrt{-1}
 \int_{X\setminus N(D)}
 \tr\bigl(
  F_{\nu}\cdot F_{\nubar}
 \bigr).
\]
We also have
\[
 \sqrt{-1}
 \int_{X\setminus N(D)}
 h_{t}(F_{\nu},F_{\nu})
 =\sqrt{-1}\int_{X\setminus N(D)}
 \tr\bigl(
 F_{\nu}\cdot (F_{\nu})^{\dagger}_{h_t}
 \bigr).
\]
By \cite[Corollary 2.6]{Decouple},
there exist $B_i>0$ $(i=3,4)$
such that
\[
\left|
(F_{\nu}^{\dagger})_{h_t}
 -F_{\nubar}
 \right|_{h_t}
 \leq
 B_3\exp(-B_4t)\|\nu\|.
\]
Hence, we obtain (\ref{eq;23.5.5.101}).
\hfill\qed

\vspace{.1in}
We obtain the $1$-cocycle
$(\tr F_{\nu})\id_E\in A^1(\Def(E,\theta))$.
\begin{lem}
 $(\delbar_E+\ad t\theta)^{\ast}_{g_X}
 ((\tr F_{\nu})\id_E)=0$.
\hfill\qed
\end{lem}

\subsubsection{Asymptotically harmonic representatives around the discriminant}
\label{subsection;23.5.5.122}

Let $U$ be an open subset of $X$
with a holomorphic coordinate $z$
such that the following holds.
\begin{itemize}
 \item The coordinate induces the holomorphic isomorphism
       $U\simeq \{|z|<2\}$.
 \item $D(\theta)\cap U\subset\{|z|<1/4\}$.
 \item Each connected component of
       $\Sigma_{\theta|U}$ is simply connected.
\end{itemize}
We fix a base point of each connected component of
$\Sigma_{\theta|U}$.
We set
$(E_U,\theta_U)=(E,\theta)_{|U}$.
Note that there exists a non-degenerate symmetric product $C_U$ of
$(E_U,\theta_U)$ compatible with
the limiting configuration $h_{\infty}$.
Let $h_{U,t}$ be the harmonic metric of $(E_U,t\theta_U)$
compatible with $C_U$.
There exists a holomorphic function $\alpha_U$
on $\Sigma_{\theta|U}$
such that $d\alpha_U=\nu_{|\Sigma_{\theta|U}}$
and that
$\alpha_U$ is $0$ at the base points
of the connected components of $\Sigma_{\theta|U}$.
We obtain
$F_{\alpha_U}\in \nbigh^0(\Def(E_U,\theta_U;C_U))$
as in \S\ref{subsection;23.4.11.3}.
We set
\[
\ttH_U(\nu,t):=
(\del_{E_U,h_{U,t}}+\ad t(\theta_U)^{\dagger}_{h_{U,t}})
F_{\alpha_U}.
\]
By Proposition \ref{prop;23.5.4.10},
we obtain the following.
\begin{lem}
\label{lem;23.5.5.40}
$\ttH_U(\nu,t)$ is a harmonic $1$-form of
 $(\End(E_U),t\ad\theta_U,h_{U,t})$ on $U$,
i.e.,
\[
 (\delbar_{E_U}+\ad t\theta_{U})\ttH_U(\nu,t)=0,
 \quad
 (\delbar_{E_U}+\ad t\theta_{U})^{\ast}_{h_{U,t},g_X}\ttH_U(\nu,t)=0.
\] 
Also, $\ttH_U(\nu,t)$ is contained in
$A^1(\Def(E_U,\theta_U;C_U))$.
Moreover, we have
$\rho(\ttH_U(\nu,t))=\pi_{\ast}(\nu)_{|U}$,
which implies
$\tr(\ttH_U(\nu,t))=\tr\pi_{\ast}(\nu)_{|U}$.
(See {\rm\S\ref{subsection;24.10.22.1}}
for $\tr\pi_{\ast}(\nu)$.)
\hfill\qed
\end{lem}

For $0<r\leq 2$,
we set $U(r):=\{|z|<r\}$.

\begin{lem}
\label{lem;23.5.5.2}
For any $1/4<r_2<r_1<2$,
there exist positive constants $B_i>0$ $(i=0,1,2,3)$
such that the following holds
on $U(r_1)\setminus U(r_2)$ for any $t\geq 1$:
\[
 |F_{\alpha_U}|_{h_{U,t}}\leq B_0\|\nu\|,
 \quad\quad
 |\ttH_U(\nu,t)-F_{\nu}|_{h_{U,t}}\leq B_1e^{-B_2 t}
 \|\nu\|,
 \quad\quad
 |\ttH_U(\nu,t)|_{h_{U,t}}
 \leq B_3\|\nu\|.
\]
\end{lem}
\pf
We obtain the first from \cite[Corollary 2.6]{Decouple},
the second from
\cite[Corollary 2.5, Proposition 2.10]{Decouple},
and the third from Lemma \ref{lem;23.5.5.1}.
\hfill\qed

\begin{lem}
\label{lem;26.1.28.1}
There exists a unique element
$\sigma_{U}(\nu,t)\in A^0(U\cap X^{\circ},\End^{\asym}(E;C_U))$
such that
\[
 (\delbar_{E_U}+\ad t\theta_U)\sigma_{U}(\nu,t)
 =\ttH_U(\nu,t)-F_{\nu}.
\]
For any $1/4<r_2<r_1<2$,
there exist $B_i>0$ $(i=1,2)$ such that
the following holds on $U(r_1)\setminus U(r_2)$
such that 
\[
 |\sigma_{U}(\nu,t)|_{h_{U,t}}
 \leq B_1e^{-B_2 t}\|\nu\|.
\]
\end{lem}
\pf
We obtain the first claim from Proposition \ref{prop;23.4.15.20}.
The second follows from Lemma \ref{lem;23.5.5.2}.
\hfill\qed

\begin{rem}
Constants $B_i$ of the lemmas in this section
are not necessarily the same.
\hfill\qed 
\end{rem}

For any $0<r<2$,
let $\|\ttH_U(\nu,t)_{|U(r)}\|_{L^2,h_{U,t}}$
denote the $L^2$-norm of $\ttH_U(\nu,t)$
on $U(r)$
with respect to $h_{U,t}$.
Note that it is independent of the choice of
a conformal metric of the Riemann surface $X$.

\begin{lem}
\label{lem;23.4.12.10}
For any $0<r<2$,
there exists $B>0$ such that
$\|\ttH_U(\nu,t)_{|U(r)}\|_{L^2,h_{U,t}}
\leq B\|\nu\|$
for any $t\geq 1$.
\end{lem}
\pf
We note the following equality:
\[
\ttH_U(\nu,t)=
(\del_{E_U,h_{U,t}}+\ad t(\theta_U)^{\dagger}_{h_{U,t}})
F_{\alpha_U}
=(\delbar_{E_U}+\del_{E_U,h_{U,t}}
 +\ad t\theta_U+\ad t(\theta_U)^{\dagger}_{h_{U,t}})
F_{\alpha_U}.
\]
By Proposition \ref{prop;23.5.4.20},
we obtain the following equality:
\begin{multline}
\label{eq;24.10.22.2}
\sqrt{-1}
 \int_{U(r)}\Bigl(
 h_{U,t}(\ttH_U(\nu,t)^{1,0},\ttH_U(\nu,t)^{1,0})
-h_{U,t}(\ttH_U(\nu,t)^{0,1},\ttH_U(\nu,t)^{0,1})
 \Bigr)
 =\\
 \sqrt{-1}\int_{\del U(r)}
 \Bigl(
 h_{U,t}(F_{\alpha_U},\ttH_U(\nu,t)^{1,0})
-h_{U,t}(F_{\alpha_U},\ttH_U(\nu,t)^{0,1})
 \Bigr).
\end{multline}
We obtain the claim of Lemma \ref{lem;23.4.12.10}
from Lemma \ref{lem;23.5.5.2}.
\hfill\qed

\begin{lem}
\label{lem;23.5.5.10}
For any $0<r<2$,
there exists $B'>0$ such that
$|\ttH_U(\nu,t)|_{h_{U,t},g_X}\leq B't\|\nu\|$ on $U(r)$
for any $t\geq 1$.
\end{lem}
\pf
We set $Y(t)=\{|w|<2t\}\subset\cnum$.
We define the isomorphism
$\varphi_t:Y(t)\to U$ by $\varphi_t(w)=t^{-1}w$.
We obtain the harmonic bundles
$(\Etilde_{U,t},\thetatilde_{U,t},\htilde_{U,t})=
\varphi_t^{\ast}(E_U,t\theta_U,h_{U,t})$
on $Y(t)$.
Let $f_{U,t}$ be the endomorphism of
$\Etilde_{U,t}$
defined by
$\thetatilde_{U,t}=f_{U,t}\,dw$.
Then, there exists $B_1>0$ such that 
any eigenvalue $\alpha$ of $f_{U,t|Q}$ 
satisfies $|\alpha|\leq B_1$
for any $t\geq 1$ and any $Q\in Y(t)$.
Let $0<r_1<r_2<1$.
By Lemma \ref{lem;23.4.12.10},
we have
\[
 \|\varphi_t^{-1}\ttH_U(\nu,t)_{|Y(r_2t)}\|_{L^2,\htilde_{U,t}}
 \leq
 B\|\nu\|.
\]
Let $g_0=dw\,d\wbar$ be the Euclidean metric of $Y(t)$.
By Proposition \ref{prop;23.5.5.31},
there exists $B_2>0$
such that
\[
|\varphi_{t}^{-1}\ttH_U(\nu,t)|_{\htilde_{U,t},g_0}
\leq B_2\|\nu\| 
\]
on $Y(r_1t)$.
It implies the claim of the lemma.
\hfill\qed

\begin{cor}
\label{cor;23.5.5.50}
There exist $B_1,B_2>0$ such that
the following holds for any $t\geq 1$:
\begin{equation}
 \label{eq;23.5.5.42}
\bigl|
(\delbar_{E_U}+\ad t\theta_U)^{\ast}_{h_t,g_X}\ttH_{U,t}
\bigr|_{h_t}
 \leq
 B_1\exp(-B_2t)\|\nu\|.
\end{equation}
\end{cor}
\pf
As remarked in Lemma \ref{lem;23.5.5.40},
we have
$(\delbar_{E_U}+\ad t\theta_U)^{\ast}_{h_{U,t},g_X}\ttH_{U,t}
=0$.
By the estimate
$h_{U,t}-h_t=O(e^{-\beta t})$
with respect to $h_{U,t}$ for some $\beta>0$
in \cite[Theorem 7.17]{Mochizuki-Szabo},
we obtain (\ref{eq;23.5.5.42})
from Lemma \ref{lem;23.5.5.10}.
\hfill\qed

\begin{lem}
\label{lem;23.5.5.112}
There exist $B_i>0$  $(i=1,2)$ 
such that the following holds.
\begin{multline}
\label{eq;23.5.5.100}
 \left|
 \sqrt{-1}
 \int_{U(1)}
 \Bigl(
 h_{t}(\ttH_U(\nu,t)^{1,0},\ttH_U(\nu,t)^{1,0})
-h_{t}(\ttH_U(\nu,t)^{0,1},\ttH_U(\nu,t)^{0,1})
 \Bigr)
-\sqrt{-1}
 \int_{\Sigma_{\theta|U(1)}}
 \nu\wedge\nubar
 \right|
 \\
 \leq B_1\exp(-B_2t)
 \|\nu\|^2.
\end{multline}
\end{lem}
\pf
By \cite[Theorem 7.17]{Mochizuki-Szabo},
there exists $\beta>0$
such that
$h_{U,t}-h_t=O(e^{-\beta t})$ on $U(1)$ as $t\to\infty$.
Hence, by Lemma \ref{lem;23.4.12.10},
there exist $B_i>0$ $(i=3,4)$
such that the following holds:
\[
 \left|
 \sqrt{-1}\int_{U(1)}
 h_{t}\bigl(\ttH_U(\nu,t)^{p,q},\ttH_U(\nu,t)^{p,q}\bigr)
-\sqrt{-1}\int_{U(1)}
 h_{U,t}\bigl(\ttH_U(\nu,t)^{p,q},\ttH_U(\nu,t)^{p,q}\bigr)
 \right|
 \\
 \leq B_3\exp(-B_4t)
 \|\nu\|^2.
\]
Recall (\ref{eq;24.10.22.2}).
By Lemma \ref{lem;23.5.5.2},
there exist $B_i>0$ $(i=5,6)$ such that
\[
\left|
 \sqrt{-1}\int_{\del U(1)}
 \Bigl(
h_{U,t}(F_{\alpha_U},\ttH_U(\nu,t)^{0,1})
 \Bigr)
 \right|\leq B_5\exp(-B_6t)\|\nu\|^2.
\]
We have
\[
 \sqrt{-1}\int_{\del U(1)}
 h_{U,t}(F_{\alpha_U},\ttH_U(\nu,t)^{1,0})
=\sqrt{-1}\int_{\del U(1)}
 \tr\bigl(
 F_{\alpha_U}
 \cdot
 (\ttH_U(\nu,t)^{1,0})^{\dagger}_{h_{U,t}}
 \bigr).
\]
We also have
\[
 \sqrt{-1}\int_{\Sigma_{\theta|U(1)}}\nu\wedge\overline{\nu}
 =\sqrt{-1}\int_{\del \Sigma_{\theta|U(1)}}
 \alpha_U\wedge \overline{\nu}
 =\sqrt{-1}\int_{\del U(1)}
 \tr\bigl(
 F_{\alpha_U}\cdot F_{\nubar}
 \bigr).
\]
By \cite[Corollary 2.6]{Decouple}
and Lemma \ref{lem;23.5.5.2},
there exist $B_i>0$ $(i=7,8)$
such that the following holds on a neighbourhood of
$\del U(1)$:
\[
 \bigl|
 (\ttH_U(\nu,t)^{1,0})^{\dagger}_{h_{U,t}}
-F_{\nubar}
 \bigr|_{h_{U,t}}
 \leq B_7\exp(-B_8t)\|\nu\|.
\]
Then, we obtain (\ref{eq;23.5.5.100}).
\hfill\qed

\subsubsection{Asymptotically harmonic representatives on $X$}
\label{subsection;23.5.8.1}

Let $U_j$ $(j=1,\ldots,m)$
be open subsets of $X$
equipped with a holomorphic coordinate $z_j$
satisfying the conditions in \S\ref{subsection;23.5.5.122},
such that
$D(\theta)\subset \bigcup_{j=1}^mU_j$
and that
$U_j\cap U_{k}=\emptyset$ $(j\neq k)$.
We set $r_1=1$ and $r_2=\frac{1}{2}$.
Let $\chi_j:X\to[0,1]$ be a $C^{\infty}$-function
such that
$\chi_j=1$ on $U_j(r_2)$
and
$\chi_j=0$ on $X\setminus U_j(r_1)$.
There exists a closed $1$-form
$\ttH'(\nu,t)\in A^1(\Def(E,t\theta))$
such that the following holds.
\begin{itemize}
 \item On $X\setminus\bigcup_{j=1}^m U_j(r_1)$,
       we have
       $\ttH'(\nu,t)=F_{\nu}$.
 \item On $U_j(r_1)\setminus \Ubar_j(r_2)$,
we have
       $\ttH'(\nu,t)=
       F_{\nu}
       +(\delbar_E+\ad t\theta)
       \bigl(
       \chi_j\cdot \sigma_{U_j}(\nu,t)
       \bigr)$.
 \item On $U_j(r_2)$, we have
       $\ttH'(\nu,t)
       =\ttH_{U_j}(\nu,t)$.
\end{itemize}

\begin{lem}
\mbox{{}}\label{lem;23.5.5.51}
\begin{itemize}
 \item $(\delbar_E+\ad t\theta)\ttH'(\nu,t)=0$,
       and
       $\rho_t(\ttH'(\nu,t))=\pi_{\ast}(\nu)$.
       In particular,
       $\tr(\ttH'(\nu,t))=\tr\pi_{\ast}(\nu)$.
       Here, $\rho_t$ denotes the morphism $\rho$
       (see Lemma {\rm\ref{lem;23.5.3.4}})
       for $(E,t\theta)$.
 \item There exist $B_i>0$ $(i=1,2)$
       such that
       $|\ttH'(\nu,t)-F_{\nu}|_{h_t,g_X}
       \leq
       B_1\exp(-B_2t)$
       on $X\setminus\bigcup U_j(r_2)$
       and that
       $|\ttH'(\nu,t)-\ttH_{U_j}(\nu,t)|_{h_t,g_X}
       \leq
       B_1\exp(-B_2t)$
       on $U_j(r_1)$.
\item
       There exists $C>0$ such that $|\ttH'(\nu,t)|_{h_t,g_X}\leq Ct\|\nu\|$
       and $\|\ttH'(\nu,t)\|_{L^2,h_t}\leq C\|\nu\|$
       for any $t\geq 1$.       
 \item There exist $B_i>0$ $(i=3,4)$ such that 
      $\bigl|
      (\delbar_E+\ad t\theta)_{h_t,g_X}^{\ast}\ttH'(\nu,t)
      \bigr|_{h_t}
       \leq B_3\exp(-B_4t)$ on $X$
       for any $t\geq 1$.
 \item $H'(\nu,t)$ is horizontal,
       i.e.,
       $[H'(\nu,t)]=\iota^{\hor}_t(\nu)$.
\end{itemize}
\end{lem}
\pf
We obtain the first claim by the construction.
We obtain the second claim by the construction
and the estimate in Lemma \ref{lem;26.1.28.1}.
We obtain the third claim
from the estimates for $F_{\nu}$ and $\ttH_{\nu}$.
The fourth claim follows from
Lemma \ref{lem;23.5.5.1}
and Corollary \ref{cor;23.5.5.50}.
We obtain the last claim
by Corollary \ref{cor;23.5.12.1}
and the construction of $H'(\nu,t)$.
\hfill\qed

\subsubsection{Harmonic representatives}

We note that
$\tr\Bigl(
(\delbar_E+\ad t\theta)^{\ast}(\ttH'(\nu,t))
\Bigr)=0$.
There exists a unique $\gamma(\nu,t)\in A^0(\End(E))$
such that $\tr\gamma(\nu,t)=0$ and
\[
 (\delbar_E+\ad t\theta)^{\ast}_{h_t,g_X}
 (\delbar_E+\ad t\theta)\gamma(\nu,t)
 =(\delbar_E+\ad t\theta)^{\ast}_{h_t,g_X}\ttH'(\nu,t).
\]

\begin{lem}
\label{lem;23.5.5.80}
There exist $B_i>0$ $(i=1,2)$ such that
the following holds on $X$ for any $t\geq 1$:
\begin{equation}
\label{eq;23.5.5.70}
 |\gamma(\nu,t)|_{h_t}
+|(\delbar_E+\ad t\theta)\gamma(\nu,t)|_{h_t,g_X}
 \leq B_1\exp(-B_2t)\|\nu\|.
\end{equation}
\end{lem}
\pf
By Proposition \ref{prop;23.4.12.1}
and Lemma \ref{lem;23.5.5.51},
there exist $B_{i}>0$ $(i=11,12)$
such that the following holds for any $t\geq 1$:
\[
 \|\gamma(\nu,t)\|_{L^2,h_t,g_X}
\leq B_{11}\exp(-B_{12}t)\|\nu\|.
\]

Let $P$ be any point of $X$,
which is not necessarily contained in $D(\theta)$.
Let $(X_P,z_P)$ be a holomorphic coordinate neighbourhood
around $P$ such that $z(P)=0$
and that $X_P\simeq\{|z_P|<2\}$.
We set $Y(t)=\bigl\{w\in\cnum\,\big|\,|w|<2t\bigr\}$.
For any $t\geq 1$, let $\varphi_t:Y(t)\to X_P$
be defined by $\varphi_t(w)=t^{-1}w$.
We obtain the harmonic bundle
$\varphi_t^{\ast}(E,t\theta,h_t)$ on $Y(t)$.
By applying Proposition \ref{prop;23.5.5.61}
to the section $\varphi_t^{\ast}\gamma(\nu,t)$,
we obtain (\ref{eq;23.5.5.70}).
\hfill\qed

\vspace{.1in}

We set
\[
 \ttH(\nu,t)
 =\ttH'(\nu,t)-(\delbar_E+\ad t\theta) \gamma(\nu,t).
\]
Then, $\ttH(\nu,t)$ is a harmonic $1$-form
satisfying
$[\ttH(\nu,t)]=[\ttH'(\nu,t)]=\iota^{\hor}(\nu)$.

\subsubsection{Estimate of the pairings}

\begin{prop}
\label{prop;23.4.12.3}
There exist $B_i>0$ $(i=1,2)$ such that
\begin{equation}
\label{eq;23.5.16.20}
 \Bigl|
 g_{H}(\iota^{\hor}_t(\nu),\iota^{\hor}_t(\nu))
 -g_{\semiflat}(\iota^{\hor}_t(\nu),\iota^{\hor}_t(\nu))
 \Bigr|
  \leq B_1\exp(-B_2t)\|\nu\|^2.
\end{equation}
\end{prop}
\pf
The estimate (\ref{eq;23.5.16.20}) is equivalent to
\begin{equation}
\label{eq;23.5.5.110}
\left|
 2\sqrt{-1}\int_X
\Bigl(
 h_{t}(\ttH(\nu,t)^{1,0},\ttH(\nu,t)^{1,0})
-h_{t}(\ttH(\nu,t)^{0,1},\ttH(\nu,t)^{0,1})
 \Bigr)
-2\sqrt{-1}\int_{\Sigma_{\theta}}
 \nu \overline{\nu}
\right| 
 \leq B_1\exp(-B_2t)\|\nu\|^2.
 \end{equation}
We set $X_1=X\setminus\bigcup_{j=1}^mU_j(1)$.
By Lemma \ref{lem;23.5.5.80},
there exist $B_i>0$ $(i=3,4)$
such that the following holds:
\[
\sup_{X_1}
 \bigl|\ttH(\nu,t)-F_{\nu}\bigr|_{h_t,g_X}
+
\sum_{j=1}^m\sup_{U_j(1)}
\bigl|\ttH(\nu,t)-\ttH_{U_j}(\nu,t)\bigr|_{h_t,g_X}
 \leq
 B_3\exp(-B_4t)\|\nu\|.
\]
Then, we obtain (\ref{eq;23.5.5.110})
from Lemma \ref{lem;23.5.5.111}
and Lemma \ref{lem;23.5.5.112}.
\hfill\qed

\subsection{Harmonic representatives of vertical infinitesimal deformations}
\label{subsection;23.5.3.31}

Let $\iota^{\ver}_t:
H^1(\Sigma_{\theta},\nbigo_{\Sigma_{\theta}})
\to
H^1(X,\Def(E,t\theta))$
denote the composite of the following morphisms:
\[
 H^1(\Sigma_{\theta},\nbigo_{\Sigma_{\theta}})
 \stackrel{c_1}{\simeq}
 H^1(\Sigma_{t\theta},\nbigo_{\Sigma_{t\theta}})
 \stackrel{c_2}{\lrarr}
 H^1(X,\Def(E,t\theta)).
\]
Here, $c_1$ is induced by the isomorphism
$\Sigma_{\theta}\simeq\Sigma_{t\theta}$,
and
$c_2$ is induced by
$\pi_{t\ast}\nbigo_{\Sigma_{t\theta}}
=\nbigh^0(\Def(E,t\theta))$.

Let $\tau$ be a harmonic $(0,1)$-form on $\Sigma_{\theta}$,
i.e., $\tau$ is a $(0,1)$-form such that
$\del_{\Sigma_{\theta}}\tau=0$.
It induces $[\tau]\in H^1(\Sigma_{\theta},\nbigo_{\Sigma_{\theta}})$.
Conversely, any element of $H^1(\Sigma_{\theta},\nbigo_{\Sigma_{\theta}})$
can be expressed uniquely in this way.
We shall study the harmonic representative of
$\ttV(\tau,t)$ of $\iota_t^{\ver}([\tau])$.
Let $\|\tau\|$ denote the $L^2$-norm of $\tau$
on $\Sigma_{\theta}$.

\subsubsection{Preliminary}

In general,
for any complex manifold $Y$,
let $Y^{\dagger}$ denote the conjugate of $Y$.
There exists the natural isomorphism
$T^{\ast}(X^{\dagger})\simeq (T^{\ast}X)^{\dagger}$.

We have the holomorphic vector bundles
$E^{\dagger}_t=(E,\del_{E,h_t})$
with the Higgs field $\theta_{h_t}^{\dagger}$ on $X^{\dagger}$.
The spectral curve of $(E,t\theta_{h_t}^{\dagger})$
equals $(\Sigma_{t\theta})^{\dagger}\subset
(T^{\ast}X)^{\dagger}=T^{\ast}(X^{\dagger})$.
Let $\pi^{\dagger}_t:(\Sigma_{t\theta})^{\dagger}\to X^{\dagger}$
denote the projection.
There exists a holomorphic line bundle
$L^{\dagger}_t$ on $(\Sigma_{t\theta})^{\dagger}$
with an isomorphism
$(\pi_t^{\dagger})_{\ast}L^{\dagger}_t\simeq E^{\dagger}_t$.

There exist the following natural morphisms:
\[
 \tr:
(\pi^{\dagger}_t)_{\ast}
 \nbigo_{\Sigma_{t\theta}^{\dagger}}
=(\pi^{\dagger}_1)_{\ast}
 \nbigo_{\Sigma_{\theta}^{\dagger}}
 \lrarr
 \nbigo_{X^{\dagger}},
\quad\quad
 \tr:
(\pi^{\dagger}_t)_{\ast}
 K_{\Sigma_{t\theta}^{\dagger}}
=(\pi^{\dagger}_1)_{\ast}
 K_{\Sigma_{\theta}^{\dagger}}
 \lrarr
 K_{X^{\dagger}}.
\]

\subsubsection{Asymptotically harmonic representatives on $X^{\circ}$}

We obtain $F_{\tau}\in A^{0,1}(X^{\circ},\End(E))$
as in \S\ref{subsection;23.4.11.3}.
For simplicity of the description,
$D(\theta)$ is also denoted by $D$.
Let $N(D)$ be a neighbourhood of $D$.

\begin{lem}
\label{lem;23.5.5.132}
There exist positive constants
$B_i$ $(i=0,1,2)$ such that the following holds
on $X\setminus N(D)$ for any $t\geq 1$:
\[
 |F_{\tau}|_{h_t,g_X}\leq B_0\|\tau\|,
 \quad\quad
 \bigl|
 (\delbar_E+\ad t\theta)^{\ast}_{h_t,g_X}
 F_{\tau}
 \bigr|_{h_t}
 \leq
 B_1\exp(-B_2t)\|\tau\|.
\]
\end{lem}
\pf
We obtain the first from \cite[Corollary 2.6]{Decouple},
and the second from \cite[Proposition 2.10]{Decouple}.
\hfill\qed

\vspace{.1in}

The following lemma is obvious by the construction.
\begin{lem}
We have
$\tr F_{\tau}=\tr(\pi^{\dagger}_1)_{\ast}(\tau)$.
In particular,
$(\delbar_E+\ad t\theta)^{\ast}_{h_t,g_X}\tr F_{\tau}=0$.
\hfill\qed
\end{lem}

\subsubsection{Asymptotically harmonic $1$-forms
around the discriminant}

Let $(U,z)$ be as in \S\ref{subsection;23.5.5.122}.
We choose a base point of
each connected component of $\Sigma_{\theta|U}$.
There exists an anti-holomorphic function $\beta_U$
such that
$\delbar\beta_U=\tau$ on $\Sigma_{\theta|U}$
and that $\beta_U=0$ at the base points
of the connected components.

We may regard $\beta_U$ as a holomorphic function on
$\Sigma_{\theta|X_P}^{\dagger}$.
We obtain the endomorphism
$F_{(E^{\dagger}_t,t\theta^{\dagger}_{h_t}),\beta_U}$
of $E^{\dagger}_t$.
It satisfies
\[
 (\del_{E,h_t}+\ad t\theta^{\dagger}_{h_t})
 F_{(E^{\dagger}_t,t\theta^{\dagger}_{h_t}),\beta_U}=0.
\]
As in Proposition \ref{prop;23.5.4.10},
we obtain the following harmonic $1$-form
of $(E_t,t\theta^{\dagger}_{h_t},h_t)$ on $U^{\dagger}$:
\[
 \ttV_U(\tau,t)
:=(\delbar_E+\ad t\theta)F_{(E^{\dagger}_t,t\theta^{\dagger}_{h_t}),\beta_U}.
\]
By Corollary \ref{cor;23.5.5.131},
$\ttV_U(\tau,t)$ is also a harmonic $1$-form of $(E,t\theta,h_t)$
on $U$.
In particular,
the following holds on $U$:
\[
 (\delbar_E+\ad t\theta)
 \ttV_U(\tau,t)=0,
 \quad\quad
 (\delbar_E+\ad t\theta)^{\ast}_{h_t,g_X}
 \ttV_U(\tau,t)=0.
\]

\begin{lem}
\label{lem;23.5.5.140}
For any $1/4<r_2<r_1<2$,
there exist $B_i>0$ $(i=0,1,2,3)$ such that
the following holds on
$U(r_1)\setminus U(r_2)$: 
\[
 \bigl|
 F_{(E^{\dagger}_t,t\theta^{\dagger}_{h_t}),\beta_U}
 \bigr|_{h_t}
\leq B_0\|\tau\|,\quad
\bigl|
 \ttV_U(\tau,t)-F_{\tau}
\bigr|_{h_t,g_X}
\leq B_1\exp(-B_2t)\|\tau\|,
 \quad
\bigl|
 \ttV_U(\tau,t)
\bigr|_{h_t,g_X}
\leq B_3\|\tau\|.
\] 
\end{lem}
\pf
We obtain the first from
\cite[Corollary 2.6]{Decouple},
the second from 
\cite[Corollary 2.5, Corollary 2.6, Proposition 2.10]{Decouple},
and the third from Lemma \ref{lem;23.5.5.132}.
\hfill\qed

\vspace{.1in}
For any $0<r<2$,
let $\|\ttV_U(\tau,t)_{|U(r)}\|_{L^2,h_{t}}$
denote the $L^2$-norm of $\ttV_U(\tau,t)$
on $U(r)$ with respect to $h_{t}$.
We obtain the following lemma
by using
Proposition \ref{prop;23.5.4.20}
and Lemma \ref{lem;23.5.5.140},
as in the proof of Lemma \ref{lem;23.4.12.10}.

\begin{lem}
For any $0<r<2$,
there exists $B>0$ such that
$\|\ttV_U(\tau,t)_{|U(r)}\|_{L^2,h_{t}}
\leq B\|\tau\|$
for any $t\geq 1$.
\hfill\qed
\end{lem}

The following lemma is similar to
Lemma \ref{lem;23.5.5.10}.
\begin{lem}
For any $0<r<2$,
there exists $B'>0$ such that
$|\ttV_U(\tau,t)|_{h_{t},g_X}\leq B't\|\tau\|$ on $U(r)$
for any $t\geq 1$.
\hfill\qed
\end{lem}

By Proposition \ref{prop;23.5.4.10},
we obtain the following lemma.
\begin{lem}
$\tr(\ttV_U(\tau,t))
=\tr (\pi_{1}^{\dagger})_{\ast}(\tau)$.
\hfill\qed
\end{lem}

\begin{rem}
From the holomorphic function $\betabar_U$ on $\Sigma_{\theta|U}$,
we obtain the holomorphic endomorphism $F_{\betabar_U}$
of $(E,t\theta)_{|U}$.
We have 
$(F_{\betabar_U})^{\dagger}_{h_t}=
F_{(E^{\dagger}_t,t\theta^{\dagger}_t),\beta_U}$.
\hfill\qed 
\end{rem}

\subsubsection{Asymptotically harmonic representatives on $X$}
\label{subsection;23.5.12.20}

Let $(U_j,z_j)$ and $\chi_j$
$(j=1,\ldots,m)$ be as in \S\ref{subsection;23.5.8.1}.
Let $\tau'$ be the $(0,1)$-form on $\Sigma_{\theta}$
defined as follows:
\[
\tau':=
 \tau-\sum_{j}
 \delbar\bigl(
  \pi^{\ast}(\chi_j)\cdot\beta_{U_j}
 \bigr).
\]
We obtain a $1$-cocycle
$F_{\tau'}\in A^1(\Def(E,t\theta))$
such that
$[F_{\tau'}]=\iota^{\ver}_t([\tau])$.
We set
\[
 \ttV'(\tau,t)
=F_{\tau'}
+\sum_{j=1}^m
(\delbar_E+\ad t\theta)
\Bigl(
 \chi_j\cdot F_{(E^{\dagger}_{h_t},t\theta^{\dagger}_{h_t}),\beta_{U_j}}
\Bigr).
\]
By the construction, the following holds.
\begin{lem}
\mbox{{}}\label{lem;23.5.11.2}
\begin{itemize}
 \item $(\delbar_E+\ad t\theta)\ttV'(\tau,t)=0$
and $[\ttV'(\tau,t)]=\iota^{\ver}_{t}([\tau])$.
 \item We have
       $\tr \ttV'(\tau,t)=\tr (\pi^{\dagger}_1)_{\ast}(\tau)$.
       In particular,
       $(\delbar_E+\ad t\theta)^{\ast}_{h_t,g_X}\tr\ttV'(\tau,t)=0$.
 \item There exist $B_i>0$ $(i=1,2)$
       such that
       $|\ttV'(\tau,t)-F_{\tau}|_{h_t,g_X}
       \leq
       B_1\exp(-B_2t)\|\tau\|$
       on $X\setminus\bigcup U_j(r_2)$
       and that
       $|\ttV'(\tau,t)-\ttV_{U_j}(\tau,t)|_{h_t,g_X}
       \leq
       B_1\exp(-B_2t)\|\tau\|$
       on $U_j(r_1)$.
\item
There exists $B>0$,
which is independent of $t$,
such that $|\ttV'(\tau,t)|_{h_t,g_X}\leq Bt\|\tau\|$
      and that $\|\ttV'(\tau,t)\|_{L^2,h_t}\leq B\|\tau\|$.
\item There exist $B_i>0$ $(i=1,2)$ such that
      $\bigl|
      (\delbar_E+\ad t\theta)^{\ast}_{h_t,g_X}
      \ttV'(\tau,t)\bigr|_{h_t}
      \leq B_1\exp(-B_2 t)\|\tau\|$.
     \hfill\qed
\end{itemize} 
\end{lem}

\subsubsection{Harmonic representatives}

There uniquely exists $\gamma(\tau,t)\in A^0(\End(E))$
satisfying
$\tr\gamma(\tau,t)=0$ and
\[
 (\delbar_E+\ad t\theta)^{\ast}_{h_t,g_X}
 (\delbar_E+\ad t\theta)\gamma(\tau,t)
 =(\delbar_E+\ad t\theta)^{\ast}_{h_t,g_X}
 \ttV'(\tau,t).
\]
The following lemma is similar to Lemma \ref{lem;23.5.5.80}.
\begin{lem}
\label{lem;23.5.11.1}
There exist $B_i$ $(i=1,2)$ such that
the following holds on $X$ for any $t\geq 1$:
\begin{equation}
 |\gamma(\tau,t)|_{h_t}
+|(\delbar_E+\ad t\theta)\gamma(\tau,t)|_{h_t,g_X}
 \leq B_1\exp(-B_2t)\|\tau\|.
\end{equation}
\hfill\qed
\end{lem}

We set
\[
 \ttV(\tau,t)
:=\ttV'(\tau,t)
-(\delbar_E+\ad t\theta)\gamma(\tau,t).
\]
Then, $\ttV(\tau,t)$ is a harmonic $1$-form
satisfying
$[\ttV(\tau,t)]=\iota_t^{\ver}([\tau])$.

\subsubsection{Estimate of the pairings}

\begin{prop}
\label{prop;23.5.5.230}
There exist $B_i$ $(i=1,2)$ such that
\begin{equation}
 \label{eq;23.5.16.21}
\Bigl|
g_{H}(\iota^{\ver}_t(\tau),\iota^{\ver}_t(\tau))
-g_{\semiflat}(\iota^{\ver}_t(\tau),\iota^{\ver}_t(\tau))
\Bigr|
\leq
B_1\exp(-B_2t)\|\tau\|^2.
\end{equation}
\end{prop}
\pf
The estimate (\ref{eq;23.5.16.21}) is equivalent to the following:
\begin{equation}
\label{eq;23.5.5.200}
\left|
 2\sqrt{-1}
 \int_X
 \Bigl(
 h_{t}(\ttV(\tau,t)^{1,0},\ttV(\tau,t)^{1,0})
-h_{t}(\ttV(\tau,t)^{0,1},\ttV(\tau,t)^{0,1})
 \Bigr)
 -2(-\sqrt{-1})
 \int_{\Sigma_{\theta}}
 \tau\overline{\tau}
\right|
\leq B_1\exp(-B_2t)\|\tau\|^2.
 \end{equation}
Because $\del\tau=0$, we have
\[
  -2\sqrt{-1}
 \int_{\Sigma_{\theta}}
 \tau\overline{\tau}
 =-2\sqrt{-1}
 \int_{\Sigma_{\theta}}
 \tau'\overline{\tau}.
\]
Because $(\delbar_E+\ad t\theta)^{\ast}_{h_t,g_X}\ttV(\tau,t)=0$,
we have
\begin{multline}
 2\sqrt{-1}
 \int_X
 \Bigl(
 h_{t}(\ttV(\tau,t)^{1,0},\ttV(\tau,t)^{1,0})
-h_{t}(\ttV(\tau,t)^{0,1},\ttV(\tau,t)^{0,1})
 \Bigr)
=-2\sqrt{-1}
 \int_Xh_t(F_{\tau'},\ttV(\tau,t)^{0,1})
 \\
=-2\sqrt{-1}
 \int_X
 \tr\bigl(
  F_{\tau'}\cdot
  (\ttV(\tau,t)^{0,1})^{\dagger}_{h_t}
 \bigr).
\end{multline}
Then, we obtain {\rm(\ref{eq;23.5.5.200})}
from 
Lemma \ref{lem;23.5.11.2} and Lemma \ref{lem;23.5.11.1}.
\hfill\qed
 
\begin{prop}
\label{prop;23.5.5.231}
There exist $B_1,B_2>0$
such that 
the following holds
for any $\nu\in H^0(\Sigma_{\theta},K_{\Sigma_{\theta}})$
and $\tau\in H^1(\Sigma_{\theta},\nbigo_{\Sigma_{\theta}})$:
\begin{equation}
 \label{eq;23.5.16.22}
\bigl|
g_H(\iota^{\ver}_t(\tau),\iota^{\hor}_t(\nu))
\bigr|
\leq B_1\exp(-B_2t)\|\nu\|\|\tau\|.
\end{equation}
\end{prop}
\pf
The estimate (\ref{eq;23.5.16.22}) is equivalent to the following:
\begin{equation}
\label{eq;23.5.5.220}
\left|
  2\sqrt{-1} 
 \int_X\Bigl(
 h_t\bigl(
 \ttV(\tau,t)^{1,0},\ttH(\nu,t)^{1,0}
 \bigr)
-h_t\bigl(
 \ttV(\tau,t)^{0,1},\ttH(\nu,t)^{0,1}
 \bigr)
  \Bigr)
\right|
\leq B_1\exp(-B_2t)\|\nu\|\|\tau\|.
\end{equation}
Because $(\delbar_E+\ad t\theta)^{\ast}_{h_t,g_X}\ttH(\nu,t)=0$,
we have
\begin{multline}
   2\sqrt{-1} 
 \int_X\Bigl(
 h_t\bigl(
 \ttV(\tau,t)^{1,0},\ttH(\nu,t)^{1,0}
 \bigr)
-h_t\bigl(
 \ttV(\tau,t)^{0,1},\ttH(\nu,t)^{0,1}
 \bigr)
 \Bigr)
 =
 2\sqrt{-1}
 \int_X
 h_t(F_{\tau'},\ttH(\nu,t)^{0,1})
 \\
 =2\sqrt{-1}
 \int_X
 \tr\bigl(
 F_{\tau'}\cdot
 (\ttH(\nu,t)^{0,1})^{\dagger}_{h_t}
 \bigr).
\end{multline}
We obtain  {\rm(\ref{eq;23.5.5.220})}
from Lemma \ref{lem;23.5.5.51} and Lemma \ref{lem;23.5.5.80}.
\hfill\qed

\end{document}

%% file: notation.tex
\newcommand{\nbiga}{\mathcal{A}}
\newcommand{\nbigb}{\mathcal{B}}
\newcommand{\nbigc}{\mathcal{C}}

\newcommand{\nbige}{\mathcal{E}}
\newcommand{\nbigf}{\mathcal{F}}

\newcommand{\nbigh}{\mathcal{H}}

\newcommand{\nbigl}{\mathcal{L}}
\newcommand{\nbigm}{\mathcal{M}}

\newcommand{\nbigo}{\mathcal{O}}
\newcommand{\nbigp}{\mathcal{P}}

\newcommand{\nbigr}{\mathcal{R}}

\newcommand{\nbigu}{\mathcal{U}}

\newcommand{\nbigx}{\mathcal{X}}

\newcommand{\proj}{\mathbb{P}}
\newcommand{\seisuu}{{\mathbb Z}}

\newcommand{\cnum}{{\mathbb C}}
\newcommand{\real}{{\mathbb R}}

\newcommand{\DD}{\mathbb{D}}

\newcommand{\gbiga}{\mathfrak A}

\newcommand{\gbigr}{\mathfrak R}

\newcommand{\gminil}{\mathfrak l}

\newcommand{\gminis}{\mathfrak s}

\newcommand{\gminiu}{\mathfrak u}

\newcommand{\vecsigma}{{\boldsymbol \sigma}}

\newcommand{\vecrho}{{\boldsymbol \rho}}

\newcommand{\vecv}{{\boldsymbol v}}
\newcommand{\vecu}{{\boldsymbol u}}

\newcommand{\vecs}{{\boldsymbol s}}

\newcommand{\llarr}{\longleftarrow}

\newcommand{\lrarr}{\longrightarrow}

\newcommand{\pf}{{\bf Proof}\hspace{.1in}}
\newcommand{\qed}{\mbox{\rule{1.2mm}{3mm}}}

\def\Hom{\mathop{\rm Hom}\nolimits}

\def\End{\mathop{\rm End}\nolimits}

\def\Cok{\mathop{\rm Cok}\nolimits}

\def\Re{\mathop{\rm Re}\nolimits}

\def\Gr{\mathop{\rm Gr}\nolimits}

\def\GL{\mathop{\rm GL}\nolimits}
\def\SL{\mathop{\rm SL}\nolimits}
\def\Tot{\mathop{\rm Tot}\nolimits}

\def\rank{\mathop{\rm rank}\nolimits}

\def\Ker{\mathop{\rm Ker}\nolimits}

\def\Gr{\mathop{\rm Gr}\nolimits}
\def\Sym{\mathop{\rm Sym}\nolimits}
\def\sym{\mathop{\rm sym}\nolimits}
\def\ad{\mathop{\rm ad}\nolimits}
\def\Res{\mathop{\rm Res}\nolimits}

\def\tr{\mathop{\rm tr}\nolimits}
\def\Tr{\mathop{\rm Tr}\nolimits}

\def\id{\mathop{\rm id}\nolimits}

\def\Supp{\mathop{\rm Supp}\nolimits}

\newcommand{\del}{\partial}
\newcommand{\delbar}{\overline{\del}}

\newcommand{\barz}{\overline{z}}
\newcommand{\zbar}{\barz}

\newcommand{\etabar}{\overline{\eta}}

\newcommand{\abar}{\overline{a}}

\newcommand{\DDlambda}{\DD^{\lambda}}

\newcommand{\lefttop}[1]{{}^{#1}\!}

\def\Def{\mathop{\rm Def}\nolimits}
\def\Harm{\mathop{\rm Harm}\nolimits}

\newcommand{\AAA}{{\boldsymbol A}}

\newcommand{\minisl}{\gminis\gminil}

\newcommand{\thetahat}{\widehat{\theta}}
\newcommand{\Ehat}{\widehat{E}}
\newcommand{\hhat}{\widehat{h}}

\newcommand{\rhotilde}{\widetilde{\rho}}

\newcommand{\Etilde}{\widetilde{E}}

\newcommand{\thetatilde}{\widetilde{\theta}}

\newcommand{\wbar}{\overline{w}}

\newcommand{\kappatilde}{\widetilde{\kappa}}

\newcommand{\htilde}{\widetilde{h}}

\newcommand{\vtilde}{\widetilde{v}}

\newcommand{\Ftilde}{\widetilde{F}}

\newcommand{\atilde}{\widetilde{a}}
\newcommand{\stilde}{\widetilde{s}}

\newcommand{\Psitilde}{\widetilde{\Psi}}
\newcommand{\nbiglhat}{\widehat{\nbigl}}

\newcommand{\Utilde}{\widetilde{U}}
\newcommand{\Dtilde}{\widetilde{D}}
\newcommand{\Xtilde}{\widetilde{X}}

\def\Proj{\mathop{\rm Proj}\nolimits}

\def\Crit{\mathop{\rm Crit}\nolimits}

\def\asym{\mathop{\rm asym}\nolimits}
\def\Pic{\mathop{\rm Pic}\nolimits}
\def\semiflat{\mathop{\rm sf}\nolimits}

\newcommand{\Ubar}{\overline{U}}

\newcommand{\tauhat}{\widehat{\tau}}

\newcommand{\Ctilde}{\widetilde{C}}

\newcommand{\betabar}{\overline{\beta}}

\newcommand{\taubar}{\overline{\tau}}

\newcommand{\kappahat}{\widehat{\kappa}}

\newcommand{\ttF}{{\tt F}}

\newcommand{\ttH}{{\tt H}}
\newcommand{\ttV}{{\tt V}}

\newcommand{\ttv}{{\tt v}}
\newcommand{\tth}{{\tt h}}

\newcommand{\vectau}{\boldsymbol{\tau}}

\newcommand{\nubar}{\overline{\nu}}
\newcommand{\betatilde}{\widetilde{\beta}}
\newcommand{\Lhat}{\widehat{L}}

\newcommand{\hor}{\tth}
\newcommand{\ver}{\ttv}

%% file: new_theorem.tex
\newtheorem{thm}{Theorem}[section]
\newtheorem{cor}[thm]{Corollary}

\newtheorem{rem}[thm]{Remark}
\newtheorem{lem}[thm]{Lemma}
\newtheorem{prop}[thm]{Proposition}
\newtheorem{df}[thm]{Definition}

\newtheorem{condition}[thm]{Condition}